\RequirePackage{etex}
\documentclass{amsart}
\usepackage{amsmath, amsthm, amssymb,  mathrsfs, epsfig}
\usepackage{dcpic, pictex, pinlabel}

\DeclareMathOperator{\vol}{vol}
\DeclareMathOperator{\tr}{tr}

\DeclareMathOperator{\SL}{\textrm{SL}}
\DeclareMathOperator{\GL}{\textrm{GL}}

\DeclareMathOperator{\Hom}{\textrm{Hom}}
\DeclareMathOperator{\stab}{stab}
\DeclareMathOperator{\Inv}{Inv}
\DeclareMathOperator{\image}{image}

\begin{document}

\newtheorem{theorem}{Theorem}[section]
\newtheorem{lemma}[theorem]{Lemma}
\newtheorem{lctheorem}[theorem]{theorem}
\newtheorem{proposition}[theorem]{Proposition}
\newtheorem{sublemma}[theorem]{Sublemma}
\newtheorem{corollary}[theorem]{Corollary}
\newtheorem{conjecture}[theorem]{Conjecture}
\newtheorem{question}[theorem]{Question}
\newtheorem{problem}[theorem]{Problem}
\newtheorem*{claim}{Claim}
\newtheorem*{criterion}{Criterion}
\newtheorem*{main_theorem}{Theorem A}
\newtheorem*{cauchy_theorem}{Topological Cauchy-Schwarz inequality}
\newtheorem*{lemma_schema}{Lemma Schema}

\theoremstyle{definition}
\newtheorem{definition}[theorem]{Definition}
\newtheorem{construction}[theorem]{Construction}
\newtheorem{notation}[theorem]{Notation}
\newtheorem{convention}[theorem]{Convention}
\newtheorem*{explanation}{Explanation}
\newtheorem*{warning}{Warning}

\theoremstyle{remark}
\newtheorem{remark}[theorem]{Remark}
\newtheorem{example}[theorem]{Example}

\numberwithin{equation}{subsection}

\newcommand{\marginal}[1]{\leavevmode\marginpar{\tiny\raggedright#1}} 

\newcommand{\innerprod}[2]{\langle {#1} , {#2} \rangle}

\newcommand{\torus}{{\rm{torus}}}
\newcommand{\area}{{\rm{area}}}
\newcommand{\length}{{\rm{length}}}
\newcommand{\id}{{\rm{Id}}}
\newcommand{\R}{{\mathbb R}}
\newcommand{\RP}{{\mathbb R}P}
\newcommand{\Z}{{\mathbb Z}}
\newcommand{\Q}{{\mathbb Q}}
\newcommand{\T}{{\mathbb T}}
\newcommand{\M}{{\mathcal M}}
\newcommand{\gn}{\textrm{gn}}   
\newcommand{\ess}{{\mathcal S}}
\newcommand{\dee}{{\mathcal D}}
\newcommand{\Mdot}{{\dot\M}}
\newcommand{\Pdot}{{\dot{{\mathcal P}}}}
\newcommand{\Owe}{{\mathcal O}} 
\renewcommand{\H}{{\mathbb H}}
\newcommand{\C}{{\mathbb C}}
\newcommand{\bbar}{\overline}
\newcommand{\hhat}{\widehat}
\newcommand{\til}{\widetilde}
\newcommand{\mmin}{\text{min}}
\newcommand{\inc}{\text{inc}}
\newcommand{\Ric}{\text{Ric}}
\newcommand{\inv}{\text{inv}}
\newcommand{\cusp}{\text{cusp}}
\newcommand{\Xb}{\stackrel{\leftrightarrow}{X}}

\def\hline{\bigskip\hrule\bigskip}  
\def\nn#1{{\it [#1]}}       
\def\cN{{\mathcal N}}
\def\cO{{\bf O}}    
\def\Oh{{\bf O}}    

\def\tn{\textnormal}

\title{Positivity of the universal pairing in $3$ dimensions}
\author{Danny Calegari}
\address{Department of Mathematics \\ Caltech \\
Pasadena CA, 91125}
\email{dannyc@its.caltech.edu}
\author{Michael H. Freedman}
\address{Microsoft Station Q \\ University of California \\
Santa Barbara, CA 93106}
\email{michaelf@microsoft.com}
\author{Kevin Walker}
\address{Microsoft Station Q \\ University of California \\
Santa Barbara, CA 93106}
\email{kevin@canyon23.net}

\date{6/4/2009, Version 1.06}

\begin{abstract}
Associated to a closed, oriented surface $S$ is the complex vector space with basis the
set of all compact, oriented $3$-manifolds which it bounds. Gluing along $S$ defines a
Hermitian pairing on this space with values in the complex vector space with basis all
closed, oriented $3$-manifolds. The main result in this paper is that this pairing
is {\em positive}, i.e.\ that the result of pairing a nonzero vector with itself is nonzero.
This has bearing on the question of what kinds of topological information can be extracted
in principle from unitary $(2+1)$-dimensional TQFTs.

The proof involves the construction of a suitable complexity function $c$ on all closed
$3$-manifolds, satisfying a gluing axiom which we call the {\em topological Cauchy-Schwarz
inequality}, namely that $c(AB) \le \max(c(AA),c(BB))$ for all $A,B$ which bound $S$,
with equality if and only if $A=B$.

The complexity function $c$ involves input from many aspects of $3$-manifold topology, and
in the process of establishing its key properties we obtain a number of results of
independent interest. For example, we show that when two finite volume hyperbolic
$3$-manifolds are glued along an incompressible acylindrical surface, the resulting
hyperbolic $3$-manifold has minimal volume {\em only} when the gluing can be done along
a totally geodesic surface; this generalizes a similar theorem for closed
hyperbolic $3$-manifolds due to Agol-Storm-Thurston.
\end{abstract}


\maketitle

\section{Introduction}

The earliest objects of study in the theory of manifolds were the fundamental
group, Poincar\'e duality, and by the 1930's, characteristic classes.
The main theme was classification: developing invariants to
distinguish one manifold from another and trying to construct
manifolds with prescribed invariants. In the mid 50's to early 60's,
Thom, Milnor, and Smale shifted the
emphasis to {\em operations} on manifolds such as cutting, gluing and especially
surgery, using powerful new structural tools such as Morse theory and cobordism.
The most spectacular successes of this
theory were confined to sufficiently high dimensions.
Especially in $3$ dimensions, arbitrary surgery or cutting and
pasting is too disruptive, and incompatible with even the coarsest
features of the classification theory such as the prime
decomposition and JSJ theorems. Three manifold topology developed
quite independently reaching a culmination, from the classification
perspective, with Thurston and Perelman.

From a completely different direction, low-dimensional topology has
been invigorated over the last two decades by ideas from physics,
especially quantum field theory. Classically, one studies {\em
fields} (e.g.\ functions, or sections of some bundle) and their
dynamics on Euclidean space, or on some smooth manifold. From
quantum physics one gets the key idea of {\em superposition} ---
that one should study complex linear combinations of fields. A
combination of these two ideas --- cobordism and superposition ---
unexpectedly resonates in low-dimensional topology, and gives rise
to a host of beautiful and subtle invariants, such as the Jones
polynomial, Reshetikhin-Turaev-Viro invariants, and Chern-Simons
partition functions.

In fact, these invariants are perhaps {\em too} subtle. After twenty
years of work, it is profoundly frustrating that we cannot say
precisely what these invariants measure or what they distinguish,
and despite some tantalizing hints (e.g.\ Kashaev's conjecture),
these invariants remain disconnected from the (hugely successful)
Thurston theory of $3$-manifolds. There are two key questions: what
information does $3$-dimensional quantum
topology distinguish in principle, particularly in the physically
motivated unitary case?  and how does it relate in detail to the
structure theory of $3$-manifolds revealed by the geometrization
program?

This paper addresses both of these questions simultaneously.
Firstly, we establish positivity for
the {\em universal manifold pairing} in $3$-dimensions, an abstract
``universal topological quantum field theory'' which is, by
construction, sensitive to all details of $3$-manifold topology.
Three dimensions is the critical case here: in two dimensions and
lower, such positivity is straightforward to establish; in four dimensions and
higher, positivity fails badly (see \cite{FKNSWW} and
\cite{Kreck_Teichner}). For instance, the partition function of
a unitary $(3+1)$-dimensional TQFT is equal on
$s$-cobordant $4$-manifolds regardless of their Donaldson
invariants (\cite{FKNSWW}). Secondly, the process of establishing positivity turns
out to involve the construction of a complexity function on closed
$3$-manifolds which involves input from every aspect of the
geometric theory of $3$-manifolds, and obeys a (highly nontrivial)
gluing axiom, which may be thought of as a kind of {\em topological
Cauchy-Schwarz inequality}.

\vskip 12pt

We now make this discussion more precise. We work exclusively in the
category of smooth, compact, oriented manifolds. Let $S$ denote a
closed $(n-1)$-dimensional manifold which is not assumed to be
connected. Let $\Mdot(S)$ denote the set of isomorphism classes of
$n$-manifolds marked by an identification of their boundary with
$S$. That is, if $A,B$ are compact $n$-manifolds with $\partial A =
S = \partial B$, then $A$ and $B$ represent the same element of
$\Mdot(S)$ if and only if there is a diffeomorphism from $A$ to $B$
extending the identity map on their boundaries. Inspired by the
superposition principle, we let $\M(S)$ denote the complex vector
space spanned by the set $\Mdot(S)$ (i.e.\ finite $\C$-linear
combinations of elements of $\Mdot(S)$). If $S$ is empty, we
abbreviate $\Mdot(\emptyset)$ and $\M(\emptyset)$ by $\Mdot$ and
$\M$ respectively.

Given $A,B \in \Mdot(S)$, let $\overline{B}$ denote $B$ with the
opposite orientation. Note that $\overline{B} \in \Mdot(\overline{S})$.
By abuse of notation, we denote by $AB$ the result of
gluing $A$ to $\overline{B}$ by the identity map on their boundaries.

Our central object of study is the complex sesquilinear pairing
$$\M(S) \times \M(S) \longrightarrow \M$$
defined by the formula
$$\sum a_i A_i \times \sum b_j B_j  \longrightarrow
\sum a_i \bbar{b}_j A_i B_j$$ This pairing is known as the {\em
universal pairing} associated to $S$, and is denoted $\langle \cdot,
\cdot \rangle_S$. A pairing on a vector space is {\em positive} if
$\langle v, v \rangle = 0$ if and only if $v=0$. With this
terminology, our main theorem is the following.

\begin{main_theorem}
For all closed, oriented surfaces $S$, the pairing
$\innerprod{\cdot}{\cdot}_S$ is positive.
\end{main_theorem}

In a necessarily non-positive pairing, a vector $v$ satisfying $v
\neq 0$ and $\langle v, v \rangle = 0$ is said to be {\em lightlike}.

The motivation for studying this abstract pairing comes from the
examples from quantum topology briefly alluded to above, most of
which are examples of {\em topological quantum field theories}, or
TQFT's for short. In Atiyah's formulation (\cite{Atiyah}) a TQFT is
a functor $Z$ from the category of $(n-1)$-manifolds and cobordisms
between them, to the category of $\C$ vector spaces and linear maps
between them, satisfying the monoidal axiom
$$Z(S_1 \coprod S_1) = Z(S_1) \otimes Z(S_2)$$
In particular, one has $Z(\emptyset) = \C$. If
$\emptyset \xrightarrow{A} S$ is a cobordism from $\emptyset$ to $S$ (i.e.\ an
element of $\Mdot(S)$) then by abuse of notation, we define $Z(A) \in Z(S)$
to be the image of $1 \in \C = Z(\emptyset)$ under the map
$Z(A):Z(\emptyset) \to Z(S)$ associated to the cobordism.

\begin{remark}
For some reason, it is common practice to use the letters $V$ and $Z$ to denote
the image of the functor on an $(n-1)$-manifold and an $n$-manifold respectively,
so that for instance $Z(M) \in V(S)$ for $M \in \Mdot(S)$. We typically do not use this
convention when discussing abstract TQFTs, but we do sometimes
when discussing specific TQFTs in order to be consistent with the wider literature.
\end{remark}

In many interesting TQFT's, the images $Z(A)$ span $Z(S)$. For
example, if $Z$ is the $SU(2)$ Chern-Simons partition function at
level $k$ where $k+2$ is prime, then Roberts \cite{Roberts} has
shown the images $Z(A)$ over all $A \in \Mdot(S)$ span $Z(S)$. In
this case, there is a natural pairing
$$Z(S) \times Z(\overline{S}) \to Z(\emptyset) = \C$$
defined on generators by composition of cobordisms:
$$\emptyset \xrightarrow{A} S \circ S \xrightarrow{B} \emptyset
= \emptyset \xrightarrow{AB} \emptyset$$
and extended by sesquilinearity to
$$\langle \cdot, \cdot \rangle_Z:\M(S) \times \M(S) \rightarrow \C$$

In physical quantum field theories, $Z(S)$ denotes the vector space
of quantum fields on a spacelike slice $S$ of spacetime.
The vector space $Z(S)$ is naturally a (positive definite) {\em Hilbert space}.
This motivates an additional axiom for
so-called {\em unitary} TQFT's, that $Z(S)$ (now usually finite
dimensional) should admit the natural structure of a Hilbert space,
and for every nonzero $A \in \M(S)$, the pairing defined above
should satisfy
$$\innerprod{Z(A)}{Z(\overline{A})}_Z = Z(A\overline{A}) = \|A\|^2 > 0$$

It follows that studying positivity of the universal pairing (in a
given dimension) is tantamount to studying what kind of topological
information in principle might be extracted from unitary TQFT's in
that dimension, since lightlike vectors must map to zero in any
unitary TQFT.  As a specific and important example, given any compact Lie
group $G$ and level $k>0$, the Reshetikhin-Turaev TQFT \cite{Reshetikhin_Turaev}, denoted by
$V_{G,k}$, and as reconstructed by \cite{BHMV} fits into the following diagram:
\[\begindc{0}[3]
    \obj(10,30)[A]{$\M(S) \times \M(S)$}
    \obj(40,30)[B]{$\M(S)$}
    \obj(10,10)[C]{$V_{G,k} \times V_{G,k}$}
    \obj(40,10)[D]{$\C$}
    \mor(4,30)(4,10){$\alpha$}
    \mor(15,30)(15,10){$\alpha$}
    \mor(40,30)(40,10){$Z_{G,k}$}
    \mor(18,30)(37,30){$\pi$}
    \mor(18,10)(37,10){$A$}
\enddc\]
where $\pi$ is the universal pairing, and $A$ is Atiyah's map for
$(G,k)$. Since the pairing $A$ makes $V_{G,k}$ into a Hilbert
space, if $v$ is lightlike then necessarily $\alpha(v)=0$.

\vskip 12pt

Given a sesquilinear pairing on a complex vector
space, positivity of the pairing is logically equivalent to the
Cauchy-Schwarz inequality. In our context,
positivity is established by a topological variant of this
inequality. How can positivity fail? Only if there is a nonzero vector
$v = \sum_i a_i A_i$ with the property that there are enough topological
``coincidences'', i.e.\ equalities of the form $A_iA_j = A_kA_l$ for
various $i,j,k,l$, so that the coefficients, when collected according
to diffeomorphism types, can cancel and the result of pairing $v$
with itself is trivial. In order to rule out such a possibility
we construct a complexity function, defined
on all closed, oriented $3$-manifolds, which is
maximized {\em only} on terms of the form $A_iA_i$. Since the
coefficients of such terms are of the form
$a_i\bar{a_i} = \|a_i\|^2 > 0$, there can be no cancellation, and the result of self-pairing
is never zero. Explicitly, the complexity function we define
satisfies the following gluing axiom:

\begin{cauchy_theorem}
There is an ordered set $\Owe$ and a complexity function $c$ defined
on the set $\Mdot$ of all closed oriented $3$-manifolds, and taking
values in $\Owe$, so that for all $A,B \in \Mdot(S)$,
$$c(AB) \le \max(c(AA),c(BB))$$
with equality {\em if and only if} $A=B$.
\end{cauchy_theorem}

\begin{remark}
It is natural to ask to what extent this analogy with the ordinary
Cauchy-Schwarz inequality can be sharpened.
\begin{question}
Is there an embedding $\Owe \to \R$ such that the resulting complexity
function $\overline{c}:\Mdot \to \R$ satisfies
$$\overline{c}(AB)^2 \le \overline{c}(AA)\overline{c}(BB)$$
for all $A,B \in \Mdot$, with equality if and only if $A = B$?
\end{question}
Note that for any finite subset $\Sigma$
of $\Mdot$, the values of $\overline{c}$ on $\Sigma$ can be chosen
to satisfy the inequality above.
The subtlety of this question concerns the order structure of $\Owe$,
particularly its accumulation points.
\end{remark}

\begin{figure}[ht]
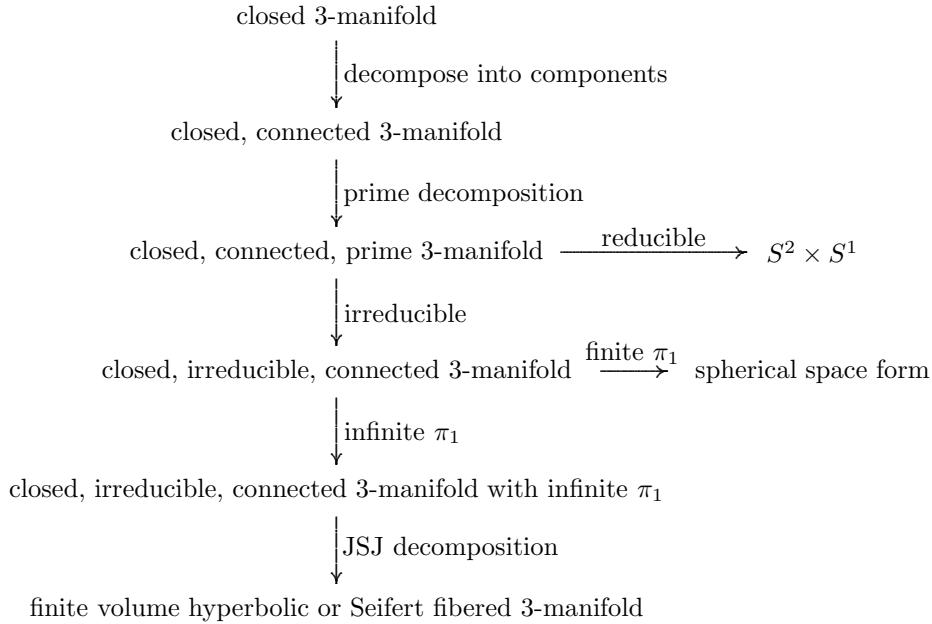

\[\begindc{0}[3]
    \obj(20,60){closed $3$-manifold}
    \obj(20,45){closed, connected $3$-manifold}
    \obj(20,30){closed, connected, prime $3$-manifold}
    \obj(20,15){closed, irreducible, connected $3$-manifold}
    \obj(20,0){closed, irreducible, connected $3$-manifold with infinite $\pi_1$}
    \obj(20,-15){finite volume hyperbolic or Seifert fibered $3$-manifold}

    \obj(80,30){$S^2 \times S^1$}
    \obj(80,15){spherical space form}

    \mor(20,60)(20,45){decompose into components}
    \mor(20,45)(20,30){prime decomposition}
    \mor(45,30)(75,30){reducible}
    \mor(20,30)(20,15){irreducible}
    \mor(50,15)(65,15){finite $\pi_1$}
    \mor(20,15)(20,0){infinite $\pi_1$}
    \mor(20,0)(20,-15){JSJ decomposition}
\enddc\]
\caption{repeated decompositions in the geometrization theorem}\label{hierarchy_flowchart}
\end{figure}

\begin{figure}[htpb]
\labellist
\small\hair 2pt
{\tiny
\pinlabel $\text{means ``proceed to next''}$ at 140 627
\pinlabel $\text{means ``make new complexity by recursing over something''}$ at 220 609
\pinlabel $\text{(e.g.\ connected components, primes, JSJ pieces, etc.}$ at 210 598
}
\pinlabel $c_0$ at 50 562
\pinlabel $\parallel$ at 50 544
\pinlabel $\#\text{ of path}$ at 50 526
\pinlabel $\text{components}$ at 50 510

\pinlabel $c_1$ at 160 562
\pinlabel $\parallel$ at 160 544
\pinlabel $\text{Dijkgraaf-Witten}$ at 158 526
\pinlabel $\text{partition functions}$ at 164 510
\pinlabel $\text{over all finite groups}$ at 170 494

\pinlabel $c_2$ at 270 562
\pinlabel $c_2'$ at 270 505
\pinlabel $\parallel$ at 270 487
\pinlabel $(r,s)$ at 270 472
\pinlabel $\#\text{irreducible primes}$ at 220 430
\pinlabel $\#S^1\times S^2\text{ factors}$ at 300 410

\pinlabel $c_3$ at 440 562
\pinlabel $\parallel$ at 440 544
\pinlabel $(c_p,$ at 360 526
\pinlabel $\sum c_p,c_\iota$ at 430 526
\pinlabel $\sum c_\iota)$ at 505 526
\pinlabel $(p_i,d_i)$ at 452 460
\pinlabel $\#\text{occurrences}$ at 420 410
\pinlabel $\text{of }i\text{th prime}$ at 420 395
\pinlabel $\text{are primes in}$ at 510 395
\pinlabel $\text{mirror pairs?}$ at 510 380

\pinlabel $c_S$ at 66 258
\pinlabel $\text{Seifert fibered}$ at 126 260
\pinlabel $\parallel$ at 66 242
\pinlabel $(m,-m',\sum -\chi^o(Q_i),\sum\text{genus}(Q_i),n,c_{CS})$ at 106 226
\pinlabel $(b,-\chi^o(Q),\#\text{sing},-|e|)$ at 110 62
{\tiny
\pinlabel $\#\text{independent tori}$ at 0 180
\pinlabel $\#\text{JSJ tori}$ at 50 140
\pinlabel $-\text{orbifold}$ at 100 185
\pinlabel $\text{Euler char.}$ at 100 175
\pinlabel $\text{base genus}$ at 135 140
\pinlabel $\#\text{SF JSJ}$ at 170 173
\pinlabel $\text{pieces}$ at 170 163

\pinlabel $\#\text{boundary}$ at 35 30
\pinlabel $\text{components}$ at 35 20
\pinlabel $-\text{orbifold Euler char.}$ at 80 -10
\pinlabel $\text{singular fiber}$ at 130 30
\pinlabel $\text{data tuple}$ at 130 20
\pinlabel $\text{rational Euler class}$ at 185 7
}

\pinlabel $c_h$ at 268 254
\pinlabel $\text{hyperbolic}$ at 315 256
\pinlabel $c_{ch}$ at 302 198
\pinlabel $\text{connected hyperbolic}$ at 380 200
\pinlabel $\parallel$ at 302 180
\pinlabel $(-\text{volume},-\text{real length spectrum})$ at 315 162

\pinlabel $c_a$ at 482 254
\pinlabel $\text{assembly}$ at 522 256
\pinlabel $(c_r,c_t)$ at 390 100
\pinlabel $\text{reflection}$ at 330 70
\pinlabel $\text{complexity}$ at 330 55
\pinlabel $\text{tensor contraction}$ at 440 70
\pinlabel $\text{complexity}$ at 440 55
\endlabellist
\centering
\vskip 5pt
\includegraphics[height=6in]{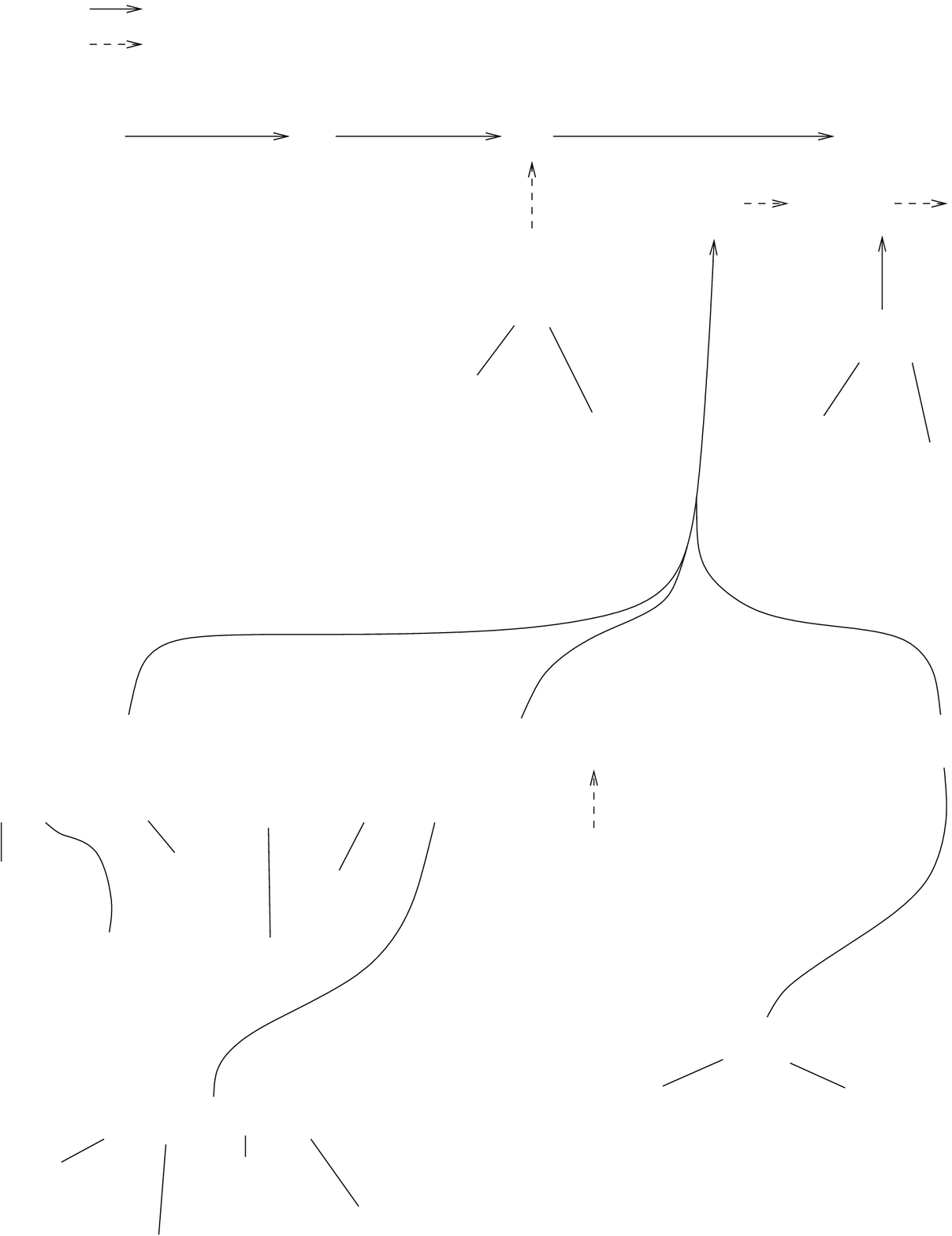}
\vskip 5pt
\caption{This flow chart may assist in navigating the proof.}
\label{complexity_flowchart}
\end{figure}

Although our main theorem has a very clean statement, it obscures an
Augean complexity in the definition of the function $c$ and the
ordered set $\Owe$. These objects involve highly nontrivial input
from many aspects of geometric $3$-manifold topology. This
complexity mirrors the complexity inherent in the geometrization
theorem for $3$-manifolds. The geometric classification theorem for
$3$-manifolds has a {\em hierarchical} structure. A manifold is
subject to a repeated sequence of decompositions, each of which only
makes sense after the manifold has been filtered through the
previous decompositions (see Figure~\ref{hierarchy_flowchart}). The
complexity function $c$ has terms which are sensitive to all
different levels of the decomposition. Only if two manifolds are
identical with respect to the first few levels of the decomposition
does $c$ discriminate between them based on features which make
sense at later levels. A flowchart
(Figure~\ref{complexity_flowchart}) summarizing the main components
of the complexity function $c$ is included to assist the reader in
understanding the logic of the argument.

Precisely because the function $c$ involves input from many
different aspects of $3$-manifold topology, several steps in the
proof of the Topological Cauchy-Schwarz inequality are of
independent interest to $3$-manifold topologists. These include:
\begin{itemize}
\item The use of Dijkgraaf-Witten finite group untwisted TQFTs to
determine handlebody and compression body splittings.
\item The construction of a $\Z[\frac 1 2]$-valued invariant of a non-closed
Seifert fibered piece $X$ in the canonical decomposition of a prime
$3$-manifold; this invariant,
a kind of relative Euler class, measures the degree to which the fibers
of $X$ along $\partial X$ are twisted as they are glued to
neighboring JSJ pieces.
\item The proof of a volume inequality for cusped hyperbolic manifolds
obtained by gluing manifolds with incompressible boundary, which
generalizes a similar inequality for closed manifolds obtained
recently by Agol-Storm-Thurston \cite{AST}.
\end{itemize}
These and other broader implications of the work are discussed
briefly as they arise, and also in an appendix.

It should be remarked that a surprisingly large fraction of the proof involves assembling the
individual complexities of hyperbolic and Seifert fibered pieces
together. In some ways the study of nontrivial JSJ decompositions is
an underestimated area of $3$-manifold topology, with many
practitioners preferring to work in either the ``pure'' hyperbolic
or Seifert fibered (or graph manifold) world. However, there are
important historical cases where the synthesis of these two worlds,
Seifert fibered and hyperbolic, is a crucial and delicate issue,
essential to obtaining the mathematical result. Especially
noteworthy cases are the resolution of the Smith conjecture, and the
cyclic surgery theorem. In both cases (as well as many others),
Cameron Gordon played the critical role both in terms of mathematics
and leadership, in this synthesis. Architecturally therefore, this
paper should be seen as part of a tradition in $3$-manifold topology
of which Cameron is perhaps the most preeminent champion, and it is
our pleasure to acknowledge his multi-generational inspiration.

Note that the assembly
of a $3$-manifold from its JSJ pieces amounts to a kind of
``decorated'' graph theory which encompasses many subtle points. In
\S\ref{acgcssect} we obtain new positivity results for certain
natural ``graph TQFTs'' as a warm-up to the $3$-manifold case.

\section{Universal pairings}

\subsection{Definition}

\begin{notation}
Throughout the paper, $S$ will usually denote a closed oriented
surface, possibly disconnected. $A,B,A_i,B_i$ and so on will denote
$3$-manifolds, usually compact with boundary equal to $S$, and
$M,N,M_i,N_i$ and so on will usually denote closed $3$-manifolds.
Sometimes these conventions are relativized, so that $S$ or $M$
might denote a complete manifold with a cusp, and $A$ or $B$ might
denote a complete manifold whose boundary is a (possibly cusped)
surface $S$.
\end{notation}

In what follows, let $S$ be a closed, oriented surface. We
explicitly allow $S$ to be disconnected. Let $\Mdot(S)$ denote the
set of smooth pairs $(A,S)$ where $A$ is a compact oriented
$3$-manifold and
$\partial A = S$, up to the equivalence relation $(A,S) \sim (B,S)$
if there is a diffeomorphism $h:A \to B$ such that $hi = j$ where
$i:S \to A$ and $j:S \to B$ are the respective inclusions. That is,
the following diagram commutes:

\[\begindc{0}[30]
    \obj(3,3){$A$}
    \obj(3,1){$B$}
    \obj(1,2){$S$}
    \mor{$S$}{$A$}{$i$}[\atright,\injectionarrow]
    \mor{$S$}{$B$}{$j$}[\atright,\injectionarrow]
    \mor{$A$}{$B$}{$h$}
\enddc\]

Let $\M(S)$ denote the $\C$ vector space spanned by the set
$\Mdot(S)$. If $S$ is empty, abbreviate $\Mdot:=\Mdot(\emptyset)$
and $\M:=\M(\emptyset)$.

\begin{definition}\label{universal_pairing_definition}
The {\em universal Hermitian pairing} is the map
$$\innerprod{\cdot}{\cdot}_S:\M(S) \times \M(S) \to \M$$
defined by the formula
$$\Big\langle \sum_i a_i A_i, \sum_j b_j B_j\Big\rangle_S = \sum_{i,j} a_i\overline{b}_j A_iB_j$$
where $A_iB_j$ denotes the closed $3$-manifold obtained by gluing
$A_i$ to $\overline{B}_j$ (i.e.\ the manifold $B_j$ with the opposite
orientation) along $S$, using the {\em canonical} inclusion of $S$
into either factor.
\end{definition}

\begin{definition}
Let $V,W$ be $\C$-vector spaces. A pairing
$\innerprod{\cdot}{\cdot}$ on $V$ with values in $W$ is {\em
positive} if $\innerprod{v}{v} =0$ if and only if $v=0$.
\end{definition}

The main theorem in this paper is

\begin{main_theorem}
For all closed, oriented surfaces $S$, the pairing
$\innerprod{\cdot}{\cdot}_S$ is positive.
\end{main_theorem}

The proof involves many ingredients, and will occupy most of the
remainder of the paper.

\subsection{Complexity and diagonal dominance}

The key to the proof of Theorem~A is the construction of a suitable
complexity function on compact, closed, oriented (not necessarily
connected) $3$-manifolds.

\begin{definition}
Let $\Owe$ be an ordered set. A function $c:\Mdot \to \Owe$ is {\em
diagonally dominant} for $S$ if for all {\em distinct} pairs $A,B \in
\Mdot(S)$ there is a strict inequality
$$c(AB) < \max(c(AA),c(BB))$$
\end{definition}

The following Lemma shows how Theorem~A follows from diagonal dominance.

\begin{lemma}[Diagonal dominance implies positive]\label{diagonal_dominant_implies_positive_lemma}
Suppose there is a function $c$ which is diagonally dominant for $S$. Then $\innerprod{\cdot}{\cdot}_S$
is positive.
\end{lemma}
\begin{proof}
Let $v = \sum_i a_i A_i$ 
with $a_i\ne 0$ and set $w = \innerprod{v}{v}_S$. The vector $w$ is a linear combination
of manifolds of the form $M_iM_j$. By diagonal dominance, the maximum of $c$ on the
set $A_iA_j$ is realized only on manifolds of the form $A_iA_i$.
Let $i$ be such that $c(A_iA_i)$ is maximum. Let $J$ be the set of indices $j$ for which
$A_jA_j = A_iA_i$, and observe that $A_jA_k \ne A_iA_i$ for $j \ne k$. It follows that
the coefficient of $w$ on the manifold $A_iA_i \in \Mdot$ is
$\sum_{j \in J} |a_j|^2 > 0$. Hence $w \ne 0$.
\end{proof}

Most of the remainder of the paper will be concerned therefore with defining a
suitable complexity function, and proving that it satisfies diagonal dominance.
In fact, the function $c$ we define does not depend on $S$.

The statement that such a diagonally dominant function exists is paraphrased
in the introduction as the {\em topological Cauchy-Schwarz inequality}, and we
will sometimes use such a term interchangeably with diagonal dominance.

\subsection{Overview}

The complexity function $c$ is constructed from a tuple of more specialized complexities
$c_0,c_1,c_2$ and $c_3$, each of which treats some narrower aspect of $3$-manifold topology.
Roughly, $c_0$ addresses connectivity, $c_1$ will tell us about the kernel of
the map on fundamental groups induced by the inclusion of $S$ into a bounding $3$-manifold,
$c_2$ is concerned with essential $2$-spheres,
while $c_3$ addresses the nature and multiplicities of the prime factors (as a mnemonic,
think of $c_i$ as addressing some aspect of  $\pi_i$ for each of $0\le i <3$).

Write $c^0 := c_0$ and $c^1 := c_0 \times c_1$ for the
lexicographic pair which first maximizes $c_0$ then $c_1$.
Similarly, $c^2 := c_0 \times c_1 \times c_2 = c^1 \times c_2$, and
finally, $c := c^3 := c_0 \times c_1 \times c_2 \times c_3 = c^2
\times c_3$.

\begin{warning}
The functions $c_1$, $c_2$, and $c_3$ are themselves
complicated, and are built up from many simpler complexity
functions.
\end{warning}

For each $0\le k\le 3$ and each complexity $c^k$ we will prove a lemma of the form

\begin{lemma_schema}
For all pairs $A,B \in \Mdot(S)$, either $$c^k(AB) < \max(c^k(AA),c^k(BB))$$
or $c^k(AB) = c^k(AA) = c^k(BB)$, and $(A,S)$ and $(B,S)$ are indistinguishable
with respect to a certain criterion $C_k$.
\end{lemma_schema}

The criteria $C_k$ are associated to the complexities $c^k$ in a certain way, and
chosen so that if $(A,S)$ and $(B,S)$ are indistinguishable with respect to criteria $C_k$ for
all $0 \le k \le 3$ then they represent the same element of $\Mdot(S)$. Actually, this
fact is absorbed into the statement of the $c^3$-Lemma Schema, which is
Theorem~\ref{c3_theorem}.

This fact together with Lemma Schema for $0 \le k \le 3$
together imply the existence of a diagonally dominant function (i.e.\ the topological
Cauchy-Schwarz inequality), and therefore
together with Lemma~\ref{diagonal_dominant_implies_positive_lemma}, they prove
Theorem~A.

The complexity $c_3$ is itself a lexicographic tuple $c_3 = c_S \times c_h \times c_a$
where the subscripts $S,h,a$ stand for {\em Seifert fibered}, {\em hyperbolic} and
{\em assembly} complexities respectively. The Lemma Schema for $c_3$ involves in its turn
a version of the Lemma Schema for $c_S,c_h,c_a$. Moreover each of these individual
internal complexity functions have their own internal structure involving simpler terms \ldots

\vskip 12pt

Here is the plan for the next few sections. In \S\ref{c_012_section} we state and prove the
Lemma Schema for $c_0,c_1,c_2$. The Seifert fibered complexity $c_S$ is treated
in \S\ref{c_S_section}, the hyperbolic complexity $c_h$ in \S\ref{c_h_section}, and
the assembly complexity (which encodes the gluing up of the Seifert fibered
and hyperbolic pieces) in \S\ref{c_a_section}.

\section{Prime decomposition}\label{c_012_section}

In this section we state and prove the Lemma Schema for the complexities
$c_0,c_1,c_2$ and the respective criteria $C_0,C_1,C_2$. Then we show how these
Lemma Schemas, together with a
complexity function $c_3$ on irreducible $3$-manifolds satisfying suitable
properties imply Lemma~\ref{diagonal_dominant_implies_positive_lemma} and therefore
Theorem~A, deferring the precise definition and proof of properties of $c_3$ to future
sections.

\subsection{The complexity $c_0$}

The definition of $c_0$ is straightforward.

\begin{definition}
Let $c_0: \Mdot \rightarrow \Z$ assign to a
closed 3-manifold $M$ the number of connected components of $M$ (an
empty manifold has zero components).
\end{definition}

\begin{definition}\label{connectivity_partition_definition}
Given $(A,S) \in \Mdot(S)$ the {\em connectivity partition} associated
to $A$ is the functor $\pi_0$ applied to the pair $(A,S)$.
\end{definition}

In other words, the connectivity partition is the pair of sets $\pi_0(S),\pi_0(A)$ and the
map $\pi_0(S) \to \pi_0(A)$ induced by inclusion $S \to A$. The data is equivalent
to the information about which components of $S$ are included into the same
component of $A$, and how many closed components $A$ has disjoint from $S$.

Criterion $C_0$ tries to distinguish elements of $\Mdot(S)$ using their connectivity partitions.
Two elements $(A,S),(B,S) \in \Mdot(S)$ are indistinguishable with respect to $C_0$
if and only if they induce the same connectivity partition.

\begin{lemma}[$c^0 = c_0$-Lemma Schema]\label{c^0_Lemma_Schema}
Either
$$c_0(AB) < \max(c_0(AA),c_0(BB))$$ or
$c_0(AB) = c_0(AA) = c_0(BB)$
and both $(A,S)$ and $(B,S)$ induce the same connectivity partition.
\end{lemma}
\begin{proof}
The key observation is that $c_0$ is maximized only when at most
one component from each of $A$ and $B$ contributes to a
component of $AB$.
\end{proof}

\subsection{The complexity $c_1$}

The complexity $c_1$ is less trivial.

To define $c_1$, begin with an ordered list ${F_1, F_2, F_3, \cdots}$
of all finite groups up to isomorphism.
According to \cite{Dijkgraaf_Witten}
there is a unique untwisted TQFT over $\mathbb{R}$ associated to
each $F_k$.

This TQFT associates to every closed surface $S$ a finite dimensional
$\R$ vector space $V_k(S)$ with an inner product $\innerprod{\cdot}{\cdot}_k$,
and to every element $(A,S) \in \Mdot(S)$ a vector
$Z_k(A,S) \in V_k(S)$ called the (relative) {\em partition function} of the pair.
When $S=\emptyset$ one has $V_k(\emptyset)=\R$, and the real number
$Z_k(M) \in \R$ is the {\em partition function} of $M$.

The {\em gluing axiom} for a TQFT states that for any two elements
$(A,S),(B,S) \in \Mdot(S)$ there is an equality
\begin{equation}\label{Z_k_equation1}
Z_k(AB) = \innerprod{Z_k(A, S)}{Z_k(B, S)}_k
\end{equation}

The crucial additional property of $\innerprod{\cdot}{\cdot}_k$ is
that it is positive definite (i.e.\ it has Euclidean signature).
For a more detailed discussion of TQFT's and a more precise definition
of $Z_k$ and $V_k$, see Appendix~\ref{TQFT_appendix}.

\begin{definition}
\label{c_1} $c_1(M) = (Z_1(M), Z_2(M), Z_3(M), \cdots)$ as a
lexicographic tuple of real numbers.
\end{definition}

\begin{definition}
Given $(A,S) \in \Mdot(S)$ the {\em compression set} is the set of isotopy classes
of embedded loops in $S$ which bound embedded disks in $A$.
\end{definition}

Criterion $C_1$ tries to distinguish elements of $\Mdot(S)$ using their compression sets.
Two elements $(A,S),(B,S) \in \Mdot(S)$ are indistinguishable with respect to $C_1$
if and only if they have the same compression sets.

\begin{remark}\label{Dehn_lemma_remark}
Dehn's Lemma and the loop theorem 
(see for example \cite{Hempel_book}, Chapter 4) imply that an essential simple
closed curve $\gamma \subset S$ is in the compression set if and only if its image in
$\pi_1(A)$ is trivial, under the homomorphism $\pi_1(S) \to \pi_1(A)$ induced by
inclusion $S \to A$. From this and from Nielsen's theorem
on presentations of free groups, one sees that
the data of a compression set is equivalent to the functor
$\pi_1$ (keeping track of basepoints, since $S$ might be disconnected)
applied to the pair $(A,S)$.
\end{remark}

\begin{remark}
As far as we know, it is possible that $c_1(A B) = c_1(A A) = c_1(B B)$
already implies the existence of 
an isomorphism $h: \pi_1(A) \to \pi_1(B)$
making the following
diagram commute:
\[\begindc{0}[30]
    \obj(3,3){$\pi_1(A)$}
    \obj(3,1){$\pi_1(B)$}
    \obj(1,2){$\pi_1(S)$}
    \mor{$\pi_1(S)$}{$\pi_1(A)$}{}
    \mor{$\pi_1(S)$}{$\pi_1(B)$}{}
    \mor{$\pi_1(A)$}{$\pi_1(B)$}{$h$}
\enddc\]
but we are not able to prove or disprove this.
Essentially, the question is whether the finite group
TQFTs, taken together, determine the induced map on $\pi_1$.
There seems to be no known technology that would address this question.
However, if it were known to be true, the remainder of the proof could be shortened considerably.
\end{remark}

\begin{lemma}[$c^1$-Lemma Schema]\label{c^1_Lemma_Schema}
Either
$$c^1(A B) < \max(c^1(A A), c^1(B B))$$
or $c^1(A B) = c^1(A A) = c^1(B B)$ and both
$(A,S)$ and $(B,S)$ induce the same compression set.
\end{lemma}

\begin{proof}
Since the inner product on $V_k(S)$ is Euclidean, the (ordinary) Cauchy-Schwarz
inequality implies that either
\begin{equation}\label{Z_k_equation2}
\innerprod{Z_k(A)}{Z_k(B)}_k < \max(\innerprod{Z_k(A)}{Z_k(A)}_k,\innerprod{Z_k(B)}{Z_k(B)}_k)
\end{equation}
or $Z_k(A)=Z_k(B)$, in which case all three inner products are equal.

Suppose there is an essential simple closed curve $\gamma \subset S$ which bounds
a disk in $A$ but not in $B$. By Dehn's Lemma (see Remark~\ref{Dehn_lemma_remark}) the loop
$\gamma$ is homotopically essential in $B$.
Since we may assume without loss of generality that
$c^0(A B) = c^0(AA) = c^0(BB)$, the inclusions $S \to A, S \to B$
induce isomorphic maps on $\pi_0$, and therefore we may restrict
attention to the components $A^0,B^0,S^0$ containing $\gamma$.
Choosing a basepoint somewhere on $\gamma$ lets us define the image of $\gamma$
in $\pi_1$ of these three spaces. By abuse of notation, we denote the image in
$\pi_1(S^0)$ by $\gamma$, and the image in $\pi_1(A^0),\pi_1(B^0)$ by $i_*(\gamma),j_*(\gamma)$
respectively.

\begin{claim}
There is a finite group $F_k$ and a homomorphism $\theta:\pi_1(B^0) \to F_k$ for which
$\theta \circ j_*(\gamma) = x \ne 1$.
\end{claim}
\begin{proof}
The fundamental group of $B^0$ is a free product of factors, corresponding to
the prime decomposition of $B^0$. The factor containing the image of $\gamma$ also
contains the image of $\pi_1(S^0)$. Since $\gamma$ is nontrivial in $S^0$, the genus of
$S$ is positive, and the rank of $H_1$ of this factor is positive. In particular,
the factor containing the image of $\gamma$ is $\pi_1$ of a Haken manifold (see \cite{Hempel_book}
Chapter 13 for a discussion of Haken manifolds).
A theorem of Hempel \cite{Hempel_residual} states that all
Haken $3$-manifolds have residually finite fundamental groups, and the claim follows.
\end{proof}

\begin{figure}[hpbt]
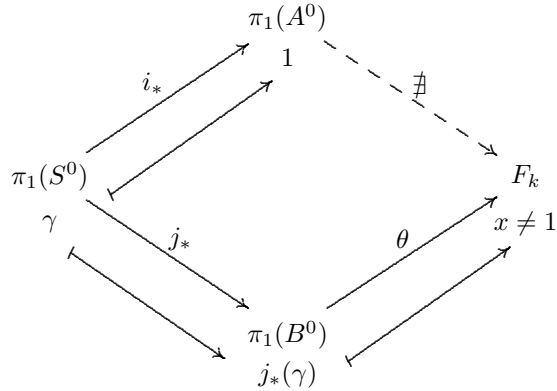

\[\begindc{0}[3]
    \obj(10,40)[S]{$\pi_1(S^0)$}
    \obj(10,34)[g]{$\gamma$}
    \obj(40,60)[Mi]{$\pi_1(A^0)$}
    \obj(40,55)[1]{1}
    \obj(40,20)[Mj]{$\pi_1(B^0)$}
    \obj(40,15)[incj]{$j_*(\gamma)$}
    \obj(70,40)[Fk]{$F_k$}
    \obj(70,34)[x]{$x \neq 1$}
    \mor{S}{Mi}{$i_*$}
    \mor{S}{Mj}{$j_*$}
    \mor(10,32)(35,15){}[\atright, \aplicationarrow]
    \mor(45,15)(71,33){}[\atright, \aplicationarrow]
    \mor{Mj}{Fk}{$\theta$}
    \mor{Mi}{Fk}{$\nexists$}[\atleft, \dasharrow]
    \mor(15,36)(41,54){}[\atright, \aplicationarrow]
\enddc\]
\caption{Since $\theta \circ j_*(\gamma) \neq 1$, the dashed arrow
cannot be filled in to a commutative diagram.}\label{unfillable_diagram_figure}
\end{figure}

Geometrically, the homomorphism $\theta$ defines a principal $F_k$ bundle $\beta'$ over
$S^0$ (by the Borel construction) which can be extended as a product to a
bundle $\beta$ over all $S$.
The bundle $\beta$ extends over $B$ but does not extend over
$A$ (see Figure~\ref{unfillable_diagram_figure}).
There is a natural basis for $Z_k(S)$ consisting of principal $F_k$ bundles of $S$ up to
isomorphism, such that the components of $Z_k(A) \in Z_k(S)$ count the number of extensions
(divided by the order of the fiber symmetry) for each principal $F_k$ bundle over $S$
(see Appendix~\ref{TQFT_appendix}). The component of $Z_k(B)$ in the direction $\beta$
is nonzero, whereas the component of $Z_k(A)$ in the direction $\beta$ is zero.
In particular,
\begin{equation}\label{Z_k_equation3}
Z_k(A) \neq Z_k(B) \in V_k(S)
\end{equation}

By equations~\ref{Z_k_equation1}, \ref{Z_k_equation2}, and \ref{Z_k_equation3} it follows that there
is an inequality
\begin{equation}\label{Z_k_equation4}
Z_k(A B) < \max(Z_k(A A),Z_k(B B))
\end{equation}
and the lemma is proved.
\end{proof}

\subsection{The complexity $c_2$}\label{c_2_section}

The complexity $c_2$ is a function of the decomposition of a closed $3$-manifold
into its irreducible factors. Similarly, two manifolds $(A,S)$ and $(B,S)$ in $\Mdot(S)$
are indistinguishable with respect to criterion $C_2$ if (roughly), they have
``identical'' sphere and disk decompositions into irreducible pieces.

There are various ambiguities in the way in which a $3$-manifold decomposes
into its prime factors, one of which arises from the $S^1 \times S^2$ factors. Such factors
are associated to non-separating spheres, but the way in which they sit in the
manifold is not unique. Roughly speaking, one $S^1 \times S^2$ factor can be
``slid'' over another; this operation is called a {\em handle slide}x.
A similar ambiguity holds for decompositions for manifolds with boundary.

If $M$ is a closed $3$-manifold, let $r(M)$ denote the number of irreducible
summands in the prime decomposition of $M$, and
let $s(M)$ denote the number of $S^1\times S^2$ summands.
The three-sphere is considered a unit, not a prime, so $r(S^3) = s(S^3)=0$.

Define $c_2'(M)$ to be the lexicographic pair $(r(M), s(M))$, and let
$M_1, M_2, \cdots$ denote the connected components of $M$, arranged in (any) order of
non-increasing complexity $c_2'(M_i)$.

\begin{definition}
$c_2(M) = (c_2'(M), c_2'(M_1), c_2'(M_2), \cdots)$ as a lexicographic tuple of
lexicographic pairs.
\end{definition}

\begin{convention}
Here and throughout, we make the convention that when we compare lexicographic lists of different
lengths, the shorter list is padded by $-\infty$ characters, where $-\infty$ is strictly
less, in the relevant partial ordering, than any other term. So if two lists agree on
their common lengths, the {\em longer} list has the {\em greater} lexicographic complexity.
\end{convention}

In order to understand how $c_2$ behaves under gluing, we adopt a more graphical way
of computing $r(M)$ and $s(M)$.

\begin{definition}[Sum graph]\label{sum_graph_definition}
Let $M$ be a closed, orientable $3$-manifold, and let
$\ess$ be a pairwise disjoint collection of embedded $2$-spheres in $M$.
The {\em sum graph} $G(\ess)$ is a finite decorated graph, defined as follows.
\begin{itemize}
\item{There is one edge for each $2$-sphere of $\ess$.}
\item{There is one vertex for each connected component of $M \setminus \ess$.}
\item{Each edge is attached to the two (not necessarily distinct)
components of $M \setminus \ess$ corresponding to
the two sides of the associated sphere.}
\item{For each vertex associated to a component $N$ of $M \setminus \ess$, 
label the vertex by $\overline{N}$, the result of capping off the 2-sphere
boundary components of $N$.
(The bar notation has a different meaning elsewhere in this paper, but the intended
meaning should always be clear from context.)}
\end{itemize}
A collection $\ess$ is {\em adequate} if all vertex labels are irreducible.
The existence of an adequate collection is due to Kneser.
\end{definition}

We define a {\em thin} vertex to be one whose label $\overline{N}$ is $S^3$;
other vertices are {\em fat}.
A {\em thin subgraph} consists of some (any) union of edges and thin vertices.

\begin{remark}
Note that adequate collections $\ess$ are allowed to have ``superfluous'' spheres --- i.e.
spheres which bound $3$-balls, or parallel collections of spheres.
\end{remark}

A given manifold $M$ admit many different adequate collections $\ess$ of spheres.
The following theorem, due to Laudenbach, describes the equivalence relation that this
ambiguity induces on sum graphs.

\begin{theorem}[Laudenbach \cite{Laudenbach}]\label{Laudenbach_theorem}
If $\ess,\ess'$ are adequate collections of spheres for $M$ then $G(\ess)$ and $G(\ess')$
are related by a sequence of edge slides and replacements of a thin subgraph with a homotopy
equivalent thin subgraph. Conversely, any sequence of such moves on a sum graph $G(\ess)$
produces a sum graph $G(\ess')$ where $\ess'$ is adequate if and only if $\ess$ is.
\end{theorem}

\begin{figure}[htpb]
\labellist
\small\hair 2pt
\pinlabel $\equiv$ at 228 35
\endlabellist
\centering
\includegraphics[scale=0.75]{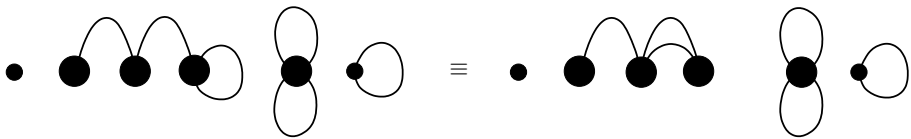}
\caption{Two different but equivalent sum graphs for a manifold of four components.  The
first component is $S^3$ and the other three components have 1, 2,
and 1 factors of $S^1 \times S^2$ respectively.} \label{sum_graph_figure}
\end{figure}

Figure~\ref{sum_graph_figure} shows two different but equivalent
sum graphs associated to the same manifold.

The numbers $r(M)$ and $s(M)$ can be recovered from any sum graph for $M$.

\begin{lemma}\label{sum_graph_properties}
Let $\ess$ be adequate for $M$. The number of fat vertices of $G(\ess)$ is equal to $r(M)$.
Moreover, $s(M) = \dim(H_1(G(\ess);\R))$.
\end{lemma}
\begin{proof}
This follows immediately from the definitions.
\end{proof}

To describe criterion $C_2$, we must relativize the definition of sum graphs
to compact $3$-manifolds with boundary.

\begin{definition}[Relative sum graph]\label{relative_sum_graph_definition}
Let $A$ be a compact, orientable $3$-manifold with boundary, and let $\dee$ be a
pairwise disjoint collection of embedded $2$-spheres and properly embedded disks in $A$.
The {\em relative sum graph} $G(\dee)$ is a finite decorated graph, defined as follows.
\begin{itemize}
\item{There is one edge for each element of $\dee$; an edge associated to a sphere is
{\em regular}; otherwise it is a {\em half edge}.}
\item{There is one vertex for each connected component of $A \setminus \dee$.}
\item{Edges join the (not necessarily distinct) components corresponding to their two sides
in $A$.}
\item{A vertex corresponding to a component disjoint from $\partial A$ is {\em regular};
otherwise it is a {\em half vertex}.}
\item{A regular vertex is thin or fat as before depending on whether its label is
$S^3$ or not.}
\item{A half vertex associated to a component $N$ determines a
label $\overline{N}$ obtained by capping off sphere components
of $\partial N$ which are disjoint from $\partial S$; a half vertex is {\em thin}
if $\overline{N}$ is $B^3$ and {\em fat} otherwise.}
\end{itemize}
The collection $\dee$ is {\em adequate} if all regular vertex labels are irreducible,
and all half vertex labels are irreducible and boundary irreducible.
\end{definition}

The analogue of Theorem~\ref{Laudenbach_theorem} holds for relative sum graphs:
\begin{theorem}[Laudenbach \cite{Laudenbach}]\label{relative_Laudenbach_theorem}
Two adequate collections $\dee,\dee'$ for $M$ determine relative sum graphs $G(\dee),G(\dee')$ which
are related by edge slides and homotopies of the thin part.
\end{theorem}

See also \cite{Casson_Gordon} for a discussion.

\begin{remark}
Theorem~\ref{relative_Laudenbach_theorem} is nontrivial even (or especially)
when $M$ is a handlebody, and implies (for instance) that any two meridian
systems $v,v'$ for a handlebody are related by a finite sequence of {\em Singer moves};
i.e.\ one disk in the system $v$ is replaced by a disk contained in the
complement of $v$ (see Singer \cite{Singer}).
\end{remark}

Criterion $C_2$ tries to distinguish elements of $\Mdot(S)$ using their relative sum graphs.
Two elements $(A,S),(B,S) \in \Mdot(S)$ are indistinguishable with respect to $C_2$
if there are adequate collections $\dee_A,\dee_B$ so that $G(\dee_A)$ and $G(\dee_B)$ are
{\em isomorphic}.

That is,
\begin{itemize}
\item{The families of disks in $\dee_A$ and $\dee_B$ intersect $S$ in the same set
of curves, up to isotopy; this determines a canonical isomorphism between the sets of
half edges in $G(\dee_A)$ and $G(\dee_B)$}
\item{There is an isomorphism between the graphs $G(\dee_A)$ and $G(\dee_B)$ which
is compatible with the canonical isomorphism on half edges, and which
takes fat (resp. thin) vertices
to fat (resp. thin) vertices and regular (resp. half) vertices
or edges to regular (resp. half) vertices or edges}
\end{itemize}

In words, two elements are indistinguishable with respect to criterion $C_2$ if they can be
cut up into irreducible pieces in a combinatorially isomorphic way (respecting the
intersection with $S$) and with the same
$S^3$ and $B^3$ terms, but without paying attention (yet) to the homeomorphism types of the
irreducible pieces which are not $S^3$ or $B^3$.

\begin{figure}[htpb]
\centering
\includegraphics[scale=0.75]{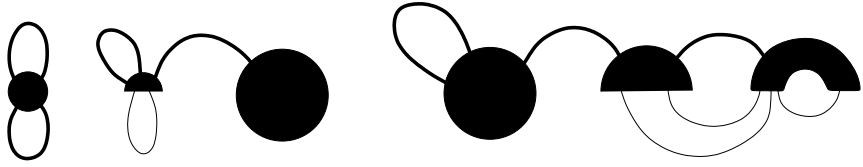}
\caption{An example of a relative sum graph $G(\dee)$. From the left, the first
component is an $S^1 \times S^2 \;\sharp\; S^1 \times S^2$ and the second
contains an $(S^1 \times S^2 \setminus D^3) \; \natural \; S^1 \times D^2$ summand.
Every edge is ``thin'' in the sense of Definition~\ref{relative_sum_graph_definition},
but in the figure, half-edges are drawn ``thinner'' than regular edges.} \label{half_sum_graph_figure}
\end{figure}

Figure~\ref{half_sum_graph_figure} shows an example of a relative sum graph.

In the sequel, for the sake of legibility, we will refer to relative
sum graphs simply as sum graphs.

\begin{lemma}[$c^2$-Lemma Schema]\label{c^2_Lemma_Schema}
Either
$$c^2(AB) < \max(c^2(AA), c^2(BB))$$
or $c^2(AB) = c^2(AA) = c^2(BB)$ and there are adequate collections $\dee_A,\dee_B$ of
disjoint spheres and disks in $(A,S),(B,S)$ respectively so that the sum graphs
$G(\dee_A)$ and $G(\dee_B)$ are isomorphic.
\end{lemma}

\begin{proof}
Let $\dee_A,\dee_B$ be adequate collections of disks and spheres for $A$ and $B$ respectively.
Without loss of generality we may assume that $c^1(AB)=c^1(AA)=c^1(BB)$ so that
$(A,S)$ and $(B,S)$ are indistinguishable with respect to criteria $C_0$ and $C_1$.
It follows that any essential simple closed curve in $S$ which
bounds a disk in $A$ also bounds a disk in $B$, so we may choose $\dee_A,\dee_B$
so that the disk components have the same boundaries in $S$, and the half edges
of $G(\dee_A)$ and $G(\dee_B)$ match up under the gluing.

Since the definition of adequate requires half vertex labels to be irreducible,
any half vertex label whose boundary contains a sphere must be $B^3$ (i.e.\ these are
exactly the set of thin half vertices). Note that
if $N$ is a component corresponding to a half vertex, then each boundary
component of $N$ is obtained from a component of $\partial N \cap S$ by filling
in boundary components of $\partial N \cap S$ with disks in $\dee$. Thus a boundary
component of $\overline{N}$ is a sphere if and only if the corresponding component of
$\partial N \cap S$ is a planar surface. In particular, thin half vertices of
$G(\dee_A)$ precisely match up with thin half vertices of $G(\dee_B)$.

Since disk components of $\dee_A$ and $\dee_B$ match up, the union $\ess = \dee_A \cup \dee_B$
is a pairwise disjoint collection of embedded spheres in $AB$.
\begin{claim}
$\ess$ is an adequate collection of spheres for $AB$.
\end{claim}
\begin{proof}
Components of $A\setminus \dee_A$ and $B\setminus \dee_B$ corresponding to regular vertices
are also components of $AB\setminus \ess$, so the labels associated to these vertices
of $G(\ess)$ are irreducible. So we just need to check irreducibility for labels of vertices
of $G(\ess)$ corresponding to unions of half vertices of $G(\dee_A)$ and $G(\dee_B)$.

Let $M_i$ be a set of components of $A\setminus \dee_A$
and $N_j$ a set of components of $B\setminus \dee_B$ corresponding to half vertices so
that the union $M = \cup_i M_i \cup_j N_j$ is a component of $AB\setminus \ess$. We need
to check that every sphere in $\overline{M}$  is inessential.
First notice that every boundary component of $N:=\cup_i \overline{M}_i \cup_j \overline{N}_j$
is a sphere, obtained from exactly one disk in each of $\dee_A$ and $\dee_B$, since
the disk components glue up exactly. Then observe that $\overline{M}$ is obtained from
$N$ by capping off each of these boundary components with a $3$-ball.

Now suppose $T$ is an embedded sphere in $\overline{M}$. We will show that $T$ bounds a
$3$-ball in $\overline{M}$.
Since $\overline{M}$ is obtained from $M$ by capping off boundary components with
$3$-balls, we may assume $T \subset M \subset AB$,
and we put $T$ in general position with respect to $S$. Since $\dee_A,\dee_B$ are adequate in $A,B$
respectively, 
we may boundary compress all components of $T\cap M_i$ in $M_i$ until each
component of $T\cap M_i$ is either boundary parallel into $S$
or parallel in $M_i$ to a disk $\delta_a$ in $\dee_A$.
Disks of the first type
can then be pushed across
$S$, reducing the number of intersections of $T \cap S$.
On the other hand disks parallel to an element $\delta_a \in \dee_A$ can be removed by
an isotopy (working innermost to outermost) across the 3-ball
in $\overline{M}$ capping off $\delta_a \cup \delta_b$, where $\delta_b \in \dee_B$
is the counterpart to $\delta_a$.
The isotopy results in
a new $T \subset M \subset AB$
with one fewer intersections with $S$. So after finitely many moves of this kind,
$T$ is disjoint from $S$, and is contained (without loss of generality) in some $M_i$.
But $\overline{M}_i$ is irreducible, so $T$ bounds a $3$-ball in $\overline{M}_i \subset \overline{M}$.
\end{proof}

It follows that the sum graphs $G(\dee_A)$ and $G(\dee_B)$ for $A$ and $B$ can be combined to give a
sum graph $G(\ess)$ for $AB$. More specifically,
\begin{itemize}
\item{Each half edge of $G(\dee_A)$ combines with its counterpart in $G(\dee_B)$ to produce
a (regular) edge of $G(\ess)$}
\item{Each thin half vertex of $G(\dee_A)$ combines with its counterpart in $G(\dee_B)$ to produce
a thin vertex of $G(\ess)$}
\item{Two or more fat half vertices from $G(\dee_A)$ and $G(\dee_B)$ combine to
produce a single fat vertex of $G(\ess)$}
\item{(Less interesting) The regular edges and vertices of $G(\dee_A)$ and $G(\dee_B)$
become edges and vertices of $G(\ess)$, preserving fatness/thinness in the case of vertices}
\end{itemize}
Note that the third bullet point is justified by the Claim above.
Similarly, $G(\dee_A)$ (respectively $G(\dee_B)$) can by combined with its mirror image
to get a sum graph $G(\ess_{AA})$ for $AA$ (respectively a sum graph $G(\ess_{BB})$ for $BB$).

\vskip 12pt

We now argue that either $G(\dee_A)$ can be chosen to be isomorphic to $G(\dee_B)$
or else the diagonal dominance inequality holds.

Let $a_1$ and $b_1$ be the number of fat vertices and fat half vertices of $G(\dee_A)$,
and define $a_2$ and $b_2$ similarly.
Then the number of fat vertices of $G(\ess_{AA})$ is $2a_1+b_1$, the number of fat vertices of
$G(\ess_{BB})$ is $2a_2+b_2$, and the number of fat vertices of
$G(\ess)$ is $a_1+a_2+b'$, where $b' \le \min(b_1, b_2)$.
It follows that the diagonal dominance inequality holds unless $a_1 = a_2$ and
$b' = b_1 = b_2$.
The latter equality forces the fat half vertices of $G(\dee_A)$ and $G(\dee_B)$
to have the same connectivity partition with respect to their common boundary.

Next we look at the second component of $c_2$, the first Betti numbers of $G(\ess_{AA})$,
$G(\ess_{BB})$ and $G(\ess)$.

Let $G_{1/2}(\dee_A)$ denote the subgraph of $G(\dee_A)$ consisting of half vertices and
half edges, and define $G_{1/2}(\dee_B)$ similarly. By the remark above, the subgraphs
$G_{1/2}(\dee_A)$ and $G_{1/2}(\dee_B)$ are isomorphic by an isomorphism which respects
the gluing (remember that half edges correspond to disks which are paired according
to Criterion $C_1$).

Because of the isomorphism of $G_{1/2}(\dee_A)$ and $G_{1/2}(\dee_B)$, the Betti numbers
$s$ can be computed  from the relative first
Betti numbers
$$s_A:=\dim(H_1(G(\dee_A), G_{1/2}(\dee_A); \R)), \;
s_B:=\dim(H_1(G(\dee_B), G_{1/2}(\dee_B); \R))$$
by the formula
$$s(AA) = s_A + s_A + C, \; s(AB) = s_A + s_B + C, \; s(BB) = s_B + s_B + C$$
where $C$ is a constant depending only on the (isomorphic) graphs $G_{1/2}(\dee_A)$
and $G_{1/2}(\dee_B)$.
If $s_A \ne s_B$, then the diagonal dominance inequality
is satisfied, so assume from now on that they are equal.

\vskip 12pt

The same argument can be applied to each connected component of $A$ and $B$ individually, using
Criterion $C_0$ to observe that the gluing of $A$ and $B$ respects their connectivity
partitions. In particular, connected components of $G(\dee_A)$ and $G(\dee_B)$
individually have isomorphic half subgraphs, the same number of fat vertices, and
the same Betti numbers. It follows that after edge slides, each component of
$G(\dee_A)$ is isomorphic to the corresponding component of $G(\dee_B)$ by an isomorphism
respecting the inclusions of $S$ into $A$ and $B$.
\end{proof}

\subsection{The complexity $c_3$}\label{c3_subsection}

We come now to the final piece of our complexity, $c_3$.  To define
$c_3$, we will leap ahead in the story and assume the following result, which is
assembled from results proved in \S\ref{c_S_section}--\ref{c_a_section}
and is itself proved in \S\ref{c_a_section}.

\begin{definition}
Let $\Pdot \subset \Mdot$ denote the set of closed, connected, oriented,
irreducible prime $3$-manifolds.
(We use ``irreducible'' here to exclude $S^1\times S^2$.)
\end{definition}

\begin{theorem}[Prime complexity]\label{prime_complexity_theorem}
There is a complexity function $c_p: \Pdot \to \Owe_p$ for some
totally ordered set $\Owe_p$, satisfying $c_p(M) = c_p(\overline{M})$,
and such that if $(A,S)$ and $(B,S)$ are {\em distinct} connected, orientable,
irreducible manifolds
with $S$ a non-empty incompressible surface not homeomorphic to $S^2$, then
$$c_p(AB) < \max(c_p(AA), c_p(BB))$$
\end{theorem}

\begin{remark}
The hypotheses on $(A,S)$ and $(B,S)$ imply that $AB,AA,BB$ are all in $\Pdot$.
\end{remark}

Deferring the definition of $c_p$ and assuming the proof of Theorem~\ref{prime_complexity_theorem},
we will now construct $c_3$ so that
$c^3 = c^2 \times c_3$ has the final property we demand of
$c := c^3$, namely:
$$c(AB) < \max(c(AA), c(BB))$$
for all {\em distinct} pairs $(A,S),(B,S) \in \Mdot(S)$.

\vskip 12pt

\begin{definition}
A {\em divisor} on $\Pdot$ is a function from $\Pdot$ to $\Z$ with finite support.
If $\Sigma$ is a divisor, write $\Sigma \ge 0$ if all the terms are non-negative.
\end{definition}

We sometimes use the circumlocution ``a finite subset of $\Pdot$ with multiplicity''
for a divisor $\Sigma \ge 0$. We extend $c_p$ to non-negative divisors, as follows:

\begin{definition}
If $\Sigma$ is a non-negative divisor, let $M_1,M_2,\ldots$ be the elements of $\Pdot$
in the support of $\Sigma$ listed with multiplicity, ordered
so that $c_p(M_i)$ is nonincreasing. Then define $\bar{c}_p$ to be the
lexicographic list
$$\bar{c}_p(\Sigma) = (c_p(M_1),c_p(M_2),\cdots)$$
\end{definition}

\begin{remark}
Note that by convention,
$\bar{c}_p(0) < \bar{c}_p(\Sigma)$ whenever $\Sigma > 0$; i.e.\ when the non-negative
divisor $\Sigma$ is nontrivial.
\end{remark}

The set $\Pdot$ admits an involution $\iota:\Pdot \to \Pdot$ given by reversing
orientations. Let $|\Pdot|$ denote the set of orbits of $\Pdot$ under this involution.
Let $p_1,p_2,\cdots$ list the elements of $|\Pdot|$. For each
divisor $\Sigma$ define $p_i(\Sigma)$ to be the number
of $M  \in \Sigma$ for which the orbit of $M$ under $\iota$ is equal to $p_i$.
Similarly, define $d_i(\Sigma)$ to be equal to $0$ if the cardinality of $p_i$
(as an orbit in $\Pdot$) is $1$, and otherwise to be equal to minus the
difference of the value of $\Sigma$ on the two
different elements of the orbit $p_i$. That is, if $P_i \in \Pdot$ is
some element in the orbit $p_i$, define
$$d_i(\Sigma) = - |\Sigma(P_i) - \Sigma(\iota(P_i))|$$
Observe that $d_i=0$ if $P_i = \iota(P_i)$, so actually this formula is correct
independent of the the cardinality of $p_i$.

\begin{definition}\label{prime_orientation_defect_definition}
For each divisor $\Sigma \ge 0$ define $c_\iota(\Sigma)$ to be the lexicographic list
$$c_\iota(\Sigma) = \text{ all the }(p_i(\Sigma),d_i(\Sigma))\text{'s listed in nonincreasing order }$$
\end{definition}

For $M \in \Mdot$, let $\Sigma(M)$ denote the non-negative divisor counting
the number of irreducible primes
in the prime decomposition of $M$ with multiplicity. Note that if $M$ is disconnected,
$\Sigma(M) = \sum_i \Sigma(M_i)$ where the $M_i$ are
the components of $M$.

With these preliminary definitions, we are in a position to define $c_3$.

\begin{definition}
$c_3(M)$ is the lexicographical tuple
$$c_3(M) = (\bar{c}_p(\Sigma(M)), \bar{c}_p(\Sigma(M_i)), c_\iota(\Sigma(M)), c_\iota(\Sigma(M_i)))$$
\end{definition}
Here the terms $\bar{c}_p(\Sigma(M_i))$ and $c_\iota(\Sigma(M_i))$ are shorthand for
the lexicographic lists of the complexities $\bar{c}_p$ and $c_\iota$ evaluated on
the components $M_i$ of $M$, ordered nonincreasingly.
Note that the order of the
$M_i$'s in the two cases might be different.

The terms $c_\iota$ become important when no non-trivial
gluing as in Theorem~\ref{prime_complexity_theorem} takes place.  In this case, we
need to look closely at prime multiplicity.

\begin{theorem}[$c^3$-Lemma Schema]\label{c3_theorem}
Let $(A,S)$ and $(B,S)$ be {\em distinct} elements
of $\Mdot(S)$. Then
\begin{equation}\label{c^3thmeqn}
c^3(AB) < \max(c^3(AA), c^3(BB))
\end{equation}
\end{theorem}

\begin{proof}
Suppose $c_0,c_1,c_2$ are all equal for $AA,AB,BB$ so that $A$ and $B$
cannot be distinguished by criteria $C_0,C_1,C_2$.

\begin{claim}
Either $\bar{c}_p(\Sigma(AB)) < \max(\bar{c}_p(\Sigma(AA),\Sigma(BB)))$ or else $A$ and $B$
pair up component by component in such a way that the fat regular vertices
of the sum graphs $G(\dee_A), G(\dee_B)$
correspond to sets of irreducible primes with pairwise equal
$c_p$-complexities, and the fat half vertices of $G(\dee_A)$ and $G(\dee_B)$ are paired by the
gluing in such a way that paired half vertices have identical labels.
\end{claim}
\begin{proof}
For each of $AA,AB,BB$, let $\Sigma^{1/2}(\cdot)$ be the divisor corresponding to
the set of manifolds appearing as labels of fat vertices in the sum graphs which
are obtained from fat half vertices in the relative sum graphs, and let
$\Sigma^r(\cdot)$ denote the other terms in $\Sigma(\cdot)$. Note that
$$\Sigma^r(\cdot) = \Sigma(\cdot) - \Sigma^{1/2}(\cdot)$$
as (non-negative) divisors.
Let $x_{AA} = \bar{c}_p(\Sigma^{1/2}(AA)),
u_{AA} = \bar{c}_p(\Sigma^r(AA))$ and $p_{AA} = \bar{c}_p(\Sigma(AA))$, and similarly
for $AA$ replaced by $AB$ or $BB$.

Ignoring orientation of factors at the moment, there is an equality
$\Sigma^r(AA) + \Sigma^r(BB) = 2\Sigma^r(AB)$
where we think of each $\Sigma^r$ as a divisor on $|\Pdot|$. It follows,
after interchanging the roles of $A$ and $B$ if necessary, that there is
an inequality
\begin{equation}\label{u_inequality}
u_{AA} \ge u_{AB} \ge u_{BB}
\end{equation}
where both inequalities are strict unless all three terms are equal.

We next turn to the terms $x_{AA}$, $x_{AB}$ and $x_{BB}$.
These correspond to gluings of fat half-vertices and satisfy
the hypotheses of Theorem~\ref{prime_complexity_theorem}.
By that theorem,
either paired fat
half vertices have labels which are diffeomorphic by a diffeomorphism preserving
the intersection with $S$, or else
\begin{equation}\label{x_inequality}
x_{AA} > x_{AB} \text{ or } x_{BB} > x_{AB}
\end{equation}
We consider these cases in turn.

\vskip 12pt
{\noindent \bf Case $x_{AA} > x_{AB}$:}
Combining with equation~\ref{u_inequality} gives $p_{AA} > p_{AB}$ and we are done.

\vskip 10pt

{\noindent \bf Case $x_{BB} > x_{AB}, \; u_{AB} = u_{BB}$:}
These two inequalities together imply $p_{BB} > p_{AB}$
and we are also done.

\vskip 10pt

{\noindent \bf Case $x_{BB} > x_{AB}, \; u_{AB} > u_{BB}$:}
Combining the second inequality with equation~\ref{u_inequality} gives $u_{AA} > u_{AB}$.

Let $u_{AA}'$ be the first term in $u_{AA}$ which differs from a term
$u_{AB}'$ in $u_{AB}$, so that there is an inequality
\begin{equation}\label{u_prime_inequality}
u_{AA}' > u_{AB}'
\end{equation}

\vskip 10pt

{\bf Subcase $x_{BB} > u_{AA}'$:}
That is, assume the first term in $x_{BB}$ is greater than $u_{AA}'$.
Combining with equation~\ref{u_prime_inequality} gives $x_{BB} > u_{AB}'$
and therefore $p_{BB} > p_{AB}$ and again we are done.

\vskip 10pt

{\bf Subcase $x_{BB} \le u_{AA}'$:}
Since this is a subcase of the case $x_{BB} > x_{AB}$, we get $x_{AB} < u_{AA}'$.
That is, the first term $x_{AB}'$ in $x_{AB}$ is already less than the first
term in $u_{AA}'$ in which case the $u$ terms dominate in $p$.
Since we are assuming throughout this case that $u_{AB} > u_{BB}$, one
inequality in equation~\ref{u_inequality} is strict, and therefore the other is too;
i.e.\ $u_{AA} > u_{AB}$. Since we have just argued that the $u$ terms dominate in $p$
in this subcase, we deduce $p_{AA} > p_{AB}$, so the claim is proved in this case too.

\vskip 12pt

Thus the only possibility is that labels of paired fat half vertices are diffeomorphic
by a diffeomorphism preserving their intersection with $S$, and
further that $u_{AA} = u_{AB} = u_{BB}$, as claimed.
\end{proof}


Repeating the preceding argument component by component, we either
obtain the desired inequality in the second term of $c_3$, or else there
are equalities
$$\Sigma^{1/2}(A_iA_i) = \Sigma^{1/2}(A_iB_i) = \Sigma^{1/2}(B_iB_i)$$ and
$$\bar{c}_p(\Sigma^r(A_iA_i)) = \bar{c}_p(\Sigma^r(A_iB_i)) = \bar{c}_p(\Sigma^r(B_iB_i))$$
for each pair of components $A_i,B_i$ which are glued up in $AB$.
Since paired fat half vertices in every component are diffeomorphic by a
diffeomorphism fixed on $S$, the only potential difference between $A$ and
$B$ is in the prime irreducible factors which do not intersect $S$; call
these the {\em free factors}.
By uniqueness of prime factorization for $3$-manifolds, $A$ and $B$ are
equal if and only if for each pair of corresponding components $A_i,B_i$
the set of free factors appearing in each are equal as divisors.

Ignoring questions of orientation of free factors for the moment,
it is clear that either the set of free factors in $A_i$ and in $B_i$ are equal,
or for some first $p_k \in |\Pdot|$, after possibly interchanging
$A$ and $B$, some component $A_i$ has at least as many copies
of $p_k$ as any component $B_j$, and more copies than $B_i$. It follows
that the maximal value of $p_k(\cdot)$ is achieved on more components
of $AA$ than of $AB$, so that $c_3(AB) < c_3(AA)$ with strict inequality
in either the third or fourth term.

However, we are working in the category of oriented $3$-manifolds and
we must take into account the fact that fat regular vertices of $B$
contribute terms in $AB$ and $BB$ which appear with reversed orientation.
The $d_i$ factors in the complexity $c_\iota$ favor a perfect balance in
the closed manifold between primes and their orientation-reverses.

Since all factors appearing in $\Sigma^{1/2}$ are amphichiral
(i.e.\ they are fixed points of $\iota$), the $d_i$ terms are maximized
only when every oriented free factor in each $A_i$ occurs in each $B_i$
with the same multiplicity. So either $c_3(AB) < \max(c_3(AA),c_3(BB))$
with strict inequality in some $d_i$ in the third or fourth term, or
each pair
$A_i$ and $B_i$ have exactly the same set of free factors with orientation
and multiplicity. In the first case, the theorem is proved. In the second case,
by what we have already proved about the $\Sigma^{1/2}$ factors,
this implies that $(A,S)$ and $(B,S)$ are equal as elements of $\Mdot(S)$, contrary
to hypothesis.
\end{proof}

This completes the proof of Theorem~A modulo the definition of $c_p$ and
the proof of Theorem~\ref{prime_complexity_theorem}. This will occupy
\S\ref{c_S_section}--\ref{c_a_section}.

\section{Seifert fibered factors}\label{c_S_section}

We now restrict attention to irreducible $3$-manifolds with incompressible boundary.

\begin{definition}
A {\em Seifert fibered space} is a compact $3$-manifold that admits
a foliation by circles. The foliation is called a {\em Seifert fibration}
and the circles are called the {\em fibers} of the Seifert fibration.
If $M$ is a Seifert fibered space, the space obtained by quotienting the
fibers to points is an orbifold $Q$ called the {\em base} of the fibration.
The orbifold points are also called {\em exceptional points}, and
the fibers lying over orbifold points are called {\em exceptional fibers}.
\end{definition}

Fundamental to the theory of Seifert fibered manifolds is the fact that
an irreducible $3$-manifold admits a natural decomposition into Seifert fibered
and atoroidal pieces:

\begin{definition}
Let $M$ be a closed, orientable, irreducible $3$-manifold. The {\em characteristic
submanifold of $M$}, denoted $\Sigma$, is a Seifert submanifold of $M$ (possibly disconnected
and with boundary) whose complement is atoroidal (and possibly disconnected), and which has
the smallest number of boundary components.
\end{definition}

Note that the boundary components of an orientable Seifert fibered manifold are all tori.

\begin{theorem}[JSJ decomposition, \cite{Jaco_Shalen}, \cite{Johannson}]\label{JSJ_theorem}
A closed, orientable, irreducible $3$-manifold has a characteristic submanifold which
is unique up to isotopy.
\end{theorem}

There is a relative version of this theorem for manifolds with incompressible
boundary, which we also need.

\begin{definition}\label{relative_JSJ_definition}
Let $M$ be a compact, orientable, irreducible $3$-manifold with incompressible boundary $S$.
The {\em characteristic submanifold of $M$}, denoted $\Sigma$, is a union of pieces, each
of which is one of the following
kinds:
\begin{enumerate}
\item{A Seifert submanifold disjoint from $S$ (call these {\em free factors})}
\item{A pair ($I$-bundle, $\partial I$-bundle) over a surface, where the intersection with $S$ is
the $\partial I$-bundle (call these {\em proper $I$-bundles})}
\item{A Seifert submanifold whose boundary intersects $S$ in a union of fibered annuli and tori}
\end{enumerate}
This union of pieces is uniquely determined by the property that its 
complement is atoroidal and acylindrical, and has the smallest possible number of frontier components.
\end{definition}

\begin{theorem}[relative JSJ decomposition, \cite{Jaco_Shalen}, \cite{Johannson}]\label{relative_JSJ_theorem}
A compact, orientable, irreducible $3$-manifold with incompressible boundary has
a characteristic submanifold which is unique up to isotopy.
\end{theorem}

\begin{remark}\label{fibering_convention}
There is some ambiguity in the fiber structure of a piece of the characteristic submanifold of
a manifold (with or without boundary). This ambiguity is discussed in the sequel, and is resolved
by certain conventions.
\end{remark}

\begin{remark}\label{fibered_annuli_condition}
The third case in Definition~\ref{relative_JSJ_definition} is somewhat nonstandard, insofar as we
insist that the fibering of a Seifert submanifold $X$ extend to a fibering of its intersection
$X \cap S$. This has the following consequence: if $T$ is a torus component of $\partial X$ 
intersecting $S$ in parallel annuli $A_i$, and the circle fibration of the $A_i$ does not extend over
$X$, we need to add a parallel copy of $T \times I$ as a piece of the characteristic submanifold
which ``insulates'' $\partial X$ from $S$. If this $T\times I$ is glued up along subannuli with nontrivial
Seifert fibered pieces on the other side of $S$, it will survive as part of a Seifert fibered component
of the JSJ decomposition of the closed manifold. Otherwise, it will be ``reabsorbed'' into the
boundary of $X$.
\end{remark}

\subsection{Surfaces of finite type}

The pieces which arise in the JSJ decomposition are canonical (up to isotopy) but
they typically have boundary, consisting of a union of tori. It is sometimes desirable
(especially when discussing the atoroidal pieces) to remove these boundary tori
and consider the open manifold which is the interior of the manifold with boundary.

It is therefore convenient to extend the definition of $\Mdot(S)$ to the case that $S$ is
a surface of finite type.

Let $S$ be an oriented surface of finite type (i.e.
$S$ is homeomorphic to a closed surface with finitely many points removed).
Let $\Mdot(S)$ denote the set of smooth pairs $(A,S)$ where
\begin{itemize}
\item{$A$ is an orientable $3$-manifold with $\partial A = S$,}
\item{$A \setminus S$
is homeomorphic to the interior of a compact $3$-manifold $\widehat{A}$, and}
\item{$\partial \widehat{A}$ decomposes as a union
$\partial \widehat{A} = \partial_v \widehat{A} \cup \partial_h \widehat{A}$
(the {\em vertical} boundary and the {\em horizontal} boundary respectively) such
that $\partial_h \widehat{A} = S$ and $\partial_v \widehat{A}$ is a finite union of tori and
annuli whose boundary compactifies $S$.}
\end{itemize}
up to the equivalence relation $(A,S) \sim (B,S)$ if there is a diffeomorphism $h:A \to B$
such that $hi = j$ where $i:S \to A$ and $j:S \to B$ are the respective inclusions.

\subsection{Classification and conventions}

In \S~\ref{Seifert_complexity_subsection} we will introduce the complexity term $c_S$ which treats both
closed and bounded oriented sufficiently large (to be defined below) Seifert fibered
spaces, but first we must discuss the classification of Seifert fibered spaces (in order to
be able to define the complexity) and describe our conventions for dealing with a few ``exceptional''
cases. As remarked above, the complexity must be defined for bounded as well as closed
manifolds. We may assume that the gluing surface $S$ is incompressible, since compressing disks are
treated by the (earlier) complexity term $c_1$.

It follows as in the proof of Theorem~\ref{c3_theorem} that
Seifert fibered pieces in the relative JSJ decomposition disjoint from the gluing surface $S$
(i.e.\ the free factors) occur in either $AA$ or $BB$ as many times as they
occur in $AB$, so it is only necessary to define $c_S$ on the kinds of
Seifert fibered pieces which arise by nontrivial gluing along an incompressible subsurface $S$.

The classification of orientable Seifert-fibered manifolds is well-known.
We follow the notation of Hatcher \cite{Hatcher}.

\begin{definition}[Notation for oriented Seifert fibered manifolds]\label{Seifert_notation}
The notation $M(\pm g,b;\alpha_1/\beta_1,\cdots,\alpha_k/\beta_k)$ denotes
a Seifert fibered manifold specified by the following properties.
The base surface $B$ has genus $g$, with sign $+$ if
$B$ is orientable and $-$ otherwise. The base surface
has $b$ boundary components. The fiber structure over the base is obtained
from a ``model'' $S^1$ bundle $E$ over $B$ with oriented total space by drilling
out $k$ fibers and gluing in solid tori with slopes $\alpha_i/\beta_i \in \Q$, and
where $\alpha_i/\beta_i \in \Z$ is allowed only if $b=0$ and $k=1$.
The ``model'' is unique if $b>0$.
If $b=0$ and $B$ is orientable set $E = S^1\times B$.
If $b=0$ and $B$ is nonorientable set
$E = (S^1\times \tilde{B})/\langle (\theta, \tilde{b}) \sim
(\bar{\theta}, {\tilde{b}}') \rangle$,
where $\tilde{B}$ is the orientation cover of $B$,
$\tilde{b}'$ is the covering translation of $\tilde{b}$,
and $\bar{\theta}$ is the complex conjugate of $\theta$.
\end{definition}

There is some redundancy in this notation, captured in the following proposition.

\begin{proposition}[Classification theorem, \cite{Hatcher} Prop. 2.1]\label{seifert_classification_theorem}
Two Seifert fiberings $$M(\pm g, b ;\alpha_1/\beta_1,\cdots,\alpha_k/\beta_k), \quad
M(\pm g, b;\alpha_1'/\beta_1',\cdots,\alpha_k'/\beta_k')$$ are isomorphic by an
orientation (and fiber) preserving diffeomorphism if and only if, after possibly
permuting indices, $\alpha_i/\beta_i = \alpha_i'/\beta_i' \pmod 1$ for each
$i$, and, in case $b=0$, if there is an equality $\sum_i \alpha_i/\beta_i = \sum_i
\alpha_i'/\beta_i'$.
\end{proposition}

If $b=0$, the sum $\sum_i \alpha_i/\beta_i$ is an invariant, called
the {\em Euler number} of the fibering.

\begin{definition}
\label{sufficiently_large_defn} A Seifert fibered space, with or without
boundary, is {\em sufficiently large} if it contains an oriented
incompressible surface (not equal to $S^2$ or $D^2$).  Otherwise, the Seifert fibered space is
called {\em small}.
\end{definition}

\begin{theorem}[Waldhausen, \cite{Waldhausen}]\label{Waldhausen_thm}
Suppose $M$ is a Seifert fibered space.
Then all incompressible surfaces in $M$ can be isotoped
to be either {\em horizontal} (i.e.\ transverse to the fibers) or
{\em vertical} (i.e.\ a union of fibers).
\end{theorem}

The Euler number is an obstruction to the existence of a horizontal surface.
In fact, the following is well-known:

\begin{lemma}[\cite{Hatcher} Prop. 2.2]\label{obstruction_to_surface}
If $b>0$, horizontal surfaces always exist (although they may be
disks if $M$ is not sufficiently large), and if $b=0$, then horizontal surfaces
exist if and only if the Euler number is zero.
\end{lemma}

Except for some special cases, the underlying $3$-manifold of
an orientable Seifert fibered space admits a unique fibering up to isomorphism:

\begin{theorem}[Exceptional list, \cite{Hatcher} Thm. 2.3]\label{exceptional_list_theorem}
Seifert fiberings of orientable Seifert manifolds are unique up to
isomorphism, with the following exceptions:
\begin{enumerate}
\item{$S^1 \times D^2$}
\item{$S^1 \tilde{\times} S^1 \tilde{\times} I$}
\item{$S^3$, $S^1 \times S^2$, lens spaces}
\item{$M(0,0;1/2,-1/2,\alpha/\beta) = M(-1,0;\beta/\alpha)$}
\item{$S^1 \tilde{\times} S^1 \tilde{\times} S^1$}
\end{enumerate}
\end{theorem}

Almost all the exceptions have base space a disk, a
sphere, or a projective plane with at most $2$, $3$, or $1$
exceptional fibers respectively. Such manifolds are all small
and will not arise by gluing irreducible, boundary irreducible components of
$A$ ($B$) along torus or annulus subcomponents of $S$ (although they might arise
as a factor in the relative JSJ decomposition). However, some cases require special treatment.

\medskip

{\noindent \bf Case $S^1 \times D^2$:} A solid torus (case (1) in Theorem~\ref{exceptional_list_theorem}) 
can arise as a piece of the relative JSJ decomposition
if its frontier consists of at least two annuli, or of one annulus whose core represents
a proper multiple of the core of the solid torus. If $A_i$ denote the (oriented) annuli in
the frontier of the solid torus, the cores of the $A_i$ represent multiples of the core of the
solid torus in homology; define the {\em degree} of the solid torus to be the sum of these numbers.
Choose orientations so that the degree is positive. We must have degree at least $2$, or else the
frontier is boundary parallel. We make the following convention about the fiber structure on a solid torus:

\begin{convention}\label{solid_torus_fiber_convention}
Let $S^1 \times D^2$ be a piece in the relative JSJ decomposition.
If the piece has degree $2$, and intersects $S$ in two annuli, each homotopic 
to the core of the solid torus, fiber the solid torus as an $I$-bundle over an annulus. Otherwise,
fiber the solid torus, using at most one exceptional fiber as the core, in the unique manner compatible
with a fibering of its frontier (or, equivalently, its intersection with $S$) by circles.
\end{convention}

It is worth remarking that the convention distinguishes between the two cases of degree $2$. 
Both cases admit the structure of an $I$-bundle or an $S^1$ bundle compatibly with their intersection
with $S$. The first case we give the structure of an $I$ bundle over an annulus (rather than a product
foliation of a solid torus by circles); the second case we give the structure of an $S^1$ bundle over a
disk with one exceptional fiber of kind $1/2$ (rather than a twisted $I$ bundle over a M\"obius strip).

\medskip

{\noindent \bf Case $S^1 \tilde{\times} S^1 \tilde{\times} I$:} A twisted $I$-bundle over 
a Klein bottle (case (2) in Theorem~\ref{exceptional_list_theorem}) can
arise in the JSJ decomposition of the closed manifold either as a free factor (in which case the
fiber structure is irrelevant), or by gluing two solid tori of degree $2$ of the same kind. It might
arise by gluing solid tori in one of two ways: either by gluing two $\text{annulus} \times I$ 
along pairs of boundary annuli with a ``twist'' (that interchanges inside and outside), 
or by gluing two $M(0,1;1/2)$ along fibered annuli in their boundaries. 
Convention~\ref{solid_torus_fiber_convention} insists in either case that the result 
has the fiber structure $M(0,1;1/2,-1/2)$. We therefore make the following:

\begin{convention}\label{twisted_bundle_Klein_convention}
Every $S^1 \tilde{\times} S^1 \tilde{\times} I$ in the closed manifold is given the fiber structure
$M(0,1;1/2,-1/2)$.
\end{convention}

An $S^1 \tilde{\times} S^1 \tilde{\times} I$ can also occur as a piece in the relative JSJ decomposition,
intersecting $S$ in a union of annuli. If these annuli are fibered compatibly with one of the fiber
structures $M(0,1;1/2,-1/2)$ or $M(-1,1;)$ we give the piece this fiber structure. Otherwise,
the piece does not really intersect $S$ at all, but is insulated by a $T^2 \times I$ as in
Remark~\ref{fibered_annuli_condition}.

Finally, 
\begin{convention}\label{twisted_Klein_torus_convention}
An $S^1 \tilde{\times} S^1 \tilde{\times} I$ that occurs as a piece in the relative JSJ
decomposition with its entire boundary on $S$ is given the fiber structure $M(0,1;1/2,-1/2)$.
\end{convention}

There is a further ambiguity, that a Seifert fibered manifold may
admit isomorphic but non-isotopic Seifert fiberings. If $M$ is
sufficiently large, the fibers determine a central $\Z$ subgroup of
$\pi_1(M)$, or a normal $\Z$ subgroup if the monodromy on fibers is
nontrivial. If the center (after passing to a double cover if
necessary) is at least $2$-dimensional, $M$ is either virtually $T^2
\times I$ or virtually $T^3$. From this one can deduce the following
standard fact:

\begin{theorem}[Big center]\label{mapping_class_exceptions}
Let $M$ be a sufficiently large
Seifert fibered space which is not on the exceptional list from
Theorem~\ref{exceptional_list_theorem}. Then either $M$ is one of
$T^2 \times I$ or $T^3$, or any automorphism of $M$ is isotopic to
an automorphism which permutes the fibers, and either fixes the orientation
of the base and fiber (if orientable), or simultaneously reverses both.
\end{theorem}

The manifolds $T^3$ and $T^2 \times I$ will be handled as a special case; for the moment we just
make the following:

\begin{convention}\label{torus_product_convention}
A properly embedded $T^2 \times I$ factor in the relative JSJ decomposition (i.e. with 
$T^2 \times \partial I \subset S$) is fibered as a product $I$-bundle over $T^2$. A
non-properly embedded factor with exactly one boundary component $\subset S$ is
thought of as being fibered by circles in an unspecified manner (since any fibration on a single
boundary component extends uniquely to a product fibration on the entire $T^2 \times I$).
\end{convention}

Other than free factors, Seifert fibered spaces may arise in the pairing
either by gluing two Seifert fibered spaces $A$
and $B$ (along part of their boundaries) so as to match the fiber
structures, or by gluing two $I$-bundles together to match end points
of intervals.  In the latter case, unless the base has positive
(orbifold) Euler characteristic, the result will be sufficiently large.
In both cases the surface $S$ along which the gluing occurs has a nontrivial fundamental
group which injects into $A$ (and $B$).

If gluing gives rise to a non-free Seifert fibered factor,
Theorem~\ref{Waldhausen_thm} implies that $S$ is either a torus or Klein
bottle, or else $S$ is horizontal. If the factor has no boundary,
then if $S$ is horizontal, the Euler number must vanish. This rules
out all closed exceptional cases in Theorem~\ref{exceptional_list_theorem}
except for $S^1 \tilde{\times} S^1 \tilde{\times} S^1$, since they either
contain no incompressible torus, or have nonvanishing Euler number.

This leaves the following exceptional manifolds:
$$T^2 \times I, \quad T^3, \quad S^1 \tilde{\times} S^1 \tilde{\times} S^1$$
We will see in the proof of Lemma~\ref{SF_DD_lemma} that with Convention~\ref{torus_product_convention},
these exceptional cases do not cause a problem.
In every other case, by Theorem~\ref{exceptional_list_theorem} and
Theorem~\ref{mapping_class_exceptions}, a non-free Seifert fibered factor arising
in $AA,AB$ or $BB$ has a unique Seifert fibered structure, up to (orientation-preserving)
fiber-permuting automorphisms which either simultaneously preserve, or
simultaneously reverse the orientations of the fibers and the base.

\subsection{Seifert fibered complexity}\label{Seifert_complexity_subsection}

We are now ready to give the definition of the complexity function $c_S$,
discussing first the complexity for {\em connected} Seifert fibered manifolds (possibly with boundary).

\begin{definition}\label{c_csf_defn}
Let $X$ be a connected, oriented, sufficiently large Seifert fibered space.
Let $Q$ denote the base orbifold of the Seifert fibering (we cannot use the
notation $B$, for obvious reasons, so think $Q$ for ``quotient'').
Define the complexity $c_{CS}(X)$ to be the ordered $4$-tuple
$$c_{CS}(X)=(b, -\chi^o(Q), \#_{\text{sing}}, -|e|)$$
where $b$ is the number of boundary components, $-\chi^o(Q)$ is (minus)
the orbifold Euler characteristic of the base (see \cite{Thurston_notes}, Definition~13.3.3), 
$\#_{\text{sing}}$ depends on the singular fibers, and $-|e|$ is (minus) the absolute value of
the Euler number, when $b=0$. 
All the orbifolds $Q$ which arise are ``good'', meaning that they have the form
$Q = \Sigma/F$, a surface $\Sigma$ modulo a finite group action.
The orbifold Euler characteristic of $Q$ is defined as $\chi^o(Q) = \chi(\Sigma)/|F|$.
The term $\#_{\text{sing}}$ is itself
an ordered list of terms defined as follows. Choose (in advance) some total
ordering of $(\Q/\Z - 0)/\pm 1$, i.e.\ non-zero elements of $\Q/\Z$ up to sign.
This ordering should not depend on $X$. For each term $\pm \alpha/\beta$
list the pair $(s(\pm \alpha/\beta),d(\pm \alpha/\beta))$ where
$s(\pm \alpha/\beta)$ is the number of singular fibers of type $\pm \alpha/\beta$,
and $d(\pm \alpha/\beta)$ is the {\em defect}, i.e.
$$d(\pm \alpha/\beta) : = - | \# (\alpha/\beta  \text{ fibers}) - \# ( -\alpha/\beta \text{ fibers}) |$$
(note that $d(\pm 1/2) = 0$).
The term $\#_{\text{sing}}$ is the ordered list of ordered pairs $(s(\pm \alpha/\beta),d(\pm \alpha/\beta))$
over all $(\Q/\Z - 0)/\pm 1$.
\end{definition}

The function $c_{CS}$ can be extended to irreducible manifolds with 
boundary a union of tori
by setting it equal to a formal minimum value ``$-\infty$'' on
any manifold which is not a sufficiently large Seifert fibered space.
For example, $c_{CS}(X) = -\infty$ if $X$ is a hyperbolic manifold.

The pairs $(s,d)$ in $\#_{\text{sing}}$ should be compared with the terms
$(p_i(\Sigma),d_i(\Sigma))$ in the Definition~\ref{prime_orientation_defect_definition}
of $c_\iota$, which serve an analogous purpose, namely to define a complexity
which is maximized on the biggest number of
objects appearing in orientation-reversed pairs.

\begin{definition}
\label{c_SF_defn}
Let $X$ be a closed, connected, oriented, irreducible $3$-manifold, and let $X_i$
be the pieces of $X$ obtained by the JSJ decomposition. Order the
$X_i$ by (decreasing) complexity $c_{CS}$, and let $c_S(X)$ be a tuple
$$c_S(X) = \left(m, -m', \sum_i -\chi^o(Q_i), \sum_i \text{genus}(Q_i), n, c_{CS}(X_1),\cdots\right)$$
where the description of the terms is as follows. The first term $m$
is the maximum number of independent tori in $X$ (i.e.\ the maximum
number of pairwise disjoint pairwise non-isotopic incompressible
tori). In the second term $m'$ is the number of tori in the JSJ
decomposition of $X$. The third term is the sum of $-\chi^o(Q_i)$
over the Seifert fibered $X_i$, and the fourth term is the sum of
the genera of the $Q_i$. The fifth term $n$ is  the number of
Seifert fibered pieces in $X$, and the remaining terms are the
ordered list of the $c_{CS}$ complexity of the $X_i$. To compare the
complexity of two lists of different length, pad the shorter list by
$- \infty$ symbols if necessary.
\end{definition}

Given an element $(A,S) \in \Mdot(S)$ where $S$ is incompressible, and
$A$ is irreducible, by Theorem~\ref{relative_JSJ_theorem},
there is a decomposition of $A$ (along essential annuli properly embedded
in $A$ and essential tori disjoint from $S$)
into submanifolds $A_F$, $A_I$ and $A_C$,
where each component of $A_F$ is Seifert fibered, each component of $A_I$ is
a properly embedded essential $I$-bundle, and
every component of $A_C$ is atoroidal and acylindrical. Further, let
$A_S := A_F \cup A_I$.

\begin{remark}
We reinterpret our conventions in this language. By Convention~\ref{solid_torus_fiber_convention}, 
a solid torus is a component of $A_I$ if and only if it is a product $\text{annulus} \times I$ 
with $\text{annulus} \times \partial I$ contained in $S$; otherwise it is a component of $A_F$.
By Convention~\ref{twisted_Klein_torus_convention}, an $S^1 \tilde{\times} S^1 \tilde{\times} I$ with
entire boundary on $S$ is a component of $A_F$. Finally, by
Convention~\ref{torus_product_convention}, a proper $T \times I$ is in $A_I$, whereas a
torus neighborhood of a torus component of $S$ which bounds
an atoroidal, acylindrical submanifold is a component of $A_F$.
\end{remark}

We are now in a position to prove diagonal dominance for $c_S$.

\begin{lemma}[$c_S$ Lemma Schema]\label{SF_DD_lemma}
Let $(A,S)$ and $(B,S)$ be elements of $\Mdot(S)$ where $S$ is incompressible,
and both $A$ and $B$ are irreducible. Then there is
an inequality
$$c_S(AB) \le \max(c_S(AA), c_S(BB))$$
which is strict unless $S \cap A_S = S \cap B_S = :S_S$,
and $(A_S,S_S) = (B_S,S_S)$ as elements of $\Mdot(S_S)$
\end{lemma}

\begin{remark}
Two copies of $\text{surface} \times I$ can be glued together to
produce a manifold with a Solv or hyperbolic structure, so it is crucial that $c_S$
occur before $c_h$ in our complexity. 
\end{remark}

\begin{proof}
Free factors in the JSJ decomposition of highest complexity are most
common in one of the doubles, by repeating the argument in the proof
of Theorem~\ref{c3_theorem}. So it suffices to restrict attention to
the complexity of non-free factors, which we now do.

\medskip

The gluing of $A_S$ to $B_S$ does not in general match $A_S \cap S$ with
$B_S \cap S$. For instance, an $I$-bundle component of $A_S$ might be
glued to an acylindrical component of $B_C$, and the union would be
part of a hyperbolic piece of the JSJ decomposition of $AB$. The terms
$m,-m'$ can almost be recovered (with the same ordering) from the later
terms in the complexity, except when distinct torus boundary
components of the pieces $X_i$ are isotopic in such a way that the fibering
is not compatible under the isotopy. More precisely, these terms prefer
isotopic tori in different $X_i$ to be fibered in the same way. These
terms are only significant when $A_F$ and $B_F$ share common boundary tori in $S$,
and are maximized when the fiberings on these tori in $A_F$ and in $B_F$ match
up, and the pieces can be glued into a single component of the JSJ decomposition
of $AB$.

In general, let $X$ be a non-free Seifert fibered component of the JSJ decomposition of
$AB$. If $S \cap X$ is horizontal, then $X$ is a maximal union of $I$-bundle
pieces of $A_I$ and $B_I$, glued along incompressible subsurfaces of $S$.
Let $X_h$ be the union of such non-free Seifert fibered components, and let
$Q$ be the base of the fibration.
Let $Q_i$ be the bases of the fibration of
the pieces in $A_I$ and $B_I$ which are glued together to make $X_h$.
Since $S$ is incompressible, the contribution
of $X_h$ to $m$ and $m'$ can be recovered from the type and multiplicity of the
$Q_i$, and the order on $(m,m')$ is compatible with the order induced from
lower complexity terms; in other words, the contribution of $X_h$ to the $m,m'$
terms is maximized when the lower order complexity terms are maximized.

There is a transverse geometric structure on the Seifert fibration of
(each component of) the manifold $X_h$,
pulled back from a geometric structure on $Q$ (see e.g.\ \cite{Thurston_notes}, \S4.8
for a discussion of transverse geometric structures on Seifert fibrations).
Away from the singular points, the base $Q$ decomposes into
subsurfaces $Q^n$ where $n$ counts the number of times the corresponding
fiber intersects $S$ (equivalently, $n$ counts the number of $I$-fibers of
$A_I$ and $B_I$ pieces which make up the given $S^1$ fiber). By Gauss-Bonnet there
is a formula
$$\chi^o(Q) = \sum n\chi^o(Q^n)$$

Write $Q = Q_h \cup Q_e \cup Q_s$ where the subscripts denote the hyperbolic, Euclidean,
and spherical components of $Q$. Since each component of $Q$ is covered by an
essential subsurface of $S$, it follows that $Q_s$ is empty.
Since each $n \ge 2$, it follows that $-\chi^o(Q_h)$ is maximized only
when $Q_h = Q^2_h$, and $Q_h$ is double covered by the hyperbolic part of the
subsurface $A_I \cap S \cap B_I$, which must be as big as possible (and therefore the
hyperbolic components of $A_I \cap S$ are equal to the hyperbolic components of $B_I \cap S$).
Finally, $\sum_i \text{genus}(Q_i)$ and $n$
are maximized only when $Q_e = Q^2_e$ and when (as above) the Euclidean
part of $A_I \cap S$ and $B_I \cap S$ is as big as possible (after the hyperbolic
parts are already matched up). In other words, $A_I \cap S = B_I \cap S$.

\medskip

If $S \cap X$ is vertical, then $X$ is a maximal union of
Seifert fibered pieces of $A_F$ and $B_F$, glued along
fibered tori or annuli contained in $S$, possibly together with a number of
$T^2 \times I$ factors in $A_F$ or $B_F$ with some $S^1$ bundle structure
(rather than with their conventional $I$-bundle structure).
Note that an annulus component
of $A_F \cap S$ cannot be contained in a torus
component of $B_F \cap S$, or else a neighborhood of this torus (in $A$)
could be glued to $A_F$, thereby enlarging it, contrary to the definition
of characteristic submanifold.

Let $X_v$ be the union of such
non-free Seifert fibered components, and let $Q$ be the base of the fibration.
Then
\begin{equation}\label{Euler_sum_formula}
\chi^o(Q) = \sum \chi^o(Q_i) - \#(\text{annuli})
\end{equation}
where the $Q_i$ are the Euler characteristics of the bases of the various
terms in $A_F$ and $B_F$,
and $\#(\text{annuli})$ counts the number of annuli in $S$ along which
they are glued. It follows that $-\chi^o(Q)$
is maximized when they are glued along as many annuli
as possible, in particular, when the fibered annuli boundary components of
the $A_F$ and $B_F$ pieces match up.

Gluing two fibered torus boundary components together does not
change $\chi^o$, but it either increases genus or decreases $m'$
(and therefore increases $-m'$). Gluing on a proper $T^2 \times I$
factor at both ends does not change $\chi^o$ or genus, but it {\em
reduces} the total number of $T^3$ factors which occur in $X_h$. So
if the leading terms in $c_S$ are maximized, it follows that all
proper $T^2 \times I$ factors are actually already contained in
$X_h$, and never appear in $X_v$, and therefore in this case $A_F \cap S = B_F \cap S$ 
with the same fiber structure.

\begin{remark}
It is worth remarking at this point that any $M(-1,1;)$ piece in the relative JSJ decomposition
that is not glued up will survive as a component of the JSJ decomposition of the closed manifold,
but with the ``wrong'' Seifert fibered structure (i.e. contrary to 
Convention~\ref{twisted_bundle_Klein_convention}). However, we have just showed that
such unglued pieces never occur in the gluing of maximal complexity, and therefore the convention may
(and will) be assumed in the sequel. The $M(-1,1;)$ fiber structure is the right one when such pieces
are glued up to other pieces in the JSJ decomposition of the closed manifold.
\end{remark}

It follows that either the inequality is strict at one of the
first four terms of $c_S$, or else we are in the case that $A_I \cap S = B_I \cap S$
and $A_F \cap S = B_F \cap S$, and moreover that the fiber structures on
$A_F \cap S$ and on $B_F \cap S$ agree. The fifth term $n$ in $c_S$ is
maximized only when the connectivity partitions of $A_I \cap S$ into
$A_I$ and of $A_F \cap S$ into $A_F$ agrees with that of $B$ (see
Definition~\ref{connectivity_partition_definition} and
Lemma~\ref{c^0_Lemma_Schema}). Therefore non-free Seifert fibered
components of $AB$ are all made from exactly one component of each
of $A_I$ and $B_I$, or one component of each of $A_F$ and $B_F$,
whose intersections with $S$ are equal.

When two $I$-bundles are glued, we have already seen that $Q = Q^2$,
in other words, that every circle in the glued up component is a
union of exactly two intervals, one in $A$ and one in $B$. This
pairing of intervals gives an identification of $A_I$ with $B_I$
rel their intersection with $S$.

It remains to check gluing of $F$-pieces. The situation is very much
like the complexity for gluings of surfaces considered in
\cite{FKNSWW}. The basic idea is that maximizing $(b,-\chi^o(Q))$
component by component (always in decreasing order) forces the most
complicated $F$ pieces in $A$ to be glued to pieces of the same
complexity in $B$ term by term, and by induction, every piece in
$A_F$ must be glued to a piece in $B_F$ with the same number of
boundary components and with base of the same Euler characteristic.
This nearly proves the Lemma, except that different topological
orbifolds can have the same number of boundary components and the
same $-\chi^o$. This ambiguity is resolved by further consideration
of the $\#_{\text{sing}}$ term, which precisely favors the double.

There is one final piece of ambiguity for closed sufficiently large
Seifert fibered spaces, namely the Euler number of the fibering. Doubled
pieces in $AA$ and $BB$ obviously satisfy $e=0$, and therefore maximize $-|e|$.
It follows that pieces in $AB$ must also have $e=0$, and therefore we
are done, by Theorem~\ref{seifert_classification_theorem}.
\end{proof}

\begin{remark}
Note that $S^1 \tilde{\times} S^1 \tilde{\times} S^1$  factors can arise either
in $X_h$ (as a union of two copies of $S^1 \tilde{\times} S^1 \tilde{\times} I$
glued along their torus boundaries) or in $X_v$ (by gluing two copies of
$M(0,1;1/2,-1/2)$ along their torus boundaries). In either case, the complexity
is maximized only when these factors pair up, i.e.\ when $(A_S,S_S) = (B_S,S_S)$.
\end{remark}

\section{Hyperbolic factors}\label{c_h_section}

We assume the reader is familiar with fundamental facts from the Thurston theory of
hyperbolic $3$-manifolds. We use this material without comment throughout this section,
justifying only facts or claims which are new or unfamiliar.
A basic reference for this material is Thurston's notes \cite{Thurston_notes}.

We also make some use of the theory of (stable) minimal surfaces, especially
in hyperbolic $3$-manifolds, and the theory of Ricci flow.
A reference for the first is \cite{Colding_Minicozzi}, and a reference for the second
is \cite{Chow_Knopf}.

\subsection{Closed hyperbolic case}

The first piece of the hyperbolic complexity is $-\vol$, where $\vol$ denotes hyperbolic
volume. For closed manifolds, the justification for using this term is the following theorem
of Agol--Storm--Thurston:

\begin{theorem}[Agol--Storm--Thurston, \cite{AST}]\label{AST_thm}
Let $S$ be a closed orientable surface so that each component has negative Euler characteristic,
and let $(A,S),(B,S) \in \Mdot(S)$ be irreducible,
atoroidal and acylindrical. Then $AA,AB,BB$ admit unique complete hyperbolic structures, and
either
$$2\vol(AB) > \vol(AA) + \vol(BB)$$
or else $\vol(AB) = \vol(AA) = \vol(BB)$ and $S$ is totally geodesic in $AB$.
\end{theorem}

\begin{remark}\label{geodesic_remark}
Since by hypothesis $A$ and $B$ are atoroidal, irreducible and acylindrical,
Thurston's hyperbolization theorem for Haken manifolds (see \cite{Otal_hyperbolization} for
an exposition)
implies that both $A$ and $B$ admit unique complete hyperbolic structures with
totally geodesic boundary. Moreover, $S$ is totally geodesic in $AB$ if and only
if the totally geodesic metrics on $S$ inherited from the hyperbolic structures on
$A$ and $B$ are isometric by an isometry isotopic to the identity.
\end{remark}

A very brief outline of the proof in \cite{AST} is as follows.
First, construct the hyperbolic structure on $AB$.
Since both $A$ and $B$ are acylindrical and atoroidal,
$AB$ is atoroidal and such a hyperbolic structure exists,
by Thurston's hyperbolization theorem for Haken manifolds.
Next, find a minimal surface representative of the surface $S$ in the hyperbolic
manifold $AB$. Since $S$ is incompressible in $AB$, it admits a stable minimal
representative in its isotopy class, by the well-known theorem of Meeks--Simon--Yau
(i.e.\ the main theorem of \cite{MSY}). Cutting $AB$ open along $S$ produces two
Riemannian manifolds with boundary $A,B$. The metrics on the
doubles $AA,BB$ are not Riemannian unless
$S$ is totally geodesic, but are $C^0$ Riemannian --- i.e.\ defined by a $C^0$ section of the
bundle of symmetric $2$-tensors. Since $S$ is minimal in $AB$, it has
vanishing mean curvature. This implies, by arguments of
Miao \cite{Miao} or Bray \cite{Bray}, arising from work on the Riemannian Penrose inequality,
that the metrics on $AA$ and $BB$ can be
uniformly approximated by $C^\infty$ Riemannian metrics satisfying uniform pointwise
scalar curvature bounds $R\ge -6$. Technology due to Miles Simon \cite{Miles_Simon}
lets one apply short time Ricci flow to such singular metrics, which preserves the pointwise scalar
curvature bounds. Then Perelman's monotonicity formula for Ricci flow with
surgery (\cite{Perelman_I, Perelman_II}) implies the global inequality in
the theorem.

In \S\ref{noncompact_volume_section} we will generalize this theorem to the case that
$S$ is a surface of finite type. There are a number of new analytic details which arise in
this case, but we otherwise stay very close to the argument of \cite{AST}.

\subsection{Isometric gluing}

The next piece of the hyperbolic complexity is the {\em length spectrum}.

\begin{definition}\label{length_spectrum_definition}
Let $M$ be a finite volume hyperbolic $3$-manifold. The {\em real hyperbolic
length spectrum} of $M$, denoted $\sigma(M)$, is a lexicographic
tuple of real numbers $\sigma(M) = (\ell_1(M),\ell_2(M),\cdots)$
which lists (with multiplicity) the length of closed geodesics in $M$ whose
complex length is real, ordered nondecreasingly by length.
If this spectrum is finite (or empty),
pad this tuple with a string of $\infty$ symbols.

Say $\sigma(M) > \sigma(N)$ if $\ell_i(M)=\ell_i(N)$ for all $i<n$ but
$\ell_n(M) < \ell_n(N)$.
\end{definition}

\begin{warning}
Note the sign convention in Definition~\ref{length_spectrum_definition}.
The manifold with the {\em smaller} geodesics has the {\em larger} value of $\sigma$.
\end{warning}

\begin{lemma}\label{spectrum_lemma}
Suppose $(A,S),(B,S) \in \Mdot(S)$ admit finite volume complete
hyperbolic structures with totally geodesic boundary $S$
such that the two induced hyperbolic structures on $S$ are isometric 
by an isometry isotopic to the identity. Then either
$$\sigma(AB) < \max(\sigma(AA),\sigma(BB))$$
or $(A,S) = (B,S)$.
\end{lemma}
\begin{proof}
Let $\ell_1$ be the shortest length in $\sigma(AB)$. We will
derive a formula for the multiplicity of $\ell_1$.

Let $\sigma(AB,\ell_1)$ denote the multiplicity of $\ell_1$ in $\sigma(AB)$ and
define $\sigma(AA,\ell_1)$ and $\sigma(BB,\ell_1)$ similarly.

The formula clearly will have the form
\begin{equation}\label{1_multiplicity_AB_formula}
\sigma(AB,\ell_1) = k_1^A + k_1^B - k_1^S + \epsilon'
\end{equation}
where $k_1^A$ is the number of closed hyperbolic geodesics of length
$\ell_1$ in $A$ (and similarly for $k_1^B$),
$k_1^S$ is the number of closed hyperbolic geodesics of length
$\ell_1$ in $S$, and $\epsilon'$ is the number of geodesics in $AB$ of
length $\ell_1$ which are not contained on either side.
The $k_1^*$ terms come from the inclusion/exclusion formula. The $\epsilon'$
term is slightly more subtle to determine.

Define $\epsilon_A$ to be the number of proper essential
arcs in $A$ of length $\ell_1/2$, and define $\epsilon_B$ similarly.
If either $A$ or $B$ contains an essential arc shorter than $\ell_1/2$,
doubling produces a geodesic in either $AA$ or $BB$ with
length strictly less than $\ell_1$ and we would be done. So without loss
of generality we can assume that every proper essential
arc in either $A$ or $B$ of length at most $\ell_1/2$ has length
exactly $\ell_1/2$, and is totally geodesic and perpendicular to the boundary.
Moreover, a geodesic loop of length $\ell_1$ in $AB$
can consist of at most one arc on each side. Hence
$\epsilon' \le \min(\epsilon_A,\epsilon_B)$.

It follows that there is a formula
\begin{equation}\label{1_multiplicity_AA_formula}
\sigma(AA,\ell_1) = 2k_1^A - k_1^S + \epsilon_A
\end{equation}
and similarly for $\sigma(BB,\ell_1)$. Comparing equation~\ref{1_multiplicity_AB_formula}
and equation~\ref{1_multiplicity_AA_formula} we deduce that either
there is an inequality
$$\sigma(AB,\ell_1) < \max(\sigma(AA,\ell_1),\sigma(BB,\ell_1))$$
or else $\epsilon_A = \epsilon_B = \epsilon'$ and therefore also $k_1^A = k_1^B$.

In other words, the shortest essential geodesic arcs in $A$ and in $B$ have the
same length and multiplicity, and are paired by the gluing map $(A,B) \to AB$. Since
the double of every such arc in $AA$ or $BB$ has complex length with zero
imaginary part, the holonomies of geodesic transport along
paired arcs in $A$ and in $B$ must be equal, or else there would be fewer arcs
in $AB$ with {\em complex} length $\ell_1$, and
$\sigma(AB) < \max(\sigma(AA),\sigma(BB))$ already. Let $\Gamma_1^A$ denote
the set of geodesic arcs in $A$ with length $\ell_1/2$ and define $\Gamma_1^B$
similarly. Then by the discussion above, without
loss of generality we can assume that the isometric identification of
$S \subset A$ and $S \subset B$ can be extended to isometric identifications
$$N_1^A = N_1^B$$
where $N_1^A$ is a regular neighborhood of $S \cup \Gamma_1^A$ and $N_1^B$ is defined similarly.

Note that the fact that the isometry can be extended over $S \cup \Gamma_1^*$ follows just by
considering the real lengths of the arcs in $\Gamma_1^*$. The fact that it can
be extended over a neighborhood depends on equality of {\em complex} lengths.
Note further that $\Gamma_1^A$ and $\Gamma_1^B$ might be empty.

\vskip 12pt

We will define inductively systems of proper geodesic arcs
$$\Gamma_1^A \subset \Gamma_2^A \subset \Gamma_3^A \subset \cdots$$
and regular neighborhoods
$$N_p^A = \text{regular neighborhood of } S \cup \Gamma_p^A$$
and similarly for $B$, where for each $p$, the set
$\Gamma_p^A \setminus \Gamma_{p-1}^A$ is the family of proper geodesic
arcs in $A$ orthogonal to $S$ having length exactly $\ell_p/2$ where
$\ell_p$ is the $p$th {\em distinct} term of $\sigma(AB)$ (i.e.\ {\em not} counted
with multiplicity).

We fix some notation which we use throughout the remainder of the proof.
Denote by $N_p^AN_p^B$ the corresponding subset of $AB$, and similarly for $AA$ and $BB$.
If $N_p^A$ is connected, let $[N_p^A]$ denote the covering space of $A$
induced by the image of $\pi_1(N_p^A)$ in $A$. 
If $N_p^A$ is disconnected, let $[N_p^A]$ denote the disjoint union of those
covering spaces.
The finite area boundary components of this disjoint union have an obvious
identification with $S$ so we continue to treat
$([N_p^A], S)$ as a pair.
Similarly,
$[N_p^AN_p^B]$ denotes the corresponding covering space of
$AB$. Note that in our notation there is an equality
$$[N_p^A][N_p^B] = [N_p^AN_p^B]$$
and similarly for $AA,BB$.

By induction, assume that $\sigma(AA)$, $\sigma(BB)$
and $\sigma(AB)$ are identical length spectra with multiplicities up to
$\ell_{p-1}$ and that the isometric identification $\partial A \to \partial B$
has been extended over $N_{p-1}^A \to N_{p-1}^B$. Since every perpendicular
geodesic arc in $A$ of length $< \ell_p/2$ can be doubled to a closed geodesic in $AA$
of length $<\ell_p$, by the induction hypothesis the length of such a geodesic
must be equal to one of the $\ell_q$ with $q<p$, and therefore $N_{p-1}^A$ contains every
perpendicular geodesic arc with length $< \ell_{p/2}$.

We write down a formula for the multiplicity $\sigma(AB,\ell_p)$ of
$\ell_p$ in $AB$ and similarly for $AA,BB$. Note that any geodesic $\gamma$ in $AB$
intersects the two sides in a collection of arcs. Each arc with length strictly
less than $\ell_p/2$ can be properly homotoped into the corresponding $N_{p-1}$
factor, as observed in the previous paragraph.

By the inclusion-exclusion formula,
\begin{equation}\label{p_multiplicity_AA_formula}
\sigma(AA,\ell_p) = 2k_p^A - k_{p-1}^{AA} + m_p^A - c_{p-1}^A
\end{equation}
where
\begin{align*}
k_p^A &= \# \lbrace \text{closed geodesics of length } \ell_p \text{ in } [N_{p-1}^A]A \rbrace \\
k_{p-1}^{AA} &= \# \lbrace \text{closed geodesics of length } \ell_p \text { in } [N_{p-1}^AN_{p-1}^A] \rbrace \\
m_p^A &= \# \lbrace \text{perpendicular geodesic arcs of length } \ell_p/2 \text{ in } A \rbrace \\
c_{p-1}^A &= \# \lbrace \text{such arcs properly homotopic into } N_{p-1}^A \rbrace
\end{align*}
and similarly for $\sigma(BB,\ell_p)$.
(We remind the reader that $[N_{p-1}^A]A$ denotes the union of $[N_{p-1}^A]$ and $\overline{A}$
along $S$, and similarly for other juxtapositions above.)

Note that $c_{p-1}^A$ is also equal to the number of doubled perpendicular arcs in $AA$
which lift to $[N_{p-1}^AN_{p-1}^A]$.
Note also that by induction, there are isometric identifications
$$N_{p-1}^AN_{p-1}^A = N_{p-1}^BN_{p-1}^B = N_{p-1}^AN_{p-1}^B$$
and therefore there are equalities
$$k^{AA}_{p-1} = k^{BB}_{p-1} = k^{AB}_{p-1}$$
where $k^{AB}_{p-1}$ is the number of closed geodesics of length $\ell_p$ in $[N_{p-1}^AN_{p-1}^B]$.

\vskip 12pt

We estimate $\sigma(AB,\ell_p)$.
Pairs of perpendicular geodesic arcs of length $\ell_p/2$ on different sides of $AB$
may not match up exactly. Moreover, even if they do match up, {\it a priori} different
arcs might have different holonomy, so the complex length of the resulting closed geodesic
might have nonzero imaginary part.

It follows that there is an {\em inequality}
\begin{equation}\label{p_multiplicity_AB_formula}
\sigma(AB,\ell_p) \le k_p^A + k_p^B - k_{p-1}^{AB} + \min(m_p^A - c_{p-1}^A,m_p^B - c_{p-1}^B)
\end{equation}
Denote $n_p^A:=m_p^A - c_{p-1}^A$.
It follows that either $\sigma(AB) < \max(\sigma(AA),\sigma(BB))$ and the theorem is proved,
or else (after possibly interchanging $A$ and $B$) there is an inequality
$$k_p^A + k_p^B + \min(n_p^A,n_p^B) \ge 2k_p^A + n_p^A \ge 2k_p^A + n_p^B$$
from which we can conclude that $k_p^A = k_p^B$ and $n_p^A = n_p^B$.

Furthermore, we deduce that the inequality in equation~\ref{p_multiplicity_AB_formula} is
actually an equality, and therefore perpendicular arcs of length $\ell_p/2$ in $A$ and $B$
match in pairs with the same holonomy.
This completes the inductive step, and lets us extend the
isometric identification $\partial A \to \partial B$ to $N_p^A \to N_p^B$.

\vskip 12pt

As we exhaust the length spectrum, we eventually obtain a complete generating
set for the fundamental group. That is, for $p$ sufficiently large,
$\pi_1(N_p^A) \to \pi_1(A)$ is an epimorphism.

It follows that the isometry on the boundaries extends to $A \to B$, and the lemma is proved.
\end{proof}

\subsection{Hyperbolic case with cusps}\label{noncompact_volume_section}

Theorem~\ref{AST_thm} and Lemma~\ref{spectrum_lemma} together let us define a complexity
function, namely the tuple $(-\vol,\sigma)$,
which is diagonally dominant for {\em closed} hyperbolic manifolds. However, we need
a generalization which is valid for {\em complete} hyperbolic manifolds of the kind which
arise in the JSJ decomposition. Lemma~\ref{spectrum_lemma} as stated applies directly
to such manifolds. The following is the required generalization of
Theorem~\ref{AST_thm} to the case of cusped manifolds.

\begin{theorem}\label{relative_AST}
Let $S$ be an orientable surface of finite type
so that each component has negative Euler characteristic,
and let $(A,S),(B,S) \in \Mdot(S)$ be irreducible,
atoroidal and acylindrical. Then $AA,AB,BB$ admit unique complete hyperbolic structures, and
either
$$2\vol(AB) > \vol(AA) + \vol(BB)$$
or else $\vol(AB) = \vol(AA) = \vol(BB)$ and $S$ is totally geodesic in $AB$.
\end{theorem}

Remark~\ref{geodesic_remark} applies equally well to the cusped case.
Note that by a hyperbolic Dehn filling argument, the only part of this theorem which does not
follow from Theorem~\ref{AST_thm} is the {\em strictness} of the inequality when $S$ is
not totally geodesic in $AB$.

\begin{proof}
Since $A,B$ are acylindrical, the manifolds $AA,BB$ and $AB$ admit
unique complete hyperbolic structures of
finite volume. The manifolds have two kinds of cusps: ``free cusps'' which do not intersect
cusps of $S$, and ``transverse cusps'' which do. The free cusps are irrelevant to the
discussion and for simplicity we assume they do not exist.

In fact, as a further simplification, we assume $AB$ has exactly one cusp $T$ which
intersects $S$ in two cusps with opposite orientations. It will be clear from the
proof in this case that multiple cusps present no additional difficulties.
Here we should think of $T = \partial \overline{AB}$ where $AB$ is homeomorphic to
the interior of $AB$, and $T_A = \partial_v \overline{A}, T_B = \partial_v \overline{B}$
are both annuli, each with two boundary components which compactify the two cusps of $S$.
The {\em meridian} on $T$ is the (unoriented) isotopy class which is isotopic to the cores of
the annuli $T_A$ and $T_B$

By \cite{MSY} the surface $S$ in $AB$ has a least area minimal representative in its proper isotopy
class.  Recall that a minimal surface is said to be {\em stable} if the second variation of
area is non-negative for all compactly supported variations. Least area surfaces are
stable. By abuse of notation, we call this surface $S$.
Since $S$ is stable, a fundamental estimate of Schoen applies.

\begin{theorem}[Schoen \cite{Schoen_estimate}]\label{stable_curvature_bound_theorem}
Let $S$ be a stable minimal surface in a Riemannian $3$-manifold $M$. Given
$r \in (0,1]$ and a point $p \in S$ such that the ball $B_r(p) \cap S$ has compact closure in
$S$, there is a constant $C$ depending only on the norms of $R_{ijk}^l$ and $\nabla^mR_{ijk}^l$
on $M$ restricted to $B_r(p)$ such that
$$|h_{ij}|^2(p) \le Cr^{-2}$$
where $h_{ij}$ is the second fundamental form.
\end{theorem}

This theorem generalized to the bounded case earlier curvature
estimates of Frankel \cite{Frankel}.

For two-sided surfaces, stability is preserved under covers (see
\cite{Colding_Minicozzi}, p.~21 for a proof).
In a complete hyperbolic $3$-manifold, the norms of
the curvature and its first derivatives are bounded by universal
constants. So Schoen's estimate implies a {\em uniform} pointwise
lower bound on the sectional curvature of a complete stable minimal
surface in a hyperbolic $3$-manifold.

\vskip 12pt

In the remainder of the proof, let $A,B$ denote the complete Riemannian manifolds with
boundary obtained from the hyperbolic $3$-manifold $AB$ by cutting along the stable minimal
surface $S$. Let $AA,BB$ denote the singular Riemannian manifolds obtained from $A$ and $B$
by metrically doubling along $S$. Our immediate goal is to prove short
time existence of the Ricci flow on the manifolds $AA,BB$, using the technology
developed by Simon \cite{Miles_Simon}.

In fact, as has become standard in discussion of Ricci flow,
following Simon \cite{Miles_Simon}, we use in place of Ricci flow a
generalization of the DeTurck gauging (DeTurck \cite{DeTurck}) called the
{\em dual Ricci-Harmonic Map flow} (see Hamilton \cite{Hamilton_singularities}, \S6
or Simon \cite{Miles_Simon}, p.~3 for a precise discussion), which
is equivalent to the Ricci flow up to a diffeomorphism. In what follows, we
refer simply to ``flow''.

In the sequel we suppress $BB$ and discuss only $AA$ for simplicity. As in
\cite{AST} and \cite{Miles_Simon} we must find suitable smooth approximations $AA^\delta$
with {\em uniform} pointwise lower bounds for scalar curvature {\em independent of $\delta$},
such that $AA^\delta \to AA$ as $\delta \to 0$. However, since the $AA^\delta$ are noncompact, we
actually employ a two-parameter family of smooth approximations $AA^\delta_k$ where
$\delta$ is a small positive real number and $k$ is one of an
infinite increasing sequence of
positive integers.
The $AA^\delta_k$ are
singular but compact approximations to $AA_k$ where
$$\lim_{k \to \infty} AA^\delta_k = AA^\delta \text{ for each }\delta, \;
\lim_{\delta \to 0} AA^\delta_k = AA_k \text{ for each }k$$
as Gromov-Hausdorff limits.
We assume the reader is familiar with Gromov-Hausdorff convergence and
Gromov-Hausdorff limits of (pointed) metric spaces; in the sequel we usually use the
term {\em Gromov limit} for brevity. For definitions and basic properties of
Gromov convergence, see \cite{Gromov_polynomial}, \S6.
The $AA_k$ turn out to be orbifolds obtained by
Dehn filling $k$ times the meridian of $AA$.

Let $AB_k$ be the closed hyperbolic orbifold
obtained by (orbifold) hyperbolic Dehn surgery on $AB$, by
filling $k$ times the meridian of $T$. The cusp $T$ becomes an
orbifold geodesic $\gamma_k \subset AB_k$ and the (topological) surface $S$ fills in
to become an orbifold $S_k$ transverse to $\gamma_k$ at two points.

Each surface $S_k$ has a least area minimal orbifold representative in its
isotopy class; indeed, following \cite{Freedman_Hass_Scott} and \cite{Hass_Scott},
by Selberg's Lemma (see e.g \cite{Raghunathan})
one may pass to a finite regular manifold cover $\til{AB}_k$ of $AB_k$, lift
$S_k$ to a topological surface $\til{S}_k$ in $\til{AB}_k$, and find a
least area minimal representative in the cover which is disjoint from or equal to
its image under every element of the deck group, or else its area could be
reduced by exchange and the roundoff trick. Then this least area representative covers
a least area orbifold representative in $AB_k$ which by abuse of notation we denote $S_k$.

Cutting $AB_k$ along $S_k$ produces $A_k,B_k$ and doubling these produces $AA_k,BB_k$.
Observe that since $\til{S}_k$ is fixed by elements of the deck group which
do not act freely, the geodesic $\gamma_k$ is perpendicular to $S_k$, and
$AA_k$ and $BB_k$ are singular orbifolds. That is, they have finite orbifold covers
which are locally isometric to the double of a hyperbolic manifold along a least area
minimal surface. Another way to see this is to cut $\til{AB}_k$ along $\til{S}_k$
to obtain $\til{A}_k$ and double, obtaining $\til{AA}_k$ which isometrically
orbifold covers $AA_k$.

\begin{lemma}
After passing to a subsequence of integers $k \to \infty$, and after possibly
replacing $S$ by another least area minimal surface, there is convergence
$$AB_k \to AB, \; S_k \to S, \; AA_k \to AA, \; BB_k \to BB$$
in the pointed Gromov-Hausdorff sense.
\end{lemma}
\begin{proof}
Thurston's hyperbolic Dehn surgery theorem (\cite{Thurston_notes} Chapter 4) implies
that for $k$ sufficiently large, $AB_k$ is hyperbolic, and $AB_k \to AB$ as $k \to \infty$.
By Schoen's Theorem~\ref{stable_curvature_bound_theorem} the surfaces
$S_k$ satisfy uniform two-sided curvature bounds. Since they have bounded geometry,
there is a convergent subsequence in $C^\infty$ whose limit is a minimal surface $S$.
The least area property is inherited by limits of minimal surfaces (\cite{Colding_Minicozzi} Chap.~1),
so $S$ is least area and may be taken to be $S$ as above.
\end{proof}

Next, we construct $AA_k^\delta$.  There are two alternate approaches.

\begin{enumerate}
\item{Bray (\cite{Bray}, \S6) uses the ODE for surface metrics
$\bar{g}_{ij}$ evolving normally in a 3-manifold
\begin{equation}
\label{bray_evolve_equation}
\frac{d}{dt}\bar{g}_{ij}(x,t) = 2\bar{g}_{ik}(x,t)h_j^k(x,t/\delta)
\end{equation}
where $-\delta \le t \le \delta$ and $h_j^k$ is antisymmetric in $t/\delta$,
to create a mirror-symmetric {\em plug} of thickness $2\delta$
interpolating between $A_k$ and its mirror image. The $h_j^k$ term is the second
fundamental form of the $t=\text{constant}$ ``slices'' of the plug.
This plug can be inserted between
the two copies of $A_k$ to build $AA_k^\delta$.}
\item{Miao (\cite{Miao}, \S3) simply mollifies the singular metric within a
tubular neighborhood of uniform thickness.}
\end{enumerate}

Since our surfaces $S_k$, $k \leq \infty$,
have no uniform size tubular neighborhoods to work with, we follow Bray \cite{Bray}
and desingularize by adding an untapered ``plug'' of thickness
$2\delta$ over $S_k$. Note that the constructions both of Bray and of Miao commute with
isometries, so that in practice we perform the desingularization in the cover
and then define $AA_k^\delta$ to be the quotient orbifold. Note for $\delta > 0$ that
$AA_k^\delta$ is a $C^\infty$ orbifold.

The metric and curvature uniformities of Bray's construction are
summarized in Lemma~\ref{metric_curvature_estimates}.
Our Lemma~\ref{metric_curvature_estimates} is parallel to \cite{AST} Prop.~4.1,
except that we include the parameter $k$ and note its uniformity.

\begin{lemma}\label{metric_curvature_estimates}
There exists a family of $C^\infty$ Riemannian
orbifolds $AA_k^\delta$, where $k \to \infty$, so that for each fixed $k$,
there is $\epsilon(k,\delta)$ so that $AA_k^\delta$ is $1+\epsilon$ bilipschitz
to $AA_k^0$, where $\epsilon \to 0$ as $\delta \to 0$. Moreover, these orbifolds
satisfy the following estimates:
\begin{enumerate}
\item{The scalar curvature satisfies $R(AA_k^\delta) \geq s(k)$, a (negative)
constant independent of $\delta$, for fixed $k < \infty$}
\item{The square norms of the full Riemann curvature tensors satisfy
an inequality $|R^l_{ijk}|^2 \le c(\delta)$. That is, they are uniformly bounded as a function of
position $x \in AA_k^\delta$ and $k \leq \infty$ (but not $\delta$).}
\end{enumerate}
\end{lemma}
\begin{proof}
The only new ingredient is the $k$-uniformity in the last assertion.
The second fundamental form $h_{ij}(x,t)$ enters into the estimates
of Bray (and Miao).  By Theorem~\ref{stable_curvature_bound_theorem} the $|h_{ij}|^2$ are
pointwise uniformly
bounded over all complete stable minimal surfaces in hyperbolic
$3$-manifolds.  So, the uniformity of $|h_{ij}|^2$ over position in
$S_k$ and over $k$ accounts for the new conclusions.
\end{proof}

Our next Lemma~\ref{conformal_modification} is parallel to \cite{AST} Prop.~4.2

\begin{lemma}\label{conformal_modification}
After a conformal modification, one can further assume that the metrics on $AA_k^\delta$
satisfy pointwise estimates for scalar curvature
$$R(AA_k^\delta) \ge -6$$
while still satisfying $AA_k^\delta \to AA_k$ in the bilipschitz sense.
\end{lemma}

Essential in Simon \cite{Miles_Simon} and Perelman \cite{Perelman_I, Perelman_II},
and implicit in the earlier work of Hamilton is the following principle:

\begin{lemma}\label{flow_commutes_with_Gromov}
Flow commutes with Gromov limit when flow on the
limit can be defined.
\end{lemma}

For Lemma~\ref{flow_commutes_with_Gromov} to be useful, it is necessary to
obtain lower bounds on the time to blow up which can be estimated uniformly.
Simon \cite{Miles_Simon} Thm.~5.2 says that for flow on a $C^0$ Riemannian
manifold $(M,g)$, the time $T$ to blow-up for flow can
be estimated in terms of $|\nabla^mR_{ijk}^l|$ of a {\em background metric}
$g'$ on $M$ which is $1+\epsilon$ bilipschitz close to $g$ for some
universal $\epsilon>0$ depending only on the dimension of $M$ (a precise statement is
the first part of Lemma~\ref{Simon_estimate} in this paper).
Simon's remarkable result, derived without
assuming that $M$ is compact or even that $g$ is better than $C^0$,
lets us prove a parallel to \cite{AST} \S6.1.

\begin{lemma}\label{flow_exists_for_T}
There is a uniform constant $T>0$ such that
for each $k$ and for all $\delta < \delta(k)$ flow exists
for time $t \in [0,T]$ on $AA_k^\delta$. Moreover, flow exists for
time $t \in [0,T]$ on $AA$.
\end{lemma}
\begin{proof}
Since, for each $k$ and each $\epsilon'>0$,
all the manifolds $AA_k^\delta$ and $AA_k$ are $1+\epsilon'$ bilipschitz when
$\delta$ is sufficiently small, it suffices to construct for each $k$
and for each $\epsilon > 0$, a $C^\infty$ metric $g(k)$ on $AA_k$ which is
$1+\epsilon$-bilipschitz to $AA_k$, and for which there are uniform
pointwise bounds on $|\nabla^mR_{ijk}^l|$, depending on $\epsilon$, but
{\em independent} of $k$.

For each individual
$k$, the existence of such a metric $g(k)$ is easy: any $C^0$ Riemannian
metric (i.e.\ a metric defined by a $C^0$ symmetric bilinear form on $TM$)
can be approximated (e.g.\ by mollifying in local co-ordinates) by some bilipschitz $C^\infty$
metric. The singular metric on $AA_k$ (away from the orbifold locus) is of
this kind, since it is obtained by doubling a genuine Riemannian metric
(also see Bray \cite{Bray}, equation~102 and the surrounding text for an
explicit estimate).
After we have obtained such an approximating metric, observe since $AA_k$
is compact, that there is some uniform pointwise bound on the
curvature and all its covariant derivatives.
However, the bound one gets in this case may well depend on $k$.

Using the fact that $AA_k \to AA$ in the Gromov sense, we see that this part
of the argument works {\em on the thick part} of $AA_k$. We need to find a
$(1+\epsilon)$-bilipschitz model for the thin part of each $AA_k$ with
square curvature bounds which depend on $\epsilon$ but not on $k$.
Then these two smooth $(1+\epsilon)$-bilipschitz models on the overlap
of the thick and thin parts can be welded together by a smooth convex
combination, at the cost of possibly increasing $\epsilon$ by a bounded
amount. So one just needs to choose the bilipschitz constants better than
necessary on each piece, so that the result of the welding is $1+\epsilon$
bilipschitz. It remains to find, for any $\epsilon>0$, a $(1+\epsilon)$-bilipshitz model for
the thin part of each $AA_k$ with uniform square curvature bounds
which depend on $\epsilon$ but not on $k$.

To do this, we must first prove a lemma about the geometry of
the cusped
least area minimal surface $S$ deep in the thin part.
We know that $S$ corresponds to a quasifuchsian group
since both $(A,S)$ and $(B,S)$ are acylindrical.

\begin{sublemma}\label{distinct_osculating}
Let $S$ be a quasifuchsian least area embedded surface in a complete cusped hyperbolic $3$-manifold $M$, and
let $\til{S}$ be a component of the preimage of $S$ in $\H^3$. Let $p \in S^2_\infty$ correspond
to a lift of the cusp. Let $B_t$ be a family of
horoballs centered at $p$ which are level sets of a Busemann function $b(t)$.
Then there is a totally geodesic plane $\pi$ through $p$ (the {\em osculating plane})
so that the restrictions 
$B_t\cap \til{S}$ and $B_t\cap\pi$ 
are Hausdorff distance $o(e^{-t})$ apart. Moreover, if $\til{S}_1,\til{S}_2$
are different components of the preimage of $S$ which both intersect $p$,
their osculating planes are distinct.
\end{sublemma}
\begin{proof}
In the upper half-space model (with $x,y$ as the horizontal co-ordinates and $z$ as the vertical co-ordinate),
put $p$ at infinity. The horoball $B_t$ corresponds to the set $z \ge e^t$, so we need to find
a vertical plane $\pi$ which is within Hausdorff distance $o(1)$ in the {\em Euclidean metric},
restricted to $z\ge t$; i.e.\ we need to show that
the Euclidean Hausdorff distance between the restrictions of
$\til{S}$ and of $\pi$ goes to $0$ as $z \to \infty$.

After composing with an isometry if necessary, we can assume that $\til{S}$ is
stabilized by $x \to x+1$. Let $\Lambda$ be the limit set of $\til{S}$. Then $\Lambda$ is the union
of $\infty$ with a proper quasiarc $\Lambda_0$ in the $(x,y)$-plane which is invariant under
$x \to x+1$. The arc $\Lambda_0$ is bounded in the slab $y \in [-C,C]$ for some constant $C$.
Since minimal surfaces are contained in the convex hulls of their boundaries, $\til{S}$ is also
contained in the slab $y \in [-C,C]$.

A stable minimal surface whose intersection with a compact set $K$ is trapped between two
barrier stable minimal surfaces which are $C^0$ close on $K$ is $C^\infty$ close
to both (in fact ``stable'' is superfluous here). Locally, this is just the Harnack inequality
for non-negative solutions of uniformly elliptic equations;
see e.g.\ \cite{Colding_Minicozzi}, pp. 20--21. This implies that $\til{S}$ is
$C^\infty$ close to the vertical planes $y=C$ and $y={-C}$ when $z$ sufficiently large.
Note that this is $C^\infty$ close in the {\em hyperbolic metric}. In the Euclidean
metric at height $z=t$, the order $n$ partial derivatives of $\til{S}$ and $y=C$ are
$o(t^{1-n})$ close; i.e.\ they differ by a term which is arbitrarily small compared
to $t^{1-n}$ as $t \to \infty$.
In particular, for $z$ sufficiently large, the tangent plane to $\til{S}$ at
each point is arbitrarily close to a vertical plane of the form $y=\text{constant}$, and
therefore $\til{S}$ is transverse to the level sets $z=\text{constant}$.
Notice that this does not yet tell us that the $y$ co-ordinate, thought of as
a function on $\til{S}$, converges to a constant
as $z \to \infty$; to establish this we must use the
{\em periodicity} of $\til{S}$, i.e.\ the fact that $\til{S}$ is invariant under
the parabolic translation $x \to x+1$.

Let $l_t$ be the intersection of $\til{S}$ with the horizontal plane $z=t$. Identifying
$z=t$ with the $(x,y)$-plane, we think of $l_t$ as the graph of a function $y=f_t(x)$
which satisfies $f_t(x+1)=f_t(x)$. Since $\til{S}$ is stable, the norm of its curvature
is bounded. Since as remarked above, the tangent plane to $\til{S}$ is uniformly close to
vertical when $z$ is sufficiently large, we can
estimate $|\partial^2f_t/\partial x^2| = O(1/t)$ (and $|\partial^2 f_t/\partial t^2| = O(1/t)$),
and so (by the periodicity of $f_t$) there is
an estimate $\max(f_t) - \min(f_t) = O(1/t)$.
But the minimal surface $\til{S} \cap B_t$ is trapped in the convex
hull of $\til{S} \cap \partial B_t = l_t \cup \infty$ for each $t$.
So $\max(f_t)$ is monotone decreasing as
$t \to \infty$ while $\min(f_t)$ is monotone increasing, and both have the same limit, which
can be taken to be $0$ after composing with an isometry. Hence setting $\pi$ equal to the
vertical plane $y=0$ satisfies the first claim of the sublemma.

\vskip 12pt

It remains to show that two different components of the preimage of $S$ have
distinct osculating planes.
Let $\til{S}_1,\til{S}_2$ be two components of the preimage of $S$, and let $\pi_1,\pi_2$ be
their osculating planes. Since the $\til{S}_i$ are disjoint, the $\pi_i$ can't cross, so
without loss of generality, we can set $\pi_1$ to be $y=0$ and $\pi_2$ to be $y=y_0$. We
want to show $y_0 \ne 0$. Let $f_{t,i}$ be as above, and suppose without loss of generality
that $f_{t,1}(x) < f_{t,2}(x)$ for each $t$ and each $x$. We want to show that $y_0>0$.

Suppose not, so that $y_0=0$ and $\pi_1 = \pi_2$.
Since both $f_{t,1}$ and $f_{t,2}$ are invariant under $x \to x+1$, there is a constant $C>0$
with $f_{t_0,2}- f_{t_0,1} > C$ for some fixed $t=t_0$. Let $\til{S}_2'$ be obtained from $\til{S}_2$
by translation $y \to y-C/2$, and let $g_{t,2} = f_{t,2} - C/2$ for each $t$. Then
$g_{t_0,2} - f_{t_0,1} > C/2>0$ but by the definition of the planes $\pi_1,\pi_2$ and the
hypothesis that they are equal, there is $t_1 > t_0$ with $g_{t_1,2} - f_{t_1,1} < 0$.
Note that both $\til{S}_2'$ and $\til{S}_1$ are invariant under $x \to x+1$, and we can
assume that $t_0,t_1$ are big enough so that the projection of $B_{t_0}$ is contained in
the Margulis tube of the cusp of $AB$. In other words, the ends of $\til{S}_1/(x \to x+1)$
and $\til{S}'_2/(x \to x+1)$ cross in an essential loop. By cut-and-paste and the roundoff trick,
we can reduce the area of $S$, contradicting the fact that $S$ was least area. This contradiction
shows $y_0 > 0$ as claimed.
\end{proof}

In our context, the surface $S$ indeed intersects the cusp in two components, giving rise
to $\til{S}_1,\til{S}_2$ as above, where the thin part of $A$ is covered by the slab
contained between $\til{S}_1$ and $\til{S}_2$ in $B_t$. Let $A_t$ denote this slab, and
let $\pi_t$ denote the slab contained between $\pi_1$ and $\pi_2$ in $B_t$. Since
there are estimates $|f_{t,i}'|= O(1/t)$, $|\partial f_{t,i}/\partial t| = O(1/t)$ for $i=1,2$,
the map $s:A_t \to \pi_t$ defined in co-ordinates by
$$s:\left( x,f_{z,1}(x)+\frac {t(y-f_{z,1}(x))}{f_{z,2}(x)-f_{z,1}(x)},z \right) \to (x,ty_0,z)$$
is $1+\epsilon$ bilipschitz, where $\epsilon \to 0$ as $t \to \infty$. In words,
$s$ fixes $x$ and $z$, and for each line $l$ which is parallel to the $y$-axis and intersects
$A_t$, it takes $l \cap A_t$ to $l \cap \pi_t$ linearly. The bilipschitz constant can
be estimated by the ratio $(f_{z,2} - f_{z,1})/y_0$ and by $1+$ the norm of the partial derivatives
of the graphs $f_1,f_2$ in the $x$ and $z$ directions.

Doubling, there is a $1+\epsilon$ bilipschitz map between the cusp of $AA$ and the hyperbolic
manifold (with horotorus boundary) obtained by doubling the quotient of the slab
$\pi_t/\langle x \to x+1 \rangle$. By Simon \cite{Miles_Simon} Thm.~5.2 flow exists for
a definite time $[0,T]$ on the singular manifold $AA$.

\vskip 12pt

The components of $\til{S}_k$ converge $C^\infty$ on compact subsets to components of $\til{S}$.
It follows that for big enough $k$, two components $T_1,T_2$ of $\til{S}_k$ which
intersect the same thin part have subdisks $D_i$ whose boundaries are
$C^\infty$ close to a pair of curves corresponding to components of
$\til{S} \cap \partial B_t$ where $t$ is
arbitrarily large (but fixed). In the notation of the previous paragraph,
we may assume $t$ is big enough so that
$|f_{t,1}| < \epsilon$, $|f_{t,2} - y_0| < \epsilon$ where $0 < \epsilon \ll y_0$.
In other words, there are round circles $C_1,C_2$ in $\H^3$ which are a constant distance
$\kappa=y_0/t$ apart so that each $\partial D_i$ is distance $\epsilon \kappa$
from each $C_i$. Let $E_i$ be a totally geodesic disk spanning each $C_i$, and let
$E_i^\pm$ be almost parallel totally geodesic disks spanning nearby round circles $C_i^\pm$
so that $D_i$ is trapped between $E_i^\pm$ for each $i$, by a barrier argument. We can choose $E_i^\pm$ so
that $E_1^+$ and $E_1^-$ are much closer than $E_1$ and $E_2$ (see Figure~\ref{trap_slab}).
The slab between the outer disks $E_1^+$ and $E_2^-$ can be foliated by totally geodesic
arcs $\gamma$ perpendicular to a bisecting plane such that for each $\gamma$,
and each choice of signs $\pm$, there is an inequality
$$\frac {d(D_1\cap \gamma,E_1^\pm\cap \gamma) + d(D_2 \cap \gamma,E_2^\pm\cap \gamma)}
{d(D_1 \cap \gamma, D_2 \cap \gamma)} \le 2\epsilon$$
(this follows just by comparing convex hulls).
Since the $D_i$ are $C^\infty$ close to the $E_i$,
each arc $\gamma$ intersects each $D_i$ and $E_i$ almost orthogonally,
with error term of order $O(\epsilon)$.
Stretching each arc $[\gamma \cap D_1,\gamma \cap D_2]$ linearly over
$[\gamma \cap E_1,\gamma \cap E_2]$ defines a
$1+O(\epsilon)$-bilipschitz map between the slab of hyperbolic space contained between
the $E_i$ and the slab contained between the $D_i$. The quotient of the second slab
by a rotation of order $k$ is half of the thin part of $AA_k$. Doubling the quotient
of the hyperbolic slab gives a suitable background metric on the thin part of $AA_k$,
so flow exists for time $[0,T]$ on each $AA_k$
and also on $AA_k^\delta$ for $\delta$ sufficiently small (depending on $k$).
\end{proof}

\begin{figure}[htpb]
\labellist
\small\hair 2pt
\pinlabel $E_1^+$ at 5 185
\pinlabel $D_1$ at -10 160
\pinlabel $E_1^-$ at 0 130
\pinlabel $E_2^+$ at 0 110
\pinlabel $D_2$ at -15 80
\pinlabel $E_2^-$ at 5 60
\endlabellist
\centering
\includegraphics[scale=0.6]{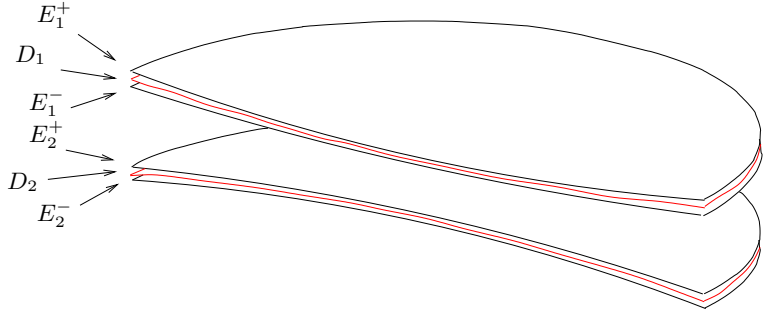}
\caption{The thin slabs between $E_1^+$ and $E_1^-$, and between
$E_2^+$ and $E_2^-$, trap $D_1$ and $D_2$, and show
that they cobound something close to the ``thick'' slab between $E_1$ and $E_2$}\label{trap_slab}
\end{figure}

Sublemma~\ref{distinct_osculating} motivates the following conjecture:

\begin{conjecture}
Let $S_1,S_2$ be complete, locally least area minimal surfaces in $\H^3$ which are either equal
or disjoint, and whose closures contain the same point $p \in S^2_\infty$.
If the Hausdorff distance of $S_1 \cap B_t$ and $S_2 \cap B_t$ is $o(e^{-t})$ then $S_1 = S_2$.
\end{conjecture}

Lemma~\ref{flow_commutes_with_Gromov} enhances Lemma~\ref{flow_exists_for_T} by
asserting that there is a diagonal sequence $\delta(k) \to 0$ so that the time $[0,T]$
flows on the manifolds $AA_k^{\delta(k)}$ converge pointwise and $C^\infty$ away from
time $0$ to the time $[0,T]$ flow on $AA$.

\vskip 12pt

We turn now to Hamilton's equation \cite{Hamilton} p.~698 for the evolution of
scalar curvature $R$ under Ricci flow ``with cosmological
constant''.  This is simply flow with a homothetic rescaling to
maintain constant volume;
denote it by $\mbox{flow}_0$.
The equation is
\begin{equation}\label{cosmological_flow}
\frac{dR}{dt} = \Delta R + 2|\Ric_0|^2 + \frac{2}{3}R(R-r)
\end{equation}
In the above, $\Ric_0$ is the traceless Ricci tensor, and $r$ is the
(spatial) average scalar curvature. We will use this equation on a
finite volume, but noncompact, manifold. Its original derivation,
for compact manifolds, still holds. Alternatively, note that we only
apply (\ref{cosmological_flow}) to compact manifolds and their
Gromov limits so the extension also follows from
Lemma~\ref{flow_commutes_with_Gromov}.

Parallel to \cite{AST}, to prove Theorem~\ref{relative_AST} it
suffices to show, using equation \ref{cosmological_flow} (or any
other method), that after short time $\mbox{flow}_0$ on $AA$, there
is an inequality $\inf(R) > -6$. Let us review this argument.
Rescaling the metric so that $\inf(R) = -6$ initially decreases
volume. Then, since Perelman's Ricci surgery is both volume
decreasing and infimum scalar curvature nondecreasing (see Perelman
\cite{Perelman_I,Perelman_II}), flow with surgery, normalized to
$\inf(R) = -6$, monotonically reduces volume. Thus
Theorem~\ref{relative_AST} is a consequence of the following:

\medskip\noindent{\bf Goal:} For sufficiently small $t>0$ show $\mbox{flow}_0 (t)(AA)$
satisfies $R > -6$.
\medskip

To obtain the goal we derive a kind of parabolic maximum principle
suitable to the cusped context.

As an important first step:

\begin{lemma}\label{non_decreasing}
For $t \in [0,T]$ $\mbox{flow}_0$ on $AA$ satisfies $R \geq -6$.
\end{lemma}
\begin{proof}
By the convergence $AA_k^\delta \to AA$ (here $\delta = \delta(k) \to 0$ as $k \to \infty$),
the inequality $R \geq -6$ for the initial
metric on each $AA_k^\delta$, and the existence and convergence of $\mbox{flow}_0$ for time
$[0,T]$ for each $k$ to flow for time $[0,T]$ on $AA$,
(i.e.\ Lemma~\ref{conformal_modification},
Lemma~\ref{flow_commutes_with_Gromov} and Lemma~\ref{flow_exists_for_T}) it
suffices to obtain this inequality termwise for each $k$, and take a limit.
From equation~\ref{cosmological_flow}, we see that at a negative spatial minimum
$\breve{R}:=\min(R)$ the right hand side is non-negative.  So, the parabolic maximum
principle implies that, on the compact manifold $AA_k^\delta$,
$\breve{R}$ is strictly increasing.  Taking limits we may lose strictness, but
obtain the weaker inequality.
\end{proof}

Assume that $AA$ is not hyperbolic (equivalently that $S$ is not totally geodesic,
equivalently that $AB$ is not obtained by a gluing isotopic to an isometry).
Under this assumption we have:

\begin{lemma}\label{identical_not_solution}
Under $\mbox{flow}_0$, $AA$ cannot satisfy $R(x,t) \equiv -6$
identically, at any finite time $t$.
\end{lemma}
\begin{proof}
Since $AA$ is not hyperbolic, it cannot become hyperbolic under flow
in any finite time.  Thus, if $R \equiv -6$, we must have $\Ric_0 \ne 0$
for some $(x,t)$ for all $t \in (0,T]$ and therefore, by analyticity, at
a set of full measure of points $x$ for each $t \in (0,T]$. At such points $x$,
equation~\ref{cosmological_flow} reads
$$0 = 0 + |\Ric_0|^2 + 0 > 0$$
which is a contradiction.
\end{proof}

We use the notation $AA(t)$ to denote the Riemannian manifold obtained by time $t$
flow on $AA$.

\begin{corollary}\label{average_increases}
$R(AA(t)) > -6$ for all $t \in (0,T]$.
\end{corollary}
\begin{proof}
The  inequality follows from $\inf(R) \ge -6$ (Lemma~\ref{non_decreasing})
and $R$ not identically $\equiv -6$ (Lemma~\ref{identical_not_solution}).
\end{proof}

Next we need a lower bound on $\Delta R$ at points $(x,t)$ where
$R(x,t)$ is near its spatial infimum which is $\ge -6$.
First, we state the principle for a single real variable.

\begin{lemma}\label{calculus_lemma}
Let $f: [-k,k] \rightarrow \mathbb{R}^+ \cup \{0\}$
be a non-negative $C^3$ function.  Let $c_3 = \max(f''')$ on
$[-k,k]$. Then
$$f(0) > -\frac{1}{4} f''(0)\left(\frac{f''(0)}{c_3}\right)^2$$
provided that
$k > 10\sqrt{\frac{f(0)}{f''(0)}}$.
Equivalently, if $$k > 10\sqrt{\frac{f(0)}{f''(0)}}$$ then
$$f''(0) > -(8c_3)^{\frac{2}{3}}f(0)^{\frac{1}{3}}$$
\end{lemma}

\begin{proof}
The worst case is when $f'(0) = 0$. In this case, $\frac{f''(0)}{2c_3}$ provides
a length scale over which $f$ decreases faster than
$-\frac{1}{2}f''(0)x^2$.
\end{proof}

The next Lemma is a multivariable version of Lemma~\ref{calculus_lemma}.

\begin{lemma}\label{covariant_calculus_lemma}
Let $f: M \rightarrow \mathbb{R}^+ \cup \{0\}$ be a
non-negative $C^3$ function on a complete Riemannian manifold with
(spatially) uniform bounds on curvatures $R_{ijk}^l$ and their first
three covariant derivatives.  Then there is a uniform constant
$c_3$, depending only on the preceding constants, so that at any
point $x \in M$, any second covariant derivative satisfies
$\nabla^2f(x) > -c_3f(x)^{\frac{1}{3}}$.
\end{lemma}

\begin{proof}
The proof is similar to that of Lemma~\ref{calculus_lemma} but in the
covariant context. One just needs to observe that the highest order terms
in $\nabla^2$ dominate.
\end{proof}

The uniform bounds required to apply Lemma~\ref{covariant_calculus_lemma} are
provided by Simon \cite{Miles_Simon} Thm~1.1, p. 1039, which in our context becomes:

\begin{lemma}[Simon \cite{Miles_Simon}, Thm~1.1]\label{Simon_estimate}
Let $M$ be a manifold.
Let $g_0$ be a complete $C^0$-metric and $h$ a
complete $C^\infty$ background metric on $M$, which is
$(1 + \epsilon)$-bilipschitz to $g_0$
(where $\epsilon$ is a universal constant depending on the dimension)
and $h$ satisfying uniform $\nabla^m R_{ijk}^l$ bounds
$|^h\nabla^m R_{ijk}^l| < k_m$.  Then there exists $T>0$, a
function of the $k$ and of dimension, and a $C^\infty$ family of
$C^\infty$ metrics $g(t)$ for $t \in (0,T]$ solving flow, such that $h$
remains $(1 + 2\epsilon)$-bilipschitz to the family and:
\begin{enumerate}
\item{$$\lim_{t \rightarrow 0} \sup_{\cdot \in M}|g(\cdot,t)-g_0(\cdot)| = 0$$}
\item{$$\sup_{x \in M}|^h\nabla^i g|^2 \leq \frac{c_i(\dim,k_0,\dots,k_i)}{t^i}$$}
\end{enumerate}
where $^h\nabla$ denotes covariant derivative in the $h$ metric.
\end{lemma}

We will use $R+6$ on $AA$ with the background metric $h$
as the non-negative function in Lemma~\ref{covariant_calculus_lemma}.
First spatial derivatives $\nabla(\Delta R)$ of
the Laplacian are estimated by third covariants $|\nabla \Delta R| = O|\nabla^3 R|$
and by Lemma~\ref{flow_exists_for_T} and Lemma~\ref{Simon_estimate} are uniformly bounded
in $x$ for any fixed interval $[t_1, t_2] \subset (0,T]$.  Now
Lemma~\ref{covariant_calculus_lemma} yields:

\begin{lemma}\label{bound_Delta_R}
For any fixed time interval $[t_1, t_2] \subset (0,T)$ and any $t \in [t_1,t_2]$
there is a lower bound for $\Delta R$ on $AA(t)$ ; precisely,
$$\Delta R(x,t) > -c_3(R(x,t)+6)^{\frac{1}{3}}$$
where $c_3$ is a positive constant independent of $x$ or $t \in [t_1, t_2]$.
\qed
\end{lemma}

Now consider Hamilton's equation~\ref{cosmological_flow} at initial time $t_1$ for
$\mbox{flow}_0$ on
$AA(t_1)$.  We set $\breve{r} = \min_{t \in [t_1, t_2]} r(t)$.
By Corollary~\ref{average_increases} and compactness of $[t_1, t_2]$,
there is a strict inequality $\breve{r} > -6$.

\begin{claim}
There is an $r_0 \in (-6, \breve{r}]$ with the property that at all
$(x,t)$ with $t \in [t_1, t_2]$ such that $R(x,t) < r_0$, there is
an inequality
\begin{equation}\label{half_cancel_equation}
\Delta R(x,t) + \frac{2}{3}R(x,t)(R(x,t) - \breve{r}) > \frac{1}{3}R(x,t)(R(x,t)-\breve{r})
\end{equation}
\end{claim}
In words, the claim says that the negative contribution of $\Delta
R$ can cancel at most half of the positive contribution of the final
term in Hamilton's equation~\ref{cosmological_flow} Assuming this
claim for the moment, we complete the proof.

From the claim, throughout $[t_1, t_2]$, there is an inequality
$dR(x,t)/dt \ge u$ where $u: = \frac{r_0}{3} (r_0 -
\breve{r})$ at those points $(x,t)$ with $R(x,t) \in [-6, r_0]$. That
is, at such points, $R(x,t)$ is increasing at a definite rate. Thus,
after $\mbox{flow}_0$ for time $t_2 - t_1$, starting at time $t_1$ we
conclude:
\begin{equation}
\label{increase_rate} R(x,t_2) \geq \min (r_0,(t_2 -t_1)u -6) > -6
\end{equation}

Equation~\ref{increase_rate} implies that the (unnormalized) flow on
$AA$ instantly reduces volume:
$$\vol(AA(t)) < \vol(AA) \text{ for all } t \in (0,T]$$

The proof of Theorem~\ref{relative_AST} in the cusped case now follows exactly as
in \cite{AST}. We now give the proof of the claim.

\begin{proof}
Here is how to construct $r_0$
so that equation~\ref{half_cancel_equation} holds
when $R(x,t) < r_0$.

We assume that we have already chosen the constant, and derive
(easily satisfied) conditions that it must satisfy. So assume that
$R(x,t) < r_0$. We have $\Delta R > -c_3(r_0 + 6)^{\frac{1}{3}}$.
Set $\epsilon' = \breve{r} + 6$ and $\epsilon = r_0 + 6$. Ultimately
we will choose $\epsilon$ (much smaller than $\epsilon'$) and
thereby choose $r_0$.

The second term of equation~\ref{half_cancel_equation} exceeds
$\frac{2r_0}{3}(r_0 - \breve{r})$.
$$\frac{r_0}{3} (r_0 - \breve{r}) > c_3 (r_0 + 6)^{\frac{1}{3}}$$
that is,
$$(\epsilon - 6)(\epsilon - \epsilon') > 3c_3 \epsilon^{\frac{1}{3}}$$
Since we want $\epsilon$ very small, $\epsilon - 6$ is negative and
bounded away from zero. So by replacing $c_3$ with a nearly
identical constant $c_3'$, we have $(\epsilon' - \epsilon) > c_3'
\epsilon^{\frac{1}{3}}$ which can be rearranged as $c_3'
\epsilon^{\frac{1}{3}} + \epsilon < \epsilon'$. Obviously we can
find such an $\epsilon'>0$ with this property, and then set $r_0 =
\epsilon - 6$.
\end{proof}

This completes the proof of Theorem~\ref{relative_AST}.
\end{proof}

\subsection{Hyperbolic complexity}

We are therefore justified in defining the hyperbolic complexity
$c_h$ as follows.

\begin{definition}
Let $M$ be a connected, complete, finite volume, orientable
hyperbolic $3$-manifold. Define $c_{ch}(M) = (-\vol(M),\sigma(M))$
as a lexicographic tuple.
\end{definition}

\begin{definition}
Let $M$ be a complete, finite volume, hyperbolic $3$-manifold.
Define $c_h(M) = \{c_{ch}(M_i)\}$ as a lexicographic tuple, where
$M_i$ are the components of $M$.  We
also use $c_h$ to denote the ``hyperbolic complexity'' of
a connected closed
irreducible $3$-manifold by applying $c_h$ to the hyperbolic
JSJ pieces, listed in decreasing order.
\end{definition}

\begin{theorem}[$c_h$-Lemma Schema]\label{hyperbolic_lemma_schema}
Let $S$ be an orientable surface of finite type with no sphere or
torus components, and let $(A,S)$, $(B,S)$ be distinct elements of
$\Mdot(S)$ which are irreducible, atoroidal and acylindrical, and such that
every component has some part of their boundary on $S$. Then
$$c_h(AB) < \max(c_h(AA), c_h(BB))$$
\end{theorem}
\begin{proof}
This follows from Theorem~\ref{relative_AST} and Lemma~\ref{spectrum_lemma}.
\end{proof}

\section{Assembly complexity}\label{c_a_section}

\subsection{Introduction}\label{introduction_subsection}

Let us pause to review where we are in the proof.

We seek a complexity function $c$, defined on homeomorphism classes
of closed 3-manifolds, such that if $A$ and $B$ have common boundary
$S$, then $c(AB) \le \max(c(AA), c(BB))$, with equality holding only
if $A$ and $B$ are homeomorphic rel $S$. We define $c$ as a
lexicographic tuple, with terms introduced throughout the course of
this paper. The input for this section will be the complexities
$c_S$ and $c_h$ of Sections $4$ and $5$.  The output of this section
will be a complexity, $c_p$, defined on prime (actually irreducible)
3-manifolds. We assume throughout this section that $AB$, $AA$ and $BB$
are connected and irreducible. The surface $S$, which need not be
connected, is further assumed to be incompressible in both $A$ and $B$.

Assuming $c_S(AB) = c_S(AA) = c_S(BB)$ and $c_h(AB) = c_h(AA) =
c_h(BB)$, we then know that these three closed 3-manifolds become
homeomorphic after cutting along JSJ tori. Lemma~\ref{SF_DD_lemma}
and Theorem~\ref{hyperbolic_lemma_schema} allow us to draw a similar
conclusion about the non-closed 3-manifolds $A$ and $B$. We recall
the framework of these arguments.

\subsubsection{Notation for pieces}

Consider the ways in which $S$ might sit with respect to the JSJ
decomposition of the glued up manifold which we refer to as $AB$
(with the same terminology applying to $AA$ and $BB$).
We distinguish between two different kinds of torus components of $S$:
\begin{definition}\label{vertical_horizontal_cusp_definition}
A {\em JSJS cusp} is a component $T$ of $S$ which is a JSJ torus.
An {\em internal cusp} is a component $T$ of $S$ which is contained in
the interior of a JSJ piece in $AB$, built from Seifert fibered
pieces in $A$ and $B$ containing $T$ in their boundary, and admitting compatible
fiberings of $T$.
\end{definition}
Internal cusps contribute to $m$ (= \# of independent tori)
but not $m'$ (= \# of JSJ tori) in the complexity $c_S$ defined
in Definition~\ref{c_SF_defn}.
In the proof of Lemma~\ref{SF_DD_lemma} it is established that
each JSJ torus is either a JSJS cusp, or (after isotopy)
can be chosen to intersect $S$ in a pair of essential circles,
which cut the torus into two essential annuli which are proper in $A$ and
$B$ respectively.

Thus the relative JSJ decompositions of $A$ and $B$ are compatible
with the JSJ decomposition of $AB$. Say that a JSJ piece of $A$ is a
{\em boundary JSJ piece}, or {\em boundary piece} for short,
if it has some boundary component on $S$ which
is not a JSJS cusp. In other words, boundary pieces
correspond to JSJ pieces of $AA$ (or $AB$ or $BB$) that are cut by $S$.

\subsubsection{Assembly complexity}

What, then, is left to do? We need to define a further complexity term $c_a$
such that $AB$, $AA$ and $BB$ as above satisfy
$$c_a(AB) \le \max(c_a(AA),c_a(BB))$$
with equality if and only if $A$ and $B$ are diffeomorphic rel $S$
(in other words, if and only if $(A,S) = (B,S)$). The term $c_a$
is sensitive to the way in which the JSJ pieces are {\em assembled}.

The way in which a $3$-manifold is assembled from its JSJ pieces can
be encoded by a decorated graph. The vertices of the graph are
labeled by (a representative of) the homeomorphism type of a JSJ
piece. The edges are labeled by the gluing homeomorphisms of the
cusps. Similarly, a relative JSJ decomposition can be encoded by a
decorated {\em relative} graph, similar in many ways to the relative sum
graphs defined in \S\ref{c_2_section}. Some of the vertices/edges in a
relative graph are special:
\begin{itemize}
\item{Some of the vertices (those corresponding
to boundary pieces) are {\em half vertices}.}
\item{Some of the edges
(those corresponding to proper JSJ annuli) are {\em half edges}.}
\item There are {\em cut points} corresponding to tori components of $S$.  We do
not consider the cut points to be vertices.
\item{There are also {\em cut edges} which connect a cut point to an ordinary vertex
(or half vertex or another cut point),
which correspond to the JSJS cusps and the JSJ pieces (on one side) which they
bound.}
\end{itemize}
In the graphical context, two half vertices glue together to make a vertex.
Two half edges are glued {\em lengthways} to make an edge, and two cut edges
are joined {\em at their cut points} to make an edge. We also allow the
case of a {\em doubly cut edge}, which is an isolated interval component
of the graph, with two cut points and no vertices.
These of course correspond to $T^2\times I$.
Compare with
the terminology in Definition~\ref{relative_sum_graph_definition},
where half vertices corresponded to components which
intersected $S$ and half edges corresponded to disks.

There is a symmetry group associated
with each vertex and half vertex, and gluing information
associated to each edge or half edge. By uniqueness of the JSJ
decomposition, two $3$-manifolds are homeomorphic if and only if there
is an isomorphism between their JSJ graphs which preserves vertex
labels and is compatible with the gluing data, up to the
action of the symmetry groups.

\subsubsection{Graph tensor TQFTs}

As a warm-up problem, we establish positivity for a few ``toy''
unitary TQFTs associated to labeled graphs with internal partial
symmetries. In the simplest case, graphs are glued only along cut
points; this is analogous to the case in which $S$ consists entirely
of JSJS cusp components. In the second case, graphs are glued along
subgraphs, corresponding to boundary pieces of $A$ and $B$.

The strategy of the proof is to build a tensor TQFT for labeled graphs
with certain desirable properties.
Chief among these are positivity (often called ``unitarity'' in the TQFT context)
and genericity.
\begin{definition}
A {\em graph tensor TQFT} consists of the following data:
\begin{itemize}
\item For every edge type $e$, a (real or complex)
vector space $V^e$ with a positive nondegenerate inner product
(symmetric in the real case, Hermitian in the complex case)
\item For every vertex type $v$ with incident edge types indexed by a (finite) set
$E(v)$, a tensor $T^v$ in $\bigotimes_{e \in E(v)} V^e$
\end{itemize}
A labeled graph tensor TQFT is {\em symmetric} if there is only one edge type,
and each vertex tensor is chosen to lie in the symmetric tensor power $T^v \in S^{n(v)} V$
where $n(v)$ is the number of edges incident on $v$.
\end{definition}

A cut point does not get a tensor.
A doubly cut edge gets the identity element of $\Hom(V^e,V^e)$,
where $e$ is the edge label.
Using the duality this can be denoted as $v_i\otimes v^i$, $v_i\otimes v_i$,  or
$v^i\otimes v^i$.

\subsubsection{Contracting tensors}

The positive nondegenerate inner product identifies $V^e$ with $(V^e)^*$ when $V^e$ is
real, or with $(\overline{V^e})^*$ when $V^e$ is complex. Thus the indices
of $T^v$ may be raised or lowered as desired. So if $v_1,v_2$ are
a pair of vertices which share an edge of type $e$, we may contract
the tensors $T^{v_1}$ and $T^{v_2}$ along their respective indices which correspond
to the factor of type $V^e$. Technically, this amounts to forming $T^{v_1} \otimes T^{v_2}$,
identifying an element of $V^e \otimes V^e$ inside this tensor, using the inner
product to think of this as an element of $V^e \otimes (V^e)^*$, and then taking a
trace to replace it with a scalar. By contracting along edges in this way, a labeled
graph with ``free edges'' $e_1,\cdots,e_n$ gives rise to a tensor in
$\bigotimes_i V^{e_i}$, and a labeled graph with no free edges gives rise to
an element of $\R$ or $\C$. Pairs of graphs with isomorphic sets of free edges
can be glued together, and the result is a graph whose ``value'' is obtained by
appropriately contracting the two tensors associated to the subgraphs.
(Note: we use the terminology {\em cut edges} in \S\ref{acgcssect} instead
of {\em free edges}, but the meaning should be clear in each case).

Typographically, this can be achieved by choosing an orthonormal basis for each
vector space $V^e$, and then expressing the components of the tensors $T^v$ in terms
of this basis, as (lower) indices, taking values in the set of basis elements.
Contracting is indicated by repeating an
index, which in the Einstein summation convention means ``sum over all possible
values of this index''. So for example $T_{ij}T_{jk}$ means $\sum_{j=e_1}^{e_n} T_{ij}T_{jk}$.

For the sake of clarity, we give a few examples. As a simplification, assume that
we are working with a symmetric labeled graph tensor TQFT. Again for
simplicity, assume that all vertices are $3$-valent and have the same type; denote the
vertex tensor by $x$. Let
$e_1,\cdots,e_n$ be an orthonormal basis for $V$. In a symmetric tensor,
the position of the indices in a given tensor are immaterial, so that for instance
$x_{ijk} = x_{kij} = x_{jki}$.

\begin{example}
\begin{figure}[htpb]
\labellist
\small\hair 2pt
\pinlabel $\Gamma_3$ at 550 130
\pinlabel $\Gamma_2$ at 255 160
\pinlabel $\Gamma_1$ at 42 160
\endlabellist
\centering
\includegraphics[scale=0.4]{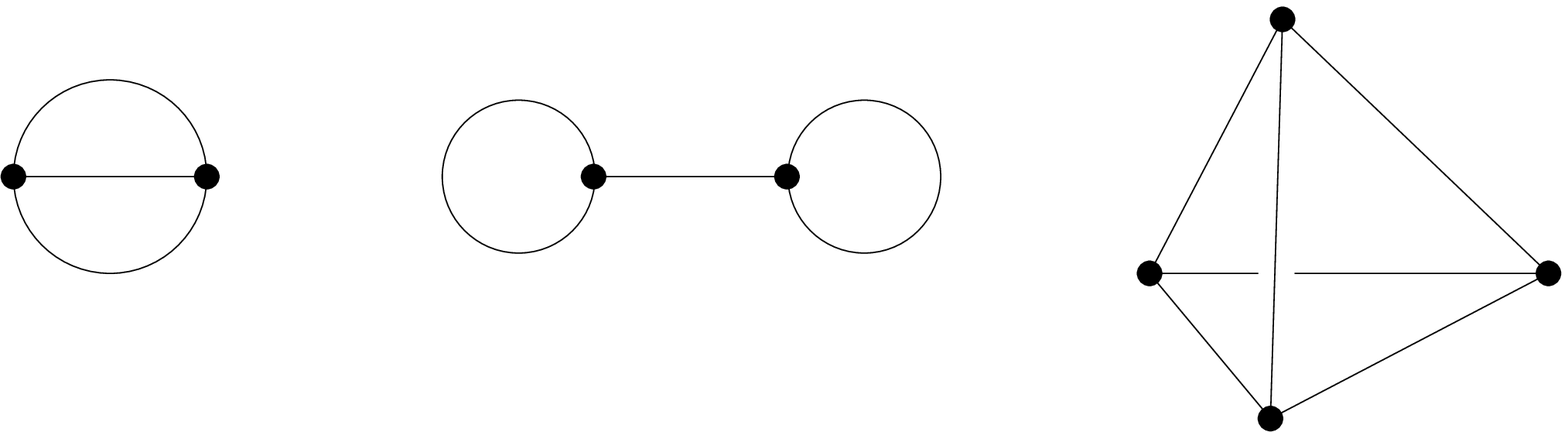}
\caption{} \label{fourgraphs}
\end{figure}
Consider the three graphs in Figure~\ref{fourgraphs}.
The theta graph $\Gamma_1$ has invariant $x_{ijk}x_{ijk}$. The barbell graph $\Gamma_2$
has invariant $x_{iij}x_{kkj}$. The tetrahedron graph $\Gamma_3$ has invariant
$x_{ijk}x_{ilm}x_{jln}x_{kmn}$.
\end{example}

Any choice of vector spaces and tensors as above defines a TQFT,
and a generic choice of tensors defines a TQFT which is powerful
enough to prove a version of positivity. If vertices are decorated
with symmetry groups which act on their incident edges, one must choose
tensors with the same symmetries. Furthermore, edges might have ``internal''
symmetries, reflected in the TQFT by choosing $V^e$ which are not just vector spaces
but $G$-modules for various groups $G$, and in this case
the tensors associated to vertices must be chosen in such a way as to respect
this $G$-action. One way to do this is to replace $V^e$ by some tensor power (on which $G$ also
acts) and choose the various tensor factors in $T^v$ to lie in the $G$-invariant subspaces.

Labeled graph tensor TQFTs of progressively greater complexity
are constructed in \S\ref{acgcssect}. The contents
of this section are not logically necessary for the rest of the paper, but
they will aid the reader considerably in following the logic in
\S\ref{actjsjsect} and subsequent sections.
Also, we have not attempted to be comprehensive in this motivational subsection;
a paper explaining the full scope of ``graph positivity'' should still be written.

\vskip 12pt

The construction of a suitable complexity function for JSJ graphs is
modeled on the template of proof developed in \S\ref{acgcssect}. It
has the following outline.
First, choose a vector space for each JSJ torus, and choose compatible
tensors $T_X$ associated to each JSJ piece $X$. The JSJ tori have ``internal'' symmetries
(parameterized by copies of $\GL(2,\Z) = \pm\SL(2, \Z)$) and these symmetries interact with
the symmetries of the JSJ pieces they bound. The major difficulty
is to find suitably symmetric $T_X$. If $T_X$ has too little symmetry,
the partition function will not be well-defined; too much and it will not distinguish
distinct glued-up manifolds. When $X$ is hyperbolic, tensors $T_X$ with
the correct (finite) symmetries can be constructed directly. Seifert
fibered pieces $X$ present more of a problem, since their symmetries
are of infinite order. So first we construct tensors $\widehat{T}_X$
for Seifert fibered $X$ with {\em too much} symmetry, which distinguishes
$AA$, $AB$ and $BB$ only up to an equivalence relation slightly weaker than homeomorphism,
which we call {\em fiber slip homeomorphism}, or {\em fish} for short. Then
we adjust $\widehat{T}_X$, replacing it by a new tensor $T_X$
which is sensitive not just to the topology of $X$, but the way in which it is glued
to its immediate neighbors in the JSJ decomposition, to control for fish equivalence.
There is one further technical point which should
be mentioned here.  When we come to define the assembly complexity $c_a$ it will be a pair
$c_a = (c_r, c_t)$. The first term, a ``reflection symmetry'' term
at the JSJ tori, regularizes the picture near $S$ and ensures the
conditions under which the tensor contraction term $c_t$ satisfies diagonal dominance.

\begin{remark}
The reader may wonder why the issue of assembly complexity cannot be
treated by simply adding a term to $c$ which ``likes'' lots of symmetry
(e.g.\ by generalizing the term $c_r$ somehow).
This is a very seductive idea as it appears to offer a
rapid finish.  Unfortunately, one may produce examples where $AB$
has ``more'' symmetry than either $AA$ or $BB$ for all notions of
``symmetry'' which we were able to quantify.
\end{remark}

\subsection{Graph positivity} \label{acgcssect}

As a warm-up, we now prove two versions of positivity for labeled
graphs. The vertices of the graph are analogous to JSJ pieces, the
edges are analogous to cusps, and a vertex label is analogous to the
homeomorphism type of a JSJ piece. Accordingly, we require that the
vertex label determine the valence of a vertex. We assume the set of
possible vertex labels is countable and that all graphs are finite.

\subsubsection{Recovering graphs from tensors}

In order for our invariants to say anything about the topology and
combinatorics of graphs, it is essential that the invariants are at
least powerful enough to {\em distinguish} graphs. We prove a lemma
to this effect initially in the context of edge-symmetric, labeled
graph tensor TQFTs, where vertices are distinguished only by the
number of incident edges, so there is exactly one vertex type for
each non-negative integer. We call this the context of ``ordinary''
graphs.

\begin{lemma} \label{enough_dimensions_distinguishes}
Let $V$ be a real vector space of dimension $m$. Let $\Gamma_1$,
$\Gamma_2$ be two (ordinary) graphs with at most $m$ edges. Assume
an identification between the cut points (boundary) of $\Gamma_1$
and $\Gamma_2$. Then for a generic assignment of tensors $T^v \in
S^{n(v)}V$ to vertex types, the invariants $T(\Gamma_1)$ and
$T(\Gamma_2)$ are not equal unless $\Gamma_1$ and $\Gamma_2$ are
isomorphic rel boundary.
\end{lemma}

\begin{proof}
Fix an orthonormal basis of $V$. For each vertex type $v$, and each
(unordered) set (with multiplicity) of indices for $T^v$, choose a
variable representing the value of the given tensor entry. We
distinguish these indeterminates by letters $x,y,z$ depending on the
cardinality of $n(v)$, and the indices by subscripts. So for
example, $x_{ijk}$ and $y_{abcd}$ are examples of these variables.
Then the result of tensor contraction can be expressed as a {\em
homogeneous polynomial} in these variables, with degree equal to the
number of vertices. In Einstein summation notation, the value of
$T(\Gamma_1)$ is represented by a monomial whose subscripts are
indeterminate indices which appear singly or in pairs: one pair for
each (uncut) edge of $\Gamma_1$ and a singleton for each cut edge of
$\Gamma_1$. Since by hypothesis the number of edges is at most equal
to the dimension of $V$, there is some monomial in the polynomial
representing $T(\Gamma_1)$ in which each basis element appears as a
subscript at most twice. By taking a vertex for each variable which
appears in the monomial, an edge for each pair of variables with a
common index, and a cut edge for each singleton index, we can
reconstruct $\Gamma_1$ from the polynomial, and similarly for
$\Gamma_2$. This monomial contains the ``gluing instructions'' for
the graph: if $b$, say, occurs as a subscript for both $x$ and $y$,
then some edge must join the vertices corresponding to $x$ and $y$.
It follows that if the entries of the various $T^v$ are chosen
algebraically independently, the values of the invariants
$T(\Gamma_1)$ and $T(\Gamma_2)$ are different unless $\Gamma_1$ and
$\Gamma_2$ are isomorphic relative to the identity on cut edges.
\end{proof}

To put ourselves in the favorable case where the number of edges of
a given graph is less than the dimension of $V$, we may let
$V_1,V_2,V_3,\cdots$ be a sequence of vector spaces with inner
products, with $\dim(V_i) = i$ for each $i$. Construct a symmetric
labeled graph tensor TQFT $T_i$ for each $i$, with algebraically
generic entries. Then for any two finite graphs $\Gamma_1,\Gamma_2$
the strings of invariants $T_1(\Gamma_1),T_2(\Gamma_1),\cdots$ and
$T_1(\Gamma_2),T_2(\Gamma_2),\cdots$ are equal {\em as strings} if
and only if $\Gamma_1$ and $\Gamma_2$ are isomorphic {\em as graphs}
(taking the cut point boundary of the graph into account, of
course). The positivity theorem below is a consequence of being able
to distinguish graphs via a family of tensor TQFTs.

Let $C$ be a finite ordered set of cut points. Let $\dot\cN_C$
denote the set of isomorphism classes of ordinary graphs with
boundary (cut points of cut edges of the graph) identified with $C$.
Let $\cN_C$ denote the vector space spanned by $\dot\cN_C$. Gluing
along cut points yields a pairing into closed graphs:
\[
    \cN_C \otimes \cN_C \to \cN .
\]

\begin{theorem}\label{ogpos}
The above pairing (for ordinary graphs) is positive.
\end{theorem}

\begin{proof}
This theorem is an easy special case of Theorem \ref{grpos2b},
so we defer the proof until then.
\end{proof}

For a second warm-up example, let us consider a class of asymmetric
graphs. These are graphs with: 1.~one edge type, 2.~arbitrarily many
vertex types, and 3.~no vertex symmetry (i.e., the edges leaving a
vertex are ordered). Two graphs are considered isomorphic if and
only if there is a isomorphism between the underlying unlabeled
graphs which preserves vertex labels and also the ordering of the
edges at each vertex. If the graphs have boundary, then the
isomorphism must preserve the ordering of the cut points. Let $\cN$
denote the real vector space generated by finite linear combinations
of isomorphism classes of such graphs.

For $C$ a finite ordered set of cut points, let $\dot\cN_C$ denote
the set of isomorphism classes of labeled graphs whose ``boundary''
is $C$. Two such graphs are considered isomorphic if there is an
isomorphism as above preserving labels and orderings, and which is
the identity on $C$. Let $\cN_C$ denote the vector space generated
by $\dot\cN_C$. Gluing along cut points gives a pairing
\[
    \cN_C \otimes \cN_C \to \cN .
\]

\begin{theorem}
\label{grpos1b} The above pairing, for asymmetric vertices and
gluing along cut points, is positive.
\end{theorem}

\begin{proof}
We will construct below a family of unitary TQFTs $\{Z_i\}$, $i = 1,
2, 3, \ldots$, for graphs which eventually distinguish
non-isomorphic (relative) graphs in $\dot\cN_C$. If $A$ is a closed
graph, define $Z(A)$ to be the infinite tuple $(Z_i(A))$, with
lexicographic ordering. We will show that $Z(\cdot)$ is a diagonally
dominant complexity function. If $A, B \in \dot\cN_C$ are
graphs with boundary $C$, then using only the Atiyah
gluing axiom (which holds trivially for graph tensor TQFTs) and the
Cauchy-Schwarz inequality we obtain (for all $i$)
\begin{eqnarray*}
    |Z_i(AB)| &=& |\langle Z_i(A), Z_i(B)\rangle| \\
        & \le & \left(|\langle Z_i(A), Z_i(A)\rangle|\cdot |\langle Z_i(B), Z_i(B)\rangle|\right)^{\frac{1}{2}} \\
        & \le & \max(Z_i(AA), Z_i(BB)) ,
\end{eqnarray*}
with equality holding only if $Z_i(A) = Z_i(B)$. If $A$ and $B$
are not isomorphic, 
we will show below that there
exists an $i$ such that $Z_i(A) \ne Z_i(B)$. It follows that $Z(AB)
< \max(Z(AA), Z(BB))$. In other words, $Z$ satisfies the diagonal
dominance inequality.

It remains to construct the TQFTs $Z_i$. Choose a real inner product
space $V = V_i$ of dimension $i$. Fix an orthonormal basis for $V$.
The vector space $V(C)$ assigned to $C$ will be $V^{\otimes C}$. For
each vertex label $\alpha$ of valence $k$, choose a tensor $T^\alpha
\in V^{\otimes k}$. With respect to the basis of $V$, this means
choosing numbers $T^\alpha_{i_1\ldots i_k} \in \R$. The key idea is
to choose these numbers (for all vertex labels and all
multi-indices) to be algebraically independent of one another. For a
labeled graph $G$ (possibly with a boundary consisting of cut
points) define $Z(G) \in V(\partial G)$ to be the tensorial
contraction described above. (See also for example Penrose
\cite{pppp}.)

Suppose now that $A$ and $B$ are graphs with boundary $C$, and that
$Z(A) = Z(B)$. An easy variation on the proof of Lemma
\ref{enough_dimensions_distinguishes} shows that if 
$Z_i(A) = Z_i(B)$ where
$i = \dim V$ is
greater that the number of edges in $A$ and $B$, then $A$ and $B$
must be isomorphic.
\end{proof}

\medskip

Next we consider the case where each vertex type (label) has a
specified symmetry group --- a subgroup of the permutation group of
the edges emanating from the vertex. These symmetry groups are
analogous to symmetries of JSJ pieces. Two graphs are considered
isomorphic if and only if there is a isomorphism between the
underlying unlabeled graphs which preserves vertex labels and such
that the induced permutation of the (ordered) edges at each vertex
lies in the symmetry group for that vertex type. If the graph has
boundary, the isomorphism is required to be the identity restricted
to the boundary. As before we define vector spaces $\cN$ (for closed
graphs), $\cN_C$ (for graphs with boundary $C$), and pairings $\cN_C
\otimes \cN_C \to \cN$.

\begin{theorem} \label{grpos2b}
The above pairing, for vertices with symmetry groups
and gluing along cut points, is positive.
\end{theorem}

\begin{proof}
The proof is similar to that of the previous theorem, except that
for each vertex type we choose a tensor that has precisely the
symmetry (neither more nor less) specified for that vertex type. We
do this this by first choosing unsymmetric tensors (with
algebraically independent entries, as before), and then
symmetrizing. (The symmetry group of a vertex of valence $k$ acts on
$\{1, \ldots, k\}$, which induces an action on $V^{\otimes k}$.)

We must show that for graphs $A, B \in \dot\cN_C$, $Z_i(A) = Z_i(B)$
for sufficiently large $i$ only if $A \cong B$ (via allowed
symmetries). $Z_i(A)$ is a contraction of symmetrized tensors. Since
symmetrizing commutes with contraction, we can also view $Z_i(A)$ as
a sum of contractions of the original unsymmetric tensors. There
will be one such summand for each gluing of the vertices of $A$
which is compatible (via vertex symmetry groups) with a gluing which
produces $A$ from its vertices. It follows that if $i$ is greater
than the number of edges in $A$ we can read off from $Z_i(A)$
(thought of as a polynomial in the unsymmetric tensor indices) the
equivalence class of the adjacency matrix of the vertices of $A$
(relative to $C$). By algebraic independence, $Z_i(A) = Z_i(B)$ only
if corresponding polynomials are identical. It follows that the
adjacency matrix of $A$ can be transformed into that of $B$ via
vertex symmetries.  In other words, $A \cong B$.
\end{proof}

\medskip

\begin{figure}[htpb]
\labellist
\small\hair 2pt
\pinlabel $\beta$ at 232 218
\pinlabel $\gamma$ at 250 50
\pinlabel $\alpha$ at 35 160
\pinlabel $A$ at 450 260
\pinlabel $H$ at 450 150
\pinlabel $B$ at 450 40
\pinlabel $\text{cuts}$ at 490 203
\pinlabel $\text{cuts}$ at 490 100
\endlabellist
\centering
\includegraphics[height=2in]{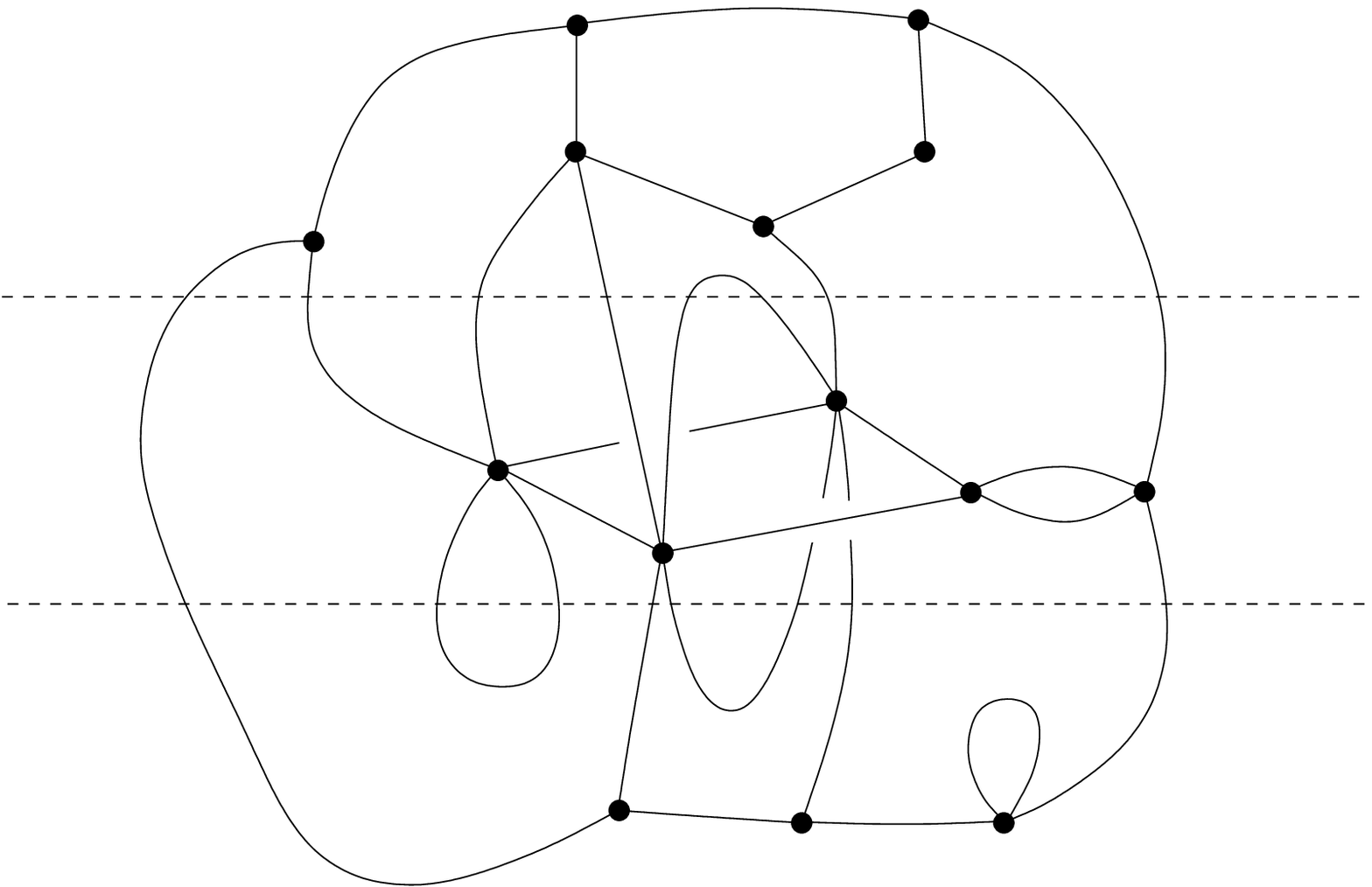}
\caption{} \label{upperlower}
\end{figure}

The third and final case of graph positivity we consider concerns
gluing graphs along an intermediate subgraph $H$ instead of gluing
along edge cut points. This is analogous to the 3-manifold case
where $S$ has no JSJS cusp components; $H$ corresponds to the unions of
the JSJ pieces which intersect $S$,
while $A$ and $B$ correspond to the remaining JSJ pieces.
It is easy to extend the
arguments given below to the case of gluing along the disjoint union
of an intermediate subgraph and cut points, analogous to the case
where $S$ has both components which are non JSJS cusps as well as components which are.
In the most general
case, we would allow arcs ($\alpha$ in Figure \ref{upperlower}) in $H$ with two cut
points  and meeting no vertices and also arcs ($\beta$ and $\gamma$)
in $\dot\cN_C$ which meet no vertices.  The basic idea is to always
associate the identity $e^i \otimes e_i$ to $\alpha$, $e^i \otimes
e^i$ to $\beta$, and $e_i \otimes e_i$ to $\gamma$.  We leave the
details of this extension as an exercise for the reader.

For simplicity we assume that all vertex types are fully symmetric;
the symmetry group of each vertex type is the full permutation group
of the edges emanating from it.
(We will see below that the symmetries of a Seifert fibered JSJ piece
act as the full permutation group on its cusps.)

We assume that the boundary of $H$ is partitioned into ``upper'' and
``lower'' cut points, both identified with $C$, and that $H$ is
equipped with an order two automorphism $\iota$ which permutes the upper and
lower cut points and fixes all vertices and non-cut edges of $H$
(see Figure \ref{upperlower}). This automorphism is analogous to
reflecting the middle level of the 3-manifold $AA$ (or $BB$) across
$S$.

We must modify the definition of $\cN_C$ to take into account
automorphisms of $H$. (So $\cN_C$ depends not just on $C$ but also
on $H$.) Two graphs in $\dot\cN_C$ are considered isomorphic if
there is an isomorphism of the underlying unlabeled graphs which
preserves vertex labels and induces a permissible permutation of the
boundary $C$. (Previously we required the isomorphism to be the
identity on $C$.) A permissible permutation is one which preserves
the adjacency relation to each vertex of $H$. Equivalently, a
permissible permutation is one which extends to an automorphism of
$H$ which fixes all vertices of $H$ and all non-cut edges of $H$
(and fixes doubly cut edges if these are permitted in $H$). Here
we're identifying $C$ with the upper (or lower) boundary of $H$. For
example, the permissible permutations of $C$ in Figure
\ref{upperlower} can exchange cut points two and three, four and
five, or six and seven (counting from the left). (If we considered
the unsymmetrical context, ``permissible'' would have a more
restricted meaning.)

Given $A, B \in \dot\cN_C$, define the closed graph $AHB$ by gluing
the boundary of $A$ to the upper boundary of $H$ and gluing the
boundary of $B$ to the lower boundary of $H$. This induces a pairing
$\cN_C \otimes \cN_C \to \cN$.

\begin{theorem}
\label{grpos3} The above pairing, for fully symmetric vertices and
gluing along a subgraph, is positive.
\end{theorem}

\begin{proof}
As in the previous two proofs, we will construct a family of graph
TQFTs whose tuple of partition functions is a diagonally dominant
complexity function on closed graphs. We can think of $Z_i(H)$ as an
operator
\[
    Z_i(H) : V_i^{\otimes C} \to V_i^{\otimes C}
\]
from the Hilbert space of the lower cut points of $H$ to the Hilbert
space of the upper cut points of $H$. The isomorphisms used to
define $\dot\cN_C$ act on $V_i^{\otimes C}$ with invariant subspace
$Z_i(\cN_C):= \Inv(V_i^{\otimes C}) \subset V_i^{\otimes C}$. We
will show below that the vertex tensors from which the TQFT is
constructed can be chosen so that the restriction of $Z_i(H)$ to
$Z_i(\cN_C)$ is strictly positive (for all $i$).  That is, there
exists an operator $U$ such that $Z_i(H) = U^\dagger U$ and $U$ is
injective on $Z_i(\cN_C)$. We then have (using ``bra-ket'' notation,
and the dagger for adjoint)
\begin{eqnarray*}
    |Z_i(AHB)| &=& |\langle Z_i(A) | Z_i(H) | Z_i(B)\rangle| \\
        &=& |\langle Z_i(A) | U^\dagger U | Z_i(B)\rangle| \\
        &=& |\langle U(Z_i(A)), U(Z_i(B))\rangle| \\
        & \le & \left(|\langle U(Z_i(A)), U(Z_i(A))\rangle|\cdot |\langle U(Z_i(B)), U(Z_i(B))\rangle|\right)^{\frac{1}{2}} \\
        & \le & \max(Z_i(AHA), Z_i(BHB)) ,
\end{eqnarray*}
with equality only if $U(Z_i(A)) = U(Z_i(B))$. Since $U$ is
injective on $Z_i(\cN_C)$, this happens only if $Z_i(A) = Z_i(B)$.
Since all entries in the symmetric tensors for each vertex type in
$A$, $B$, and $H$ may be assumed to be algebraically independent, a
slight extension of the previous proof shows that this happens for
all $i$ only if $A \cong B$.

To illustrate this extension, imagine that $Z_i(\cN_C)$ is $S^2(V)$,
the symmetric square of $V$.  This would result from two upper cut
points emanating from one vertex of $H$.  Then the symmetrized
relative partition functions in $Z_i(\cN_C)$ are $$Z_i(A) =
Z_i(A)_{p,q} + Z_i(A)_{q,p} \text{ and } Z_i(B) = Z_i(B)_{p,q} +
Z_i(B)_{q,p}$$ where $\{p,q\}$ is a multi-index for $S^2(V)$.
Algebraic independence ensures that these expressions are equal if and only if
they have the same terms, up to permutation.  In general,
$Z_i(A) = Z_i(B)$ if and only if their {\em unsymmetrized} partition functions
agree after some permissible permutation of cut edges.  Then for $i$
large enough to label all edges of $A$ ($B$) differently, we
conclude $A \cong B$

It remains to show how to construct the (graph) TQFT so that
$Z_i(H)$ is a positive operator on $Z_i(\cN_C)$. In what follows, we
suppress the index $i$.

Consider a single vertex type $a$ with $n:=n(a)$ emanating cut
edges. Denote by $V^+ \subset V$ the non-negative cone, consisting
of non-negative sums of elements from a fixed orthonormal basis
${e_1, \dots, e_d}$. We must find a tensor $T^a \in S^nV$ with the
following property. Suppose $v$ is a vertex of type $a$ which
appears as a vertex in $H$ with $k$ upper and lower cut edges which
are interchanged by $\iota$, and $n-2k$ ``horizontal edges'' (i.e.\
edges joining the vertex to other vertices in $H$). Associated to
$k$ lower cut edges (resp. upper cut edges) there is a vector space
$V_l$ (resp. $V_u$) which parameterizes the possible inputs along
these edges. Explicitly, each such vector space is isomorphic to the
$k$th symmetric tensor power of $V$.  The involution $\iota$ induces
an isomorphism
$$\iota:V_u \to V_l$$
and therefore it makes sense to say that an operator from $V_l$ to $V_u$ is
{\em positive symmetric}.

The property we want is that after contracting $T^a$ with $n-2k$
vectors in $V^+$ (corresponding to input from the horizontal edges
emanating from $v$), the result should be a positive symmetric
operator from $V_l$ to $V_u$. That is, for any $k$ satisfying $2k\le
n$, and any vectors $v_j^+ \subset V^+$ for $2k+1 \leq j \leq n$,
\begin{equation}\label{positive_after_horizontal_contraction_eqn}
T^a(v^+_{2k+1}, \cdots, v^+_n):V_l \to V_u
\text{ is positive symmetric}
\end{equation}

Roughly, the entries of $T^a$ (in the fixed basis $e_1,\cdots,e_d$)
are chosen by any formula which extravagantly weights index
coincidences, and then perturbed to be algebraically independent. An
example of what we mean by such an extravagant formula is given in
the following equation:
\begin{equation}
\label{eqn**} T_{i_1, \dots, i_n}^a \approx 1\cdot
(2^{d^n})^{\#\text{ of pairs}}\cdot (2^{d^{n^n}})^{\#\text{ of triples}} \cdot (2^{d^{n^{n^n}}})^{\#\text{ of quadruples}}\cdots
\end{equation}
for any assignment of indices $i_1,\cdots,i_n$ chosen from the set
$\lbrace e_1,\cdots,e_d\rbrace$. Here ``pairs'' is the number of
unordered pairs of indices with the same label, ``triples'' is the
number of unordered triples of indices with the same label, and so
on, so for example $T^a_{e_1,e_2,e_1,e_3,e_2,e_1,e_3,e_4}$ contains
two pairs and one triple, and its value should have size on the
order of $(2^{d^n})^2\cdot 2^{d^{n^n}}$. There is a convenient
criterion, originally due to Sylvester, to verify that a symmetric
matrix $M_{ij}$, $1 \leq j \leq m$ acts as a positive operator.

\begin{lemma}
\label{lemma0.7} $M_{ij}$ is positive if and only if the partial
determinants are all positive, i.e.\ $\det(M_{ij}) > 0$ for $i, j \leq
l$, where $1 \leq l \leq n$.
\qed
\end{lemma}

Easier than checking that $T^a$ of the form suggested in \ref{eqn**}
satisfies the criterion of Lemma~\ref{lemma0.7}, is to observe that the
lexicographic tuple:
\begin{equation}
\label{eqn***} (\#  2k \text{ tuples}, \#2k-1 \text{ tuples}, \cdots,
\# \text{pairs})
\end{equation}
is diagonally dominant under concatenating (symmetrized) strings of
$k$ symbols drawn from a finite set. Thus, if the finite ordered set
consisting of values of the tuple (\ref{eqn***}) is embedded in
$\R^+$ with sufficiently (depending on $d$) rapid decay, then
positivity follows from estimating determinants (as simply the
``order of magnitude'' of the product of diagonal elements) and then
applying Lemma~\ref{lemma0.7}.  In more detail, consider what we
need to know to safely estimate a determinant of a matrix $M$ simply
as the product of the diagonal entries.  It would suffice if all of
the off-diagonal entries are much smaller than all of the diagonal
entries.  Unfortunately, we are not in this situation but we can
easily arrange (normalizing all entries to be real and $\geq 1$)
that
$$M_{ij} < \frac{1}{n!}\sqrt{M_{ii} M_{jj}}$$ for all $i \neq j$ and
$n$ equal to the dimension of the matrix $M$.  Given a permutation
$\sigma \neq id$, it easily follows that $$\prod_i M_{i, \sigma(i)}
< \frac{1}{n!} \prod_i M_{ii}$$ Thus, the diagonal term dominates
the calculation of the determinant.

Similarly, when finitely many operators which weight index
coincidences sufficiently extravagantly (depending on $d$ and $H$)
are (partially) contracted according to the combinatorics of $H$,
the result is a {\em positive} operator $T(H)$ with precisely the
desired symmetries.  Since the property of being positive is open,
we can perturb the tensors to be algebraically independent. This
completes the proof.
\end{proof}

\subsection{Tensors for JSJ pieces} \label{actjsjsect}

The main ingredient $c_t$ in the assembly complexity will be similar
to the graph complexity defined above. In broad outline, we will
assign vector spaces to JSJ tori and appropriate tensors to JSJ
pieces, then define the complexity of a closed manifold to be the
result of a big tensorial contraction.

The gluing graph of JSJ pieces has three new features not present in
the naked graph context:
\begin{enumerate}
\item{The vertices (JSJ pieces) are now
oriented manifolds, so we will use complex scalars $\C$ and assign
to complex conjugation the role of keeping track of orientation.}
\item{Each vertex has its own group of symmetries.  Ideally, we should
employ tensors with exactly the right symmetry:  Too little symmetry
and the partition function will not be well defined, too much and
subtle gluing distinctions that should distinguish $A$ from $B$ will
be lost.  In the case of hyperbolic vertices, the symmetries
constitute a finite group $G$ and we will build appropriately
symmetric tensors from many smaller tensors or ``micro-vertices''
situated within each fundamental domain of $G$. Infinite symmetry in
the Seifert fibered cases is a problem and requires a separate
trick, namely passing to modular quotients.}
\item{The edges of the JSJ graph carry interesting
information.  If we assume each vertex type comes with a marking
(i.e.\ a choice of co-ordinates) on (the homology of) its cusps,
then ``orienting an edge'' is tantamount to choosing an element in
$g \in \GL(2,\Z)$ which specifies the gluing (reversing the
orientation sends $g$ to $g^{-1}$).  Thus, the vector space $V$ on
which we build tensors should have an action of $\GL(2,\Z)$ (but
also a large $GL(2, \Z)$-invariant factor).}
\end{enumerate}

The fact that mapping class groups of Seifert fibered spaces are
typically infinite needs to be addressed first since its solution
determines the vector space $V$ on which we build tensors. Actually,
we will construct a family of vector spaces $V_{m,m'}$ and tensor
complexities $c_t^{m,m'}$, for $m,m' = 2, 3, 4$ and so on, where the
vector space $V_{m,m'}$ is isomorphic to $\C^{m^2m'}$ and admits a
natural action of $\GL(2,\Z)$ factored through the finite group
$\GL(2, \Z_m)$. Explicitly, a basis for $V_{m,m'}$ is given by the
elements in the finite Abelian group $H_1(T^2;\Z_m)$ cross the set
$\{1,2,\dots,m'\}$ (here and throughout we use the ``topologists
notation'' $\Z_m$ for the finite cyclic group $\Z/m\Z$). The natural
action of $\GL(2,\Z)$ on $H_1(T;\Z)$ induces an action of
$\GL(2,\Z)$ on $H_1(T^2;\Z_m)$ which factors through the quotient
$\GL(2, \Z_m)$. While no single representation of $\GL(2,\Z)$ is
faithful, the intersection of the kernels is trivial, which is
sufficient for our purposes.
The action on $\{1,2,\dots, m'\}$ is trivial.
The numbers $m$ and $m'$ have no logical relation except
that they both will be sent to infinity (for different reasons). For
notational convenience, we set $m=m'$ and denote $V_{m,m'}$ as $V_m$
and similarly $c_t^{m,m'}$ as $c_t^m$.

Unfortunately, a technical problem prevents us from getting the
symmetry group exactly right. We will build in a little too much
symmetry for the tensors associated to Seifert fibered JSJ pieces.
The resulting partition function would not be able to distinguish
manifolds which differ by a relation which we call ``fiber
slipping''. Finally, this deficiency is anticipated and corrected by
a local modification of the tensor associated to each JSJ piece.

Let $m \geq 2$ be an integer and let $V_m$ be the $\C$ vector space
spanned by the finite set $H_1(T^2; \Z_m) \times \{1,2,\dots,m\}$.
Let $X$ be a JSJ piece, which in our context is either finite volume
hyperbolic or Seifert fibered.  We will associate to $X$ a carefully
chosen but suitably generic tensor $T_X$ with $n$ covariant indices
where $n$ is the number of cusps of $X$. Reversing the orientation
of $X$ conjugates the entries of $T_X$. That is, $T_{\overline{X}} =
\overline{T}_X$. Let $h:X \rightarrow X$ be a homeomorphism
(orientation preserving or reversing). The induced action of $h$ on
$T_X$ will factor through the action of $h$ on $H_1(\partial X;\Z_m)
\times \{1,2,\dots,m\}^c$, where $c$ is the number of cusps.  The
action on the set permutes factors according to the action on the
end $end(X)$, which, of course, embeds in the symmetric group $S_c$.
Thus, the mapping class group of $X$ acts on $T_X$ but factored
through the finite group $\GL(2, \Z_m) \times S_c$ for each $m$. If
$l \in H_1(\partial X, \Z_m) \times \{1,2,\dots,m\}$ is a
multi-index then symmetry demands:
\begin{equation}
\label{eqn4} h(T_l) =
\begin{cases}
T_{h_*(l)} \text{ if } h \text{ preserves orientation} \\
\overline{T}_{h_*(l)} \text{ if } h \text{ reverses orientation}\\ \end{cases}
\end{equation}

For example, consider $X$ with three cusps reflected by $h$ as shown
in Figure~\ref{3cusps} below. We use $h$ to identify the groups at
the upper and lower cusps and the ``bar'' on $k$ to denote the
induced involution on the middle level group to obtain:
\begin{equation}
\label{eqn5} h(T_{ijk}) = \overline{T}_{ji\overline{k}}
\end{equation}

\begin{figure}[htpb]
\labellist
\small\hair 2pt
\pinlabel $X$ at 90 110
\pinlabel $h$ at -10 80
\pinlabel $T_i^2$ at 135 150
\pinlabel $T_j^2$ at 135 10
\pinlabel $T_k^2$ at 220 80
\endlabellist
\centering
\includegraphics[height=1.5in]{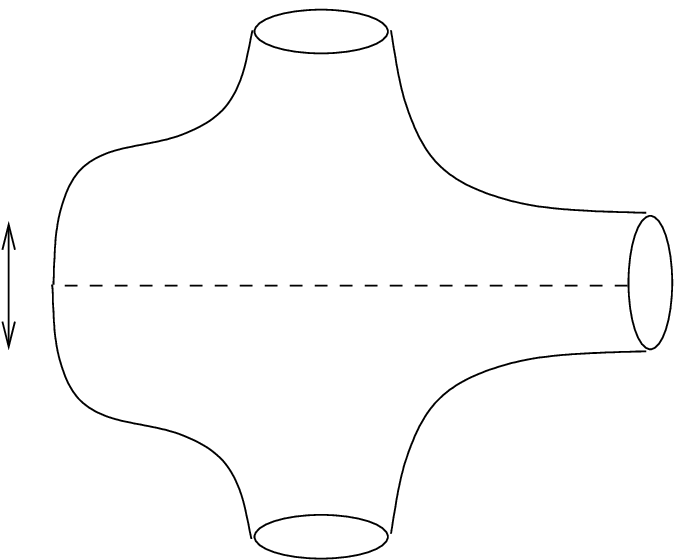}
\caption{} \label{3cusps}
\end{figure}

In what follows, a {\em reflection symmetry} $r$ (or just {\em
reflection} for short) refers to an orientation reversing involution
whose fixed set $\Sigma \subset X$ is $2$-dimensional. Note that
when $X$ is hyperbolic, $r$ may be taken to be an isometry.  For
each SF $X$ and with respect to any reflection symmetry $r$, the
tensors $T_X$ will be constructed to have the positivity property
$P$ below. Positivity is an open (convex) condition, so we may first
concentrate on producing positive Hermitian tensors for each JSJ
piece $X$ and then perturb the construction to achieve genericity.

\vskip 12pt

We fix $m$ and abbreviate $V_m$ by $V$. Note that $V$ has a
canonical basis, namely the set $H_1(T^2;\Z_m) \times
\{1,\dots,m\}$. We let $V^+ \subset V$ denote the cone of
non-negative (real) linear combinations of these basis elements.
Suppose $r: X \rightarrow X$ is a reflection.  Then $r$ partitions
the cusps of $X$ into three groups: the {\em horizontal} cusps which
the fixed set $\Sigma$ meets, and two subsets, the {\em upper} and
{\em lower} cusps which are disjoint from $\Sigma$, and interchanged
by $r$. (We arbitrarily designate one side of $\Sigma$ as ``upper''
and the other as ``lower''.) Now suppose $X$ has $n$ cusps: $k$
upper, $k$ lower, and $n-2k$ horizontal (all with respect to $r$).
If supplied with vectors $v_{2k+1}, \cdots, v_n \in V^+$ at the
horizontal cusps, $T_X$ becomes an operator
$$\displaystyle T_X(v_{2k+1}, \dots, v_n): \bigotimes_k V \rightarrow
\bigotimes_k V$$ from the vector space associated to the lower
cusps, to the vector space associated to the upper cusps.
Furthermore, the homeomorphism $r$ induces a natural {\em
complex-linear} isomorphism between these two vector spaces, so it
makes sense (using this identification) to ask if this operator is
positive, for all reflections $r$ and all positive $\{v_i\}$
inserted at the horizontal cusps of $r$.

We explain the overall strategy.  The gluings $AB$, $AA$, and $BB$
along $S$ are analogous to the graph case of Theorem \ref{grpos3},
where vertices lie on the gluing set.  (Both hyperbolic and SF JSJ
pieces may be divided in two by $S$.)  So, we will use some kind of
$\langle A^- | H | B^- \rangle$ format to calculate $c_t^m$ and
establish a topological Cauchy-Schwartz inequality from a (weak)
positivity property of $H$.  $H$ is an operator associated to the
``middle third'' while $A^-$ ($B^-$) is, by abusing notation
slightly by dropping the $Z$, the vector (covector) associated to
the upper and lower thirds.  We define the {\em geometric middle
third} (GMT) as the gluing along $S$ of the JSJ pieces incident on
$S$ together with a copy of $T^2 \times I$ (analogous to $\alpha$ in
Figure \ref{upperlower}) for each JSJS cusp.  Thus, the GMT is made
by gluing boundary pieces and adding a $T^2 \times I$ for each
horizontal cusp.  Recall from Definition
\ref{vertical_horizontal_cusp_definition} that the JSJS cusps are
the torus components of $S$ which are also JSJ tori of the glued up
manifold $AB$.

There is a very nice conjecture (appendix \ref{appendix_conjecture})
that implies the existence of ``positive'' tensors for hyperbolic
$3$-manifolds.  With that input, $H$ could be assembled from tensors
corresponding to JSJ pieces of the GMT. Since we lack a proof for
this conjecture we will actually define $H$ from a smaller {\em
algebraic middle third} (AMT).  The AMT, although primarily a device
for constructing the operator $H$, does have a geometric
realization:  $\text{AMT} = S \cup (\text{SF pieces of GMT glued
along } S)$.  Thus, the difference $\displaystyle \text{GMT}
\setminus \text{AMT} = \bigcup_{\text{hyperbolic } X \in \text{GMT}}
(X \setminus S)$.  All that remains in the AMT of the hyperbolic
part of the GMT is $S$.

\begin{figure}[htpb]
\labellist
\small\hair 2pt
\pinlabel $\beta$ at 160 245
\pinlabel $\gamma$ at 160 70
\pinlabel $\alpha$ at 295 160
\pinlabel $\text{one}$ at 345 167
\pinlabel $\text{cusp}\times I$ at 345 153
\pinlabel $6\text{ boundary pieces}$ at -65 167
\pinlabel $\text{glued in pairs}$ at -65 153
\pinlabel $\text{upper third }A^-$ at 450 260
\pinlabel $\text{GMT}$ at 450 160
\pinlabel $\text{lower third }B^-$ at 450 60
\endlabellist
\centering
\includegraphics[height=2.7in]{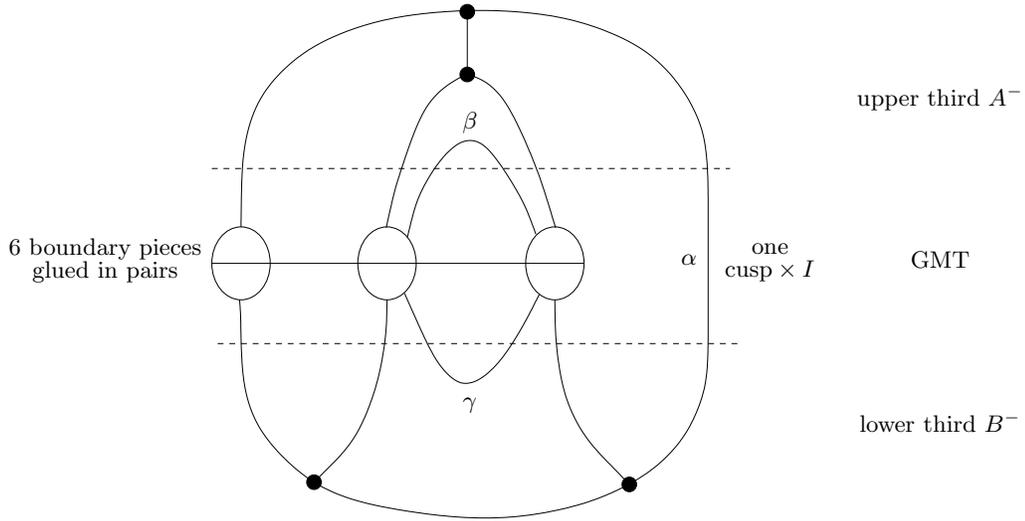}
\caption{$AB$ is decomposed into thirds $A^- \cup M \cup B^-$} \label{thirds}
\end{figure}

We will also need a notation for the complement of the AMT: $A \cup
B \setminus \text{AMT} = A^\sim \cup B^\sim$, so $A^- \subset
A^\sim$ and $B^- \subset B^\sim$.  Thus, the formal scheme for
calculating $c_t^m$ is more accurately written:
\begin{equation} \label{c_t^m_eqn}
c_t^m(AB) = \langle c_t^m(A^\sim) | H^m | c_t^m(B^\sim) \rangle
\end{equation}
Note the close analogy with Theorem \ref{grpos3} (see Figure
\ref{thirds} for clarification).

What makes a calculation of this form possible is that the tensors
$T_X$ we build for hyperbolic $X$ are actually a partial contraction
of a very large tensor network $N_X$ (enjoying exactly the same
symmetry group $G$ as $X$) which crosses $S$ orthogonally.  This
allows $T_X$, for $X$ hyperbolic and in the GMT, to be represented
partly in the bra and partly in the ket of (\ref{c_t^m_eqn}).

The key will be to show that $H^m$ is strictly (``block'') positive
(to be defined shortly) on the space of bras (kets) that it will
ultimately be paired with in (\ref{c_t^m_eqn}).  This is where the
symmetry issues arise.  The final step will be to recover the
information lost by settling for tensors $T_X$, $X$ SF which are
slightly too symmetric and hence only ``block'' positive.

We take a moment to clarify the correspondence to the graph case.
JSJS cusps correspond to isolated edges in $H$ with a cut point at
either end, like the arc $\alpha$ in Figure~\ref{upperlower}. Note
that two cusps in the boundary of the middle third may cobound a
single $T^2\times I$ component in $A^-$ or $B^-$; these correspond
to isolated edges, such as the arcs $\beta$ and $\gamma$ in
Figure~\ref{upperlower} in the graph case. Also as in the graph
context, when we come to assign tensors to vertices, we adhere to
the convention that doubly cut edges (of type $\alpha,\beta,\gamma$)
get identity tensors (up to index raising and lowering).

We now construct the tensors $T_X$ for $X$ SF. In contrast to the
rigidity of hyperbolic manifolds, Seifert fibered JSJ pieces $X$
with nontrivial boundary have infinite mapping class groups (whose
elements we also refer to as {\em symmetries}). The restriction of
the group of symmetries to the boundary of a JSJ piece is the
semidirect product of two groups. The first, $\til{P}$, is an
extension of $\Z_2$ by the full permutation group $P$ of the
boundary components of the quotient orbifold $Q$:
\[
    1 \to P \to \tilde{P} \to \Z_2 \to 1 .
\]
The image in $\Z_2$ measures the action on orientation.

The second group is generated by {\em arcs of Dehn twist} (hereafter
abbreviated to ADT) restricted to the boundary. An ADT connects a
positive Dehn twist on a boundary component to a negative Dehn twist
on a (possibly) distinct boundary component.  The picture to have in
mind is to start with a proper embedded arc in the quotient orbifold
$Q$. Sitting above this arc is a proper annulus, which cobounds
circles in two (possibly) distinct boundary components. Cut along
this annulus and reglue after twisting once around each circle fiber
in the annulus.

Reparameterizing a component $X$ (e.g.\ by an ADT) changes the way
in which $X$ is glued to its neighbors, and how it is glued to
itself (if at all).  It would be desirable to find tensors $T_X$ for
Seifert fibered pieces $X$ which are positive on exactly the
subspace of vectors invariant under all such parameterizations;
however, we were not able to find such tensors directly. Instead, we
first define ``blocky'' tensors which are completely insensitive to
the fiber summand of each $H_1(\partial_i X; \Z_m)$ index
(effectively ignoring the troublesome ADT's).  Second, we determine
that the resulting tensors, while not leading to a strictly
diagonally dominant complexity function on homeomorphism classes of
manifolds, nevertheless would yield a diagonally dominant complexity
function on manifolds up to the weaker relation of ``fiber slip
homeomorphism'' (fish). Finally, we enhance our tensors
$\widehat{T}_X \leadsto T_X$ by letting them feel their immediate
neighbors in the JSJ decomposition to resolve the fishy ambiguity.
Morally, fiber slip homeomorphism is the equivalence relation
obtained by failing to keep track of a certain relative Euler class
(living in certain $\Z$ torsors). The remedy is to find canonical
co-ordinates on these $\Z$ torsors, and favor relative Euler classes
which arise by doubling.

\vskip 12pt

The first step is to define ``block positive'' for Seifert fibered
pieces.

For the moment, fix a fibered structure on the JSJ piece $X$. The
choice of fiber structure picks out a canonical $\Z_m$ summand in
$H_1(\partial_i X;\Z_m)$. We let $\pi$ be the homomorphism which
quotients out this summand:
\begin{align*}
\pi: H_1(\partial_i X ; \Z_m) \times \{1,2,\dots,m\} &\to
\\ &H_1(\partial _i X; \Z_m)/\langle[\text{fiber}]\rangle \times
\{1,2,\dots,m\} \\ &\cong \Z_m \times \{1,2,\dots,m\}
\end{align*}
The map $\pi$ induces a map from $V_m \cong \C^{m^3}$ to the complex
vector space, isomorphic to $\C^{m^2}$, spanned by $H_1(\partial _i
X; \Z_m)/\langle[\text{fiber}]\rangle \times \{1,2,\dots,m\}$. The
``blocky'' tensors we will construct have entries which only depend
on the image of each co-ordinate in $\C^{m^2}$ under this
projection. Fix a reflection $r$, and partition the cusps of $X$
into upper, lower and horizontal cusps as before, where the
horizontal cusps are those meeting the fixed set $\Sigma$ of $r$. We
want $\widehat{T}_X$ to be some totally symmetric $n$-tensor
depending only on $\pi(\text{index})$ for each index. The operator
$\Oh := \widehat{T}_X(v_{2k+1},\dots,v_n)$ maps between associated
symmetric powers of $S_{\text{lower}}$, $S_{\text{upper}}$ of $V$
and descends to $\Oh^\pi$ mapping between symmetric powers $S^\pi$
of $\pi(V) := \C[\pi(H_1(\partial_i X; \Z_m))] \times
\{1,2,\dots,m\}$.

\begin{figure}[hbpt]
\[\begindc{0}[15]
    \obj(0,4){$S_{\text{lower}}$}
    \obj(0,8){$V_{\text{lower}}$}
    \obj(14,4){$S_{\text{upper}}$}
    \obj(14,8){$V_{\text{upper}}$}
    \obj(0,0){$S^{\pi}_{\text{lower}}$}
    \obj(14,0){$S^{\pi}_{\text{upper}}$}
    \mor(0,8)(0,4){$\text{symmetrize}$}
    \mor(14,8)(14,4){$\text{symmetrize}$}
    \mor(1,8)(13,8){}
    \mor(1,4)(13,4){$\Oh := \widehat{T}_X(v_{2k+1}, \dots, v_n)$}
    \mor(0,4)(0,0){$\pi$}
    \mor(14,4)(14,0){$\pi$}
    \mor(1,0)(13,0){$\Oh^\pi$}
\enddc\]
\caption{} \label{quotientdiagram}
\end{figure}

We now define block positivity as follows:
\begin{definition}
\label{block_pos_def} A tensor $\widehat{T}_X$ assigned to a Seifert
fibered JSJ piece X satisfies {\em property $P$} if, for all
reflections $r$ and positive inputs $\lbrace v_{2k+1},\cdots,v_n
\rbrace$ at the horizontal cusps, the operator
$\widehat{T}_X(v_{2k+1}, \dots,v_n)$ is {\em block positive}; i.e.
if $\Oh^\pi$ as above is positive.
\end{definition}

Actually, property $P$ will only be important for JSJ pieces which
lie the middle third $M$ of $AB$, $AA$, or $BB$, so it is harmless
to restrict attention to those Seifert fibered $X$ which arise by
gluing two boundary pieces along (part of) $S$. Thus the pieces in
question must be sufficiently large, and (in order to participate in
forming $H$), must have nonempty boundary.
According to the
analysis of exceptional cases in \S\ref{c_S_section}
(Theorem~\ref{exceptional_list_theorem},
Theorem~\ref{mapping_class_exceptions}, and
Lemma~\ref{SF_DD_lemma}), the only sufficiently large JSJ piece
$X$ with boundary that might arise and that lacks a unique Seifert fibered structure 
is $S^1 \tilde{\times} S^1 \tilde{\times} I \cong M(+0, 1; 1/2, 1/2)$. But
although this manifold admits a second nonisomorphic Seifert fibered
structure, this latter structure is excluded by Convention~\ref{twisted_bundle_Klein_convention}.
Thus without loss of generality we are
justified is assuming a unique Seifert fibered structure on $X$ for
the purpose of defining $P$. This is crucial, as it justifies our
use (above) of a fixed fibered structure on $X$, and thus makes the
quotient operator $\Oh^\pi$ well defined as an operator from
$S^\pi_{\text{lower}}$ to $S^\pi_{\text{upper}}$.

\vskip 12pt

For $X$ Seifert fibered, we must construct a generic
$\widehat{T}_X$ satisfying property $P_S$:
$$\widehat{T}_X \in \C[H_1(\partial X ; \Z_m)] \cong
\bigotimes_{i=1}^n \C[H_1(\partial_i X ; \Z_m)] =
\bigotimes_{i=1}^n V_m$$ Recycling a trick from the proof of
Theorem~\ref{grpos3}, block positivity will be ensured by choosing
the entries of $\widehat{T}_X$ according to a function which
extravagantly weights coincidences. The $\widehat{T}_X$ entry should
be (approximately) given by composing the string
\begin{equation}
\label{eqn****} (\#n\text{-tuples, \#}(n-1)\text{-tuples, } \dots
\text{ ,\#pairs)}
\end{equation}
into the reals $\R$ by a function such as \ref{eqn**}. Coincidence
tuples are diagonally dominant functions on multi-indices under
concatenation, and so are such compositions into $\R$. Here we mean
$k$-tuples of indices which become the same in the quotient group
$H_1/\text{fiber}$. We also insist that the entries of
$\widehat{T}_X$ are algebraically generic subject to the relation
$\widehat{T}_{\overline{X}} = \overline{\widehat{T}}_X$ and the
insensitivity to fiber coordinates.

To repeat, the entries $(\widehat{T}_X)_{i_1,\cdots,i_n}$ of
$\widehat{T}_X$, where the $i_j$ denote elements of the various
$H_1(\cusp;\Z_m)$, depend only on the images of these $i_j$ under
the projection $\pi$. It follows that the operator $\Oh$, obtained
by contracting $\widehat{T}_X$ with positive vectors
$v_{2k+1},\cdots,v_n \in V^+$ associated to horizontal cusps, is not
positive but merely non-negative. It is the induced operator
$\Oh^\pi$ which is positive. The operator $\Oh$ itself has null
directions spanned by vectors of the form:
\begin{equation} \label{eqn*****}
    e_{b_1,f_1} \otimes \cdots \otimes e_{b_{i-1},f_{i-1}}
        \otimes (e_{b_i, f_i} - e_{b_i, f_i'}) \otimes \cdots \otimes e_{b_{i+1},f_{i+1}}
        \otimes e_{b_k,f_k}
\end{equation}
where the pair $f_i, f_i'$ are distinct values of the
fiber co-ordinates in some basis for $H_1(\cusp;\Z_m)$.

We have proven:
\begin{lemma} \label{prop_p_lemma}
Any Seifert fibered $X$ with unique Seifert fibering admits a tensor
$T_X$ satisfying property $P$.
\end{lemma}

\vskip 12pt

Next we construct the tensor $T_X$ when $X$ is hyperbolic. It is
important that $T_X$ be constructed from a (partial) contraction of
a tensor network $N$ defined on a fundamental domain $W$ of $X$ with
respect to its full symmetry group $G$, and transported around $X$
by the group.  This way, $N$ will enjoy the symmetries of $X$.  We
must take a little further care to ensure that $N$ has no additional
symmetry as a labeled graph.  We will use our usual $V_m$ to label
all edges.

The construction is not difficult, but involves a number of steps
which are described in the next several pages. Recall that the goal
of this paper is to establish positivity of the universal
$(2+1)$-dimensional TQFT.  If one knew a particular unitary TQFT
that was sufficiently rich and easy to analyze, it might be possible
simply to check positivity there. A genuine TQFT has the {\em
factorization property}: whenever $X$ can be decomposed along a
surface into pieces $Y$ and $Z$ one has a factorization $T_X = T_Y
T_Z$, a compact notation for the appropriate tensorial contraction.
However, one can construct weaker substitutes for a full TQFT which
have the factorization property only for certain prescribed
decompositions. In a sense, we are defining a ``weak TQFT,'' in the
style of Penrose, for hyperbolic manifolds. But it is only
factorizable with respect to decompositions along the fixed sets of
reflections, and not with respect to arbitrary topological
decompositions.

\begin{figure}[htpb]
\labellist
\small\hair 2pt
\pinlabel $\text{cusps}$ at 170 60
\pinlabel $\frac{1}{2}\text{ cusp}$ at 380 130
\pinlabel $\frac{1}{|H|}\text{ cusp}$ at 460 10
\endlabellist
\centering
\includegraphics[height=1.5in]{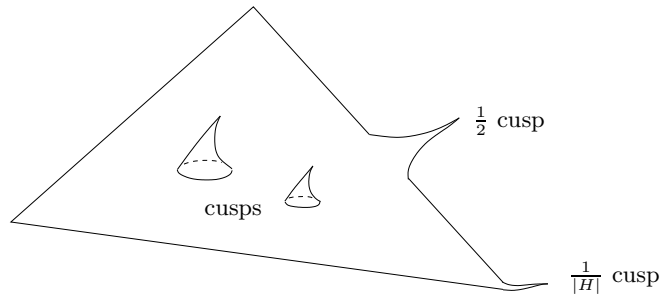}
\caption{Schematic fundamental domain $W$, with cusps and fractional cusps} \label{cusps}
\end{figure}

Let $G$ be the group of isometries of $X$. Since $X$ is a finite
volume hyperbolic manifold, $G$ is finite. Let $W \subset X$ be a
piecewise smooth fundamental domain for the action of $G$ on $X$.
The fundamental domain $W$ may contain cusps in its interior and
also {\em fractional cusps}, meaning a cusp $C$ which is ``paved'' by
$k$ copies of $W$. We think of each copy of $W$ as containing a
$1/h$ fraction of $C$, where $h$ is the order of the symmetry
subgroup $H \subset G$ stabilizing the cusp in question.  See
Figure~\ref{cusps}.

\begin{figure}[htpb]
\labellist
\small\hair 2pt
  \pinlabel $V_m$ at 55 173
  \pinlabel $V_m$ at 20 148
  \pinlabel $V_m$ at 82 195
  \pinlabel $V_m$ at 104 215
  \pinlabel $V_m$ at 127 238
  \pinlabel $V_m$ at 155 261
  \pinlabel $V_m$ at 220 256
  \pinlabel $V_m$ at 280 243
  \pinlabel $V_m$ at 250 223
  \pinlabel $V_m$ at 315 220
  \pinlabel $V_m$ at 285 191
  \pinlabel $V_m$ at 365 196
  \pinlabel $V_m$ at 335 135
  \pinlabel $V_m$ at 367 101
  \pinlabel $V_m$ at 455 51
  \pinlabel $V_m$ at 320 35
  \pinlabel $V_m$ at 275 40
  \pinlabel $V_m$ at 221 47
  \pinlabel $V_m$ at 167 54
  \pinlabel $V_m$ at 113 63
  \pinlabel $V_m$ at 57 70
\endlabellist
\centering
\includegraphics[height=2.2in]{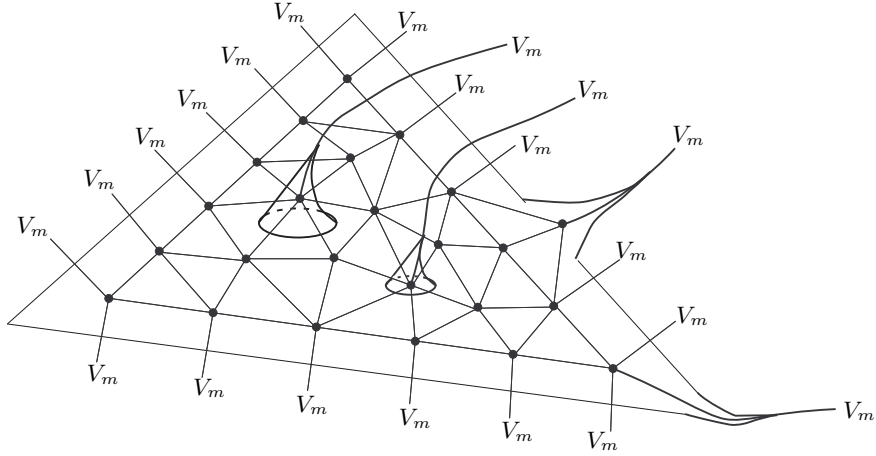}
\caption{The tensor network $N$: one leg exits every cusp and cusp
fraction, several legs exit every face.} \label{voronoi}
\end{figure}

To get the symmetries of $T_X$ right, we initially build a tensor
$T_W$ for $W$ (see Figure \ref{voronoi}) with indices or ``legs''
exiting cusps, cusp fractions, and codimension one faces of $W$.
$T_X$ is assembled from $|G|$ copies of $T_W$ (or $\bbar{T}_W$ for
orientation reversing elements of $G$) and additional $\tau$-tensors
described below.

We build a tensor $T_W$ for $W$ as follows:  Begin with a finite set
$v$ of vertices geometrically positioned in the interior of $W$ so
that every point of $x \in W$ either has injectivity radius$(x) \leq
\epsilon$ or $dist(x,v) < \epsilon$.  Consider the graph $\dot{N}$
consisting of the dual $1$-cells to the faces of the Voronoi
decomposition associated to $v$.  $\dot{N}$ is a fine triangulation
of thick($W$) with (many) edges orthogonal to each face of $W$.  Any
vertex $p \in v$ whose Voronoi cell contains a cusp or cusp fraction
is supplied an additional edge exiting that cusp or cusp fraction.
We turn $\dot{N}$ into a tensor network, now denoted $N$ by choosing
the vector space $V$ for each edge and a distinct generic tensor
$T_p$ for each vertex $p \in v$.  We take all the entries of all the
$T_p$ to be algebraically independent.

\begin{figure}[htpb]
\labellist
\small\hair 2pt
\pinlabel $W$ at 385 175
\endlabellist
\centering
\includegraphics[height=1.8in]{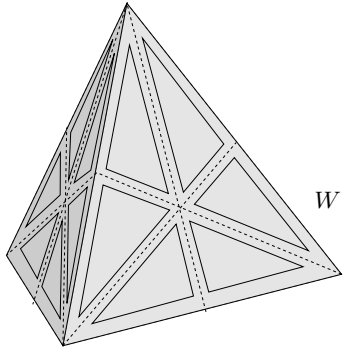}
\caption{$24$ copies of $W$ are glued together by
identifying $72$ faces in pairs. $l$ copies of $V_m$ are contracted along
each such face.} \label{txasm}
\end{figure}

For $\epsilon$ sufficiently small, the labeled network $N$ recovers
the combinatorics of $W$ with each face($W$) marked by numbered legs
from unique tensors $T_p$.  Now let $G$ label $|G|$ copies of $N$.
Actually if $q \in G$ reverses orientation, we conjugate the entries
of each $T_p$ to get $\bbar{N}$.  These copies of $N$ and $\bbar{N}$
may be assembled into a larger tensor network $GN$.  There is now a
tensor $\til{T}_X$ associated to $X$ by contracting the internal
edges (those not exiting cusps or cusp fractions) of $GN$.  It is
not difficult to see that $\til{T}_X$ has precisely the same
symmetry group (for small $\epsilon$) as $X$.  In fact, such a
labeled network encodes the homeomorphism type of the manifold (or
orbifold) $X$ on which it lies.  Also very important to our approach
is that given a reflection symmetry $r \in G$, $\til{T}_X$ may be
written as a contraction $\til{T}_X^{top} \til{T}_X^{bot}$ along the
legs orthogonal to the fixed set of $r$.  This ``decomposability''
of $T_X$ allows us to write part of $\til{T}_X$ in the bra and part
in the ket when pairing to the AMT.  Consult Figure \ref{voronoi}
for a visualization of $T_W$ as a contraction of $N$ and Figure
\ref{txasm} for a visualization of $|G|=24$ copies of $N$ assembled
into $GN$.

To repeat, $\widetilde{T}_X$ is assembled from $|G|$ copies of $T_W$
and its complex conjugate $\overline{T}_W$ by contraction along the
face-crossing indices. Figure~\ref{txasm} illustrates how $24$
copies of $W$ are assembled to make $X$, in a special case where $G$
is equal to the symmetric group $S_4$. $l$ legs cross each
codimension-1 face and additional indices, not indicated in the
figure, exit each cusp and cusp fraction.

In more detail, label the copies of $W$ making up $X$ as $\{gW | g
\in G\}$. $G$ acts on the indices in the tensor associated to cusps
or fractional cusps in the obvious way. That is, $H_1(C; \Z_m)$
serves separately as a basis for {\it each} fraction of the cusp
$C$, and $g$ acts as $g_*$. For instance, if we fix a basis $\{e_1,
\dots, e_{m^2}\} \times \{1,2,\dots,m\}$ for cusps and fractional
cusps in the copy labeled $(\id) W$, then the corresponding indices
in the copy labeled $gW$ are the basis elements $\{g_*e_1, \dots,
g_*e_{m^2}\} \times \{1,2,\dots,m\}$.  (Note the inaction on the
second factor.) On each copy $gW$ we place a copy of either $T_W$ or
$\overline{T}_W$ depending on whether $g$ preserves or reverses
orientation.

Our preliminary tensor $\widetilde{T}_X$ associated to $X$ is the
contraction of these $|G|$ copies of $T_W$ (or $\overline{T}_W$)
along the $l$ legs spanning the codimension-1 faces of $\{gW\}$. For
every reflection symmetry $r$ of $X$ with codimension one fixed set
$\Sigma$, the tensor $\widetilde{T}_X$ contracted with positive
vector inputs for each horizontal (with respect to $r$) cusp or cusp
fraction can be written in the form
$$\widetilde{T}_X(v_{2k+1},\cdots,v_n) = \widetilde{T}_Y^\dagger
\widetilde{T}_Y$$
where $\widetilde{T}_Y$ is the operator corresponding to the ``lower half'' of $X$
(again with respect to $r$), and therefore is non-negative Hermitian.

The tensor $\widetilde{T}_X$ is not yet the tensor we are looking
for because $T_W$ (and similarly $\overline{T}_W$) has one index
(ranging in $V_m$) for each {\em fractional cusp} in $W$, whereas we
want our tensor $T_X$ to have one index for each {\em cusp} of $X$.
A cusp $C$ stabilized by $H < G$ will contribute $|H|$ indices to
$\widetilde{T}_X$ whereas {\it one} is the desired number. To create
a tensor with the correct index set, we need to contract
$\widetilde{T}_X$ with suitably chosen tensors $\tau_{H(C)}$ for
each cusp $C$, where $H < G$ is the stabilizer.

\begin{figure}[htpb]
\labellist
\small\hair 2pt
\pinlabel $\widetilde{T}_X$ at 220 -10
\pinlabel $\tau_{H(C_1)}$ at 60 -8
\pinlabel $\tau_{H(C_2)}$ at 375 80
\pinlabel $\text{cusp }C_1\text{ with }|H|=4$ at 110 70
\pinlabel $\text{cusp }C_2\text{ with }|H|=6$ at 390 25
\endlabellist
\centering
\includegraphics[scale=0.5]{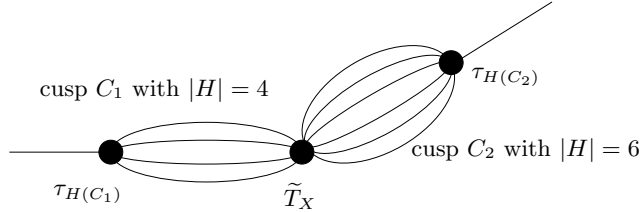}
\caption{The tensor $T_X$ is obtained by contracting $\tau_{H(C_1)}$ and $\tau_{H(C_2)}$
with $\widetilde{T}_X$} \label{figure8}
\end{figure}

We define $\tau_{H(C)}$ to be an element of $\Hom(V_m,V_m^{\otimes
|H|})$ satisfying the {\em reflection identity property} (RIP),
explained below, and so that the tensor $T_X$ obtained by
contracting $\widetilde{T}_X$ with a copy of $\tau_{H(C)}$ for each
cusp $C$ satisfies:
\begin{enumerate}
\item $\tau_{H(C)}$ is natural, so that $T_X$ may be constructed from the homeomorphism
type of $X$ alone, without resort to any additional data such as a base
point on the cusps.
\item $\tau_{H(C)}$ makes $T_X$ $G$-invariant.
(Recall that $G$ acts by permuting cusps internally on $V$.)
Succinctly one may say that $G$ acts via the action on $H_1(\partial
X; \Z_m)$. Actually, this is a corollary of item (1) above.
\end{enumerate}

Conditions 1 (and 2) hold provided $\tau_{H(C)}$ satisfies (a)
invariance under the action of $H$ on the $V_m^{\otimes |H|}$ factor
and (b) naturality, in the sense that $\tau_{H(C)}$ is defined using
only the action of $H$ on $C$ without additional choices such as a
base point on $C$.

All the $V$ factors in $\Hom(V_m, V_m^{\otimes |H|})$ are
canonically identified. If we fix a base point $p$ in $C$ marking
one fundamental domain for the $H$ action, then the formula
\[
    \tau^p_{H(C)}(e_i) = \bigotimes_{h \in H} h_*(e_i) , \quad\quad e_i=(f_j,s_k)\in H_1(C;
    \Z_m) \times \{1,2,\dots,m\}
\]
clearly satisfies condition (a).
To achieve (b), set
\[
    \tau_{H(C)}(e_i) = \sum_{h\in H} \tau^{h(p)}_{H(C)}(e_i) .
\]

\begin{figure}[htpb]
\labellist
\small\hair 2pt
\pinlabel $e_i$ at -30 133
\pinlabel $\tau^p_H$ at 100 190
\pinlabel $=$ at 320 133
\pinlabel $h_*e_i$ at 450 220
\pinlabel $h'_*e_i$ at 480 165
\pinlabel $h''_*e_i$ at 480 110
\pinlabel $h'''_*e_i$ at 450 55
\pinlabel $r$ at 880 133
\endlabellist
\centering
\includegraphics[scale=0.3]{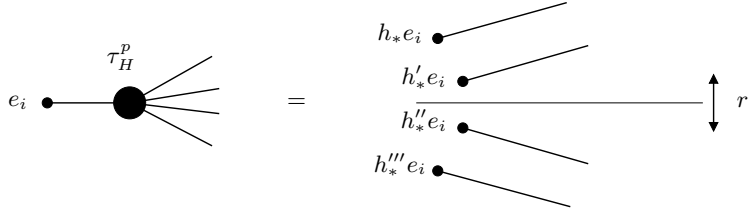}
\caption{The definition of $\tau^p_H$. Note that $rh=h'''$ and $rh'=h''$.} \label{taufig}
\end{figure}

Note that ``inputting'' vectors into tensors is a special case of
tensor contraction --- we regard elements of $V$ as tensors with one
index. This motivates the notation in Figure \ref{taufig}: the small
dots are 1-tensors. As illustrated in Figure \ref{taufig}, the
vector inputs induced by $e_{i(C)}$ are $r$-symmetric.

Specifically, if a basis element $e_i$ is inserted into the first
index of $\tau_H$ then the resulting $|H|$-tensor is converted to an
operator $\mathcal{O}$ by segregating the remaining $|H|$ index
slots into $r$-upper and $r$-lower groups with respect to a
reflection symmetry $r$.  RIP states that $\mathcal{O}$ is the
identity operator from the appropriate ``upper space'' $\mathcal{U}
\cong \C^{h/2}$ to the appropriate ``lower space'' $\mathcal{L}
\cong \C^{h/2}$, $h = |H|$, with respect to prescribed bases. Fixing
a base point fundamental domain $W$ on the $r$-upper side of the
cusp $C$, let $W = W_{id}$, $W_{g_2},\dots, W_{g_{h/2}}$ be the
domains on the upper side.  Then $$basis(\mathcal{U}) = \{h_*e_i
\otimes h_*g_{2*}(e_i) \otimes \cdots \otimes h_*g_{h/2*}(e_i) | h
\in H\}$$ and similarly
$$basis(\mathcal{L}) = \{rh_*e_i \otimes rh_*g_{2*}(e_i) \otimes \cdots
\otimes rh_*g_{h/2*}(e_i) | h \in H\}$$ with respect to domains
$W_r$, $W_{rg_2}, \dots, W_{rg_{h/2}}$ on the lower side.  Figure
\ref{taufig} illustrates this symmetry.

Our final tensor $T_X$ for hyperbolic JSJ pieces is obtained from
$\til{T}_X$ by contracting a copy of $\tau_H$ at each cusp $C$ (with
symmetry group $H$).

\begin{figure}[htpb]
\labellist
  \small\hair 2pt
  \pinlabel $\text{upper hyp.}$ at 132 144
  \pinlabel $\text{lower hyp.}$ at 132 58
  \pinlabel $\text{l.h.}$ at 45 58
  \pinlabel $\text{u.h.}$ at 45 144
  \pinlabel $\tau$ at 57 100
  \pinlabel $\tau$ at 90 100
  \pinlabel $\tau$ at 185 100
  \pinlabel $\text{SF}$ at 268 100
  \pinlabel $\text{SF}$ at 370 100
  \pinlabel $\text{SF}$ at 350 164
  \pinlabel $\text{SF}$ at 350 36
  \pinlabel $\text{hyp}$ at 245 174
  \pinlabel $\text{hyp}$ at 245 28
  \pinlabel $\text{hyp}$ at 412 176
  \pinlabel $\text{hyp}$ at 412 26
  \pinlabel $H^m$ at 442 100
  \pinlabel $r$ at 467 100
  \pinlabel $\text{vertical}$ at 485 174
  \pinlabel $\text{cusp}$ at 489 162
\endlabellist
\centering
\includegraphics[scale=0.8]{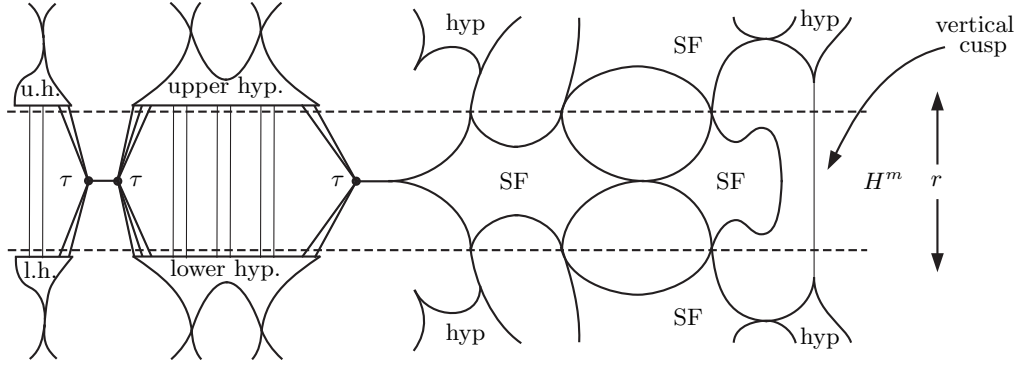}
\caption{The AMT operator $H^m$} \label{amtopfig}
\end{figure}

Having constructed $T_X$ for $X$ SF and $\tau$, these pieces can now
be assembled (together with identity operators on the legs of $T_W$
(and $\bbar{T}_W$) orthogonal to $fix(r)$ and also on vertical
cusps) into the AMT operator $H^m$.  Figure \ref{amtopfig} below
organizes the definition of $H$ into a picture.

In more detail, $H^m$ is obtained by contracting geometrically
adjacent tensors at the AMT.  These include $T_X$ for SF $X$, a
$\tau_H$ for each horizontal (hyperbolic) cusp, and then (without
any actual contraction) all the legs (identity operators on copies
of $V$) which run orthogonally to $S \cup (\text{hyperbolic JSJ
piece})$ and $id_V$ for vertical cusps.  The $\tau_H$ will have its
first index (slot) contracted to the first index (slot) of another
$\tau_{H'}$ when two hyperbolic pieces are glued at $S$. Where a
hyperbolic piece is glued to a SF $X$, $\tau_H$ will have its
first index (slot) contracted with the ``cusp'' of the SF $X$.

We now identify the subspace where $H^m$ acts as a positive
operator.  For each mid-level SF $X$, let $S_{\text{SF } X}^\pi$
denote the reduced symmetric product of the upper ``cusps'' of $X$.
For each mid-level hyperbolic $X$, let $\Z_{X,r}$ denote the tensor
product of the (space $V$ associated to each of the) legs of $GN$
orthogonal to $fix(r) \cap X = S \cap X$.  For each horizontal
hyperbolic cusp $C$, let $\mathcal{U}_C$ be the ``upper space''
defined above. Let:
\[
Y_{up} = \left(\bigotimes_{\substack{\text{mid level} \\ \text{SF }
X}} S_{\text{SF } X}^\pi \right) \otimes \left(
\bigotimes_{\substack{\text{mid level} \\ \text{hyperbolic } X}}
\Z_{X,r} \right) \otimes \left( \bigotimes_{\substack{\text{horizontal} \\
\text{cusps } C}} \mathcal{U}_C \right) \otimes \left(
\bigotimes_{\substack{\text{vertical}
\\ \text{cusps } C}} V_C \right)
\]
Make the analogous definition for $Y_{down}$, replacing
$\mathcal{U}_C$ with $\mathcal{L}_C$.  Perfect $r$-symmetry of the
GMT is obvious for $AA$ and $BB$ and will be established near the
end of the proof for $AB$ (unless $c(AB) < \max{c(AA),c(BB)}$).
Assuming such symmetry of the GMT, the reflection $r$ identifies
$Y_{up}$ and $Y_{down}$, so simply call them both $Y$.

\begin{lemma} \label{GMT_sym_lemma}
If the GMT is $r$-symmetric then $H^m$ is positive on $Y$ for all
$m=2,3,\dots$.
\end{lemma}

\begin{proof}
The proof is extremely simple.  $H^m$ is a sum over contraction
variables $e_i$ of tensor products of operators which are either
strictly positive (on $S_{\text{SF } X}^\pi$) or the identity (after
our identification of $Y_{up}$ and $Y_{down}$).  Positivity is
closed under both tensor product and convex combination.
\end{proof}

{\em Note:}  If $S_{\text{SF } X}^\pi$ is replaced by the ordinary
symmetric power $S_{\text{SF } X}$ we can only conclude {\em block}
positivity above.  This point is addressed in the final paragraph of
section 6.

\vspace{12pt}

We are almost in a position to define the complexity function $c_t$.
Suitable tensors $T_X$ for JSJ pieces $X$ have been constructed when
$X$ is hyperbolic. When $X$ is Seifert fibered, a final tweak is
required because of block positivity of these tensors. This is the
subject of the next subsection.

\subsection{Better tensors for Seifert fibered JSJ pieces}\label{better_tensor_subsection}

As we noted earlier, the tensors $\widehat{T}_X$ for Seifert fibered
JSJ pieces constructed above possess too much symmetry --- they are
insensitive to arcs of Dehn twists. In this subsection we repair
this defect at the expense of making the tensor $T_X$ for a Seifert
fibered JSJ piece slightly less local. $T_X$ will depend on the way
in which $X$ meets its neighboring JSJ pieces, not merely on the
homeomorphism type of $X$.

First we define the notion of {\em fiber slip homeomorphism}, or
{\em fish} for short. A Seifert fibered JSJ piece $X$ with non-empty
boundary will (by Theorem~\ref{exceptional_list_theorem}) have a
unique Seifert fibered structure unless $X \cong S^1
\widetilde{\times} S^1 \widetilde{\times} I$ in which case the
quotient orbifold is either a disk with two orbifold points of order
$2$, or a M\"obius band. To avoid ambiguity, by convention we always
give $S^1 \tilde{\times} S^1 \tilde{\times} I$ the Seifert fibered
structure for which the base orbifold is a disk with two order $2$
points.

\begin{definition} \label{fiber_slip_defn}
For a Seifert fibered JSJ piece $X \subset A$ (respectively $X
\subset M$), a {\em fiber slip} is a regluing of $A$ (resp. $M$) by
a fiber preserving homeomorphism of one or more cusps of $X$.  Note
that if two cusps of $X$ glue to each other so as to mismatch a
fiber, then there is a $\Z \oplus \Z$-family of fiber slips at the
common torus.
\end{definition}

\begin{definition} \label{fish_equiv_defn}
$(A,S)$ and $(B,S)$ (respectively $M$ and $M'$) are {\em fish
equivalent} if after regluing the Seifert fibered JSJ pieces of
$(A,S)$ (resp. $M$) by fiber slips, the result is homeomorphic to
$(B,S)$ rel $S$ (resp. homeomorphic to $M'$).
\end{definition}

To resolve the issue of fiber slips, we modify $c_t$, the partition
function induced by our preliminary choice of tensors for JSJ
pieces, into a more powerful tensor complexity. The tensors
$\widehat{T}_X$ previously assigned to a Seifert fibered JSJ piece
$X$ will now be perturbed, and in some cases re-scaled, in a way
which depends on the immediate neighbors of $X$ in the JSJ
decomposition. The perturbation $T_X$ will depend on the
homeomorphism types of the neighboring $(X',C')$ where $X'$ is a JSJ
piece glued along $C'$ to $(X,C)$ and also on the way in which $C$
is glued to $C'$. Note that $X'$ may equal $X$.

Consider a closed irreducible $3$-manifold $P$ and its JSJ decomposition.
Corresponding to the combinatorics of the decomposition is a gluing
graph $\mathcal{G}$ whose vertices are labeled by diffeomorphism
types of JSJ pieces and whose edges correspond to gluings of cusps.

\begin{lemma}
\label{lemma***1} For each Seifert fibered piece $X_i$ of the JSJ
decomposition, let $Q_i$ denote the underlying orbifold. A fiber
slip homeomorphism of $P$ which covers the identity on each $Q_i$ is
determined by a single integer $t_i$ for each $X_i$.
\end{lemma}

\begin{proof}
Fiber slips are of infinite order and hence cannot extend as
symmetries of an adjacent hyperbolic piece.  By the definition of
the JSJ decomposition, when neighboring (or self-\!\! neighboring)
Seifert fibered JSJ pieces meet along a cusp $C$, their fibers do
not match in $C$.  As a result, a fiber slip on one side will not
extend to {\it any} diffeomorphism on the other.  As a consequence,
a fiber slip is determined by what it does to each Seifert fibered
JSJ piece relative to its cusps. After reparameterizing a Seifert
fibered piece $X_i$ by arcs of Dehn twist, the effect of the slip
may be localized to a single cusp, where it differs from the
identity by some power of a Dehn twist along a fiber. Note that only
an orientation of $X_i$, not $Q_i$, is required to identify the sign
of $t_i$.
\end{proof}

\begin{remark}
If we think of a Seifert fibered space as a circle bundle over an
orbifold, a choice of trivialization of this bundle over each
boundary component determines a relative Euler class, which can be
paired with the fundamental class of the base orbifold to get a
number. Two Seifert fibered spaces which differ by a fiber slip
differ (after choosing trivializations on their boundaries) only in
the value of this number.
\end{remark}

It follows that the set of ways of gluing $X_i$ to its neighbors
which are in a fixed fish equivalence class is an oriented
$\Z$-torsor, which we denote $\tau_i$. It turns out that it is
(almost) possible to choose a canonical basepoint for this torsor.
Either one element can be unambiguously labeled $0$, or else two
consecutive elements can be labeled as $\{-\frac{1}{2},
+\frac{1}{2}\}$.

Let $\{C\}$ denote the set of cusps of $X$, a Seifert fibered JSJ
piece of the irreducible $3$-manifold $P$, and assume that $X$ has a
fixed Seifert fibered structure. We need to come up with a set
\{bases\} of canonical base curves, one for each $C$ which pairs
nontrivially with the fiber in $C$. In fact, the ``set'' \{bases\}
will consist, not of a single set, but of an orbit of sets under the
action of the mapping class group of $X$.  In particular, an ADT
acts on the set of choices \{bases\} of canonical base curves. The
elements in \{bases\} are determined homologically. Assume $X$ has
$q$ boundary components. Let $K$ be the kernel of the map
$H_1(\partial X ; \Q) \rightarrow H_1(X;\Q)$ induced by inclusion.
Then $K$ has dimension $q$. 
Even if base and fiber of $X$ are not
separately orientable, we can, and do, consistently orient the
fibers on $\partial X = \bigcup \{C\}$. 
For each $i \leq q$, let
$f_i$ denote the class  in $H_1(\partial X;\Q)$ of the oriented fiber
of the $i$th boundary component of $X$. 
Each difference $f_i - f_{i+1}$ is in $K$, and
together these differences span a $q-1$ dimensional subspace of $K$.
Let $b \in H_2(X, \partial X ; \Q)$ be any class such that $\partial
b$, together with the $f_i - f_{i+1}$, forms a basis for $K$, and
scale $b$ to be primitive in $H_2(X, \partial X ; \Z)$. 
$\partial_i b$ may be divisible in $H_1(\partial_i X; \Z)$, so define
$b_i \in H_1(\partial_i X; \Z)$ to be a primitive fraction of $\partial_i b$.
Take $\{b_i\}$ to be the desired set of basic curves. 
For each cusp $C_i$, 
$b_i\cap f_i = r_i \ne 0$.
Replacing $b_i$ by $-b_i$ changes $r_i$ to $-r_i$; adding a multiple of 
the preimage of
some $f_j - f_{j+1}$ to $b$ obviously leaves $r_i$ fixed. Once we fix an
orientation on $f$ and $\partial X$, we can resolve the ambiguity in
the signs of the $r_i$'s by possibly changing the sign of $b_i$ so
that each $r_i$ is positive.
This makes each $b_i$ well defined.

The set of isotopy classes of simple essential curves in $C$ whose
intersection number with the oriented fiber $f$ is a fixed $r>0$ is a
$\Z$-torsor $\tau_C$. 
An orientation-preserving homeomorphism $\alpha$ of $X$ either preserves or reverses
the class of the fiber $f$.
Since the preceding construction is homological, if $\alpha(C_i) = C_j$ then
$\alpha(b_i) = \pm b_j$.
Since $\alpha$ is orientation-preserving, either $\alpha(f) = f$ and 
$\alpha(b_i) = b_j$ or $\alpha(f) = -f$ and 
$\alpha(b_i) = -b_j$.
If a cusp $C$ is fixed by $\alpha$ then $\alpha(b+kf) = \pm(b+kf)$
in the two cases respectively.
Hence the action of $\Z$
on the torsor $\tau_C$ commutes with the action of the mapping class
group of $X$.

We would like to identify a canonical basepoint (i.e.\ $0$) in the
torsor $\tau_C$.  We almost succeed; either we find a unique element
which we call $0$, or we find two elements which differ by $1$ which
we call $-\frac 1 2$ and $\frac 1 2$ (consistent with the action of
$+1$). We call this {\em marking the torsor $\tau_C$}.

To mark $\tau_C$, look at $C'$, the cusp to which it is glued in the
JSJ decomposition. Here $C'$ is a cusp of another JSJ piece $X'$,
where $X=X'$ is possible. If $C'$ lies in a Seifert fibered JSJ
piece, then consider the fiber $f'$ in $C'$ and how it intersects
with the fiber $f$ of $C$ under the gluing. By the definition of the
JSJ decomposition, $f$ and $f'$ cannot be isotopic.  Let $b$ be a
``base'' class of $C$ as constructed above. The function
$$g_{SF}(k) = |(b+kf) \cap f'|$$
admits a unique minimum for $k \in \Q$, and therefore achieves its
minimal value for $k \in \Z$ on at most two values of $k$, which
differ by $\pm1$. We mark $\tau_C$ by $\{b + k_0 f \; | \;
g_{SF}(k_0) \text{ is minimal}\}$.

If $C'$ lies in a hyperbolic piece, consider
$$g_{hyp}(k) = \frac{\length^2_{C'} (b+kf)}{\area(C')}$$
where lengths and areas are computed with respect to some fixed
Euclidean structure on $C'$ in the similarity class induced by the
hyperbolic structure on $X'$. Since $\length^2/\area$ is
dimensionless, the result only depends on the similarity class. The
quadratic form $g_{hyp}(k)$ is strictly convex for $k \in \R$ and
therefore has either one or two minima $k_0 \in \Z$. Again, mark
$\tau_C$ by $\{b + k_0 f \; | \; g_{hyp}(k_0)$ is minimal$\}$.

The previous discussion allows us to define an invariant of an
oriented Seifert fibered JSJ piece $X$ in an oriented, irreducible
$3$-manifold $P$, called the {\em gluing number}, and denoted
$\gn(P,X)$. The gluing number takes values in $\Z[\frac 1 2]$, and
is computed as follows. The marking identifies each $\tau_C$ with a
copy of $\Z$ or $\frac 1 2 + \Z$. Therefore the set of base classes
in each $\tau_C$ can be thought of as a set of integers and
half-integers, whose sum is $\gn(P,X)$. Note that $\gn$ does not
change under an orientation-preserving homeomorphism which changes
$f$ and $b$ to $-f$ and $-b$. Moreover, adding an element $f_i -
f_{i+1}$ to $b$ adds $1$ to one term in the sum and $-1$ to the
other, and therefore does not change $\gn$. Hence $\gn$ is
well-defined, independent of choices. Furthermore, a fiber slip
alters $\gn$ by adding exactly the integer defined in
Lemma~\ref{lemma***1}.

We use the notation $\gn(\til{X},X)$ interchangeably with
$\gn(P,X)$, where $\til{X}$ denotes the union of $X$ with its
immediate neighbors in the irreducible prime $P$. This notation is
convenient, since it underscores the fact that $\gn$ depends on
$\til{X}$ (and not just $X$) but not necessarily on all of $P$.

\begin{figure}[htpb]
\labellist
\small\hair 2pt
\pinlabel $g_{SF}(k)$ at 145 270
\pinlabel $\text{for }X$ at 145 250
\pinlabel $g_{SF}(k)$ at 505 270
\pinlabel $\text{for }\overline{X}$ at 505 250
\pinlabel $-3$ at 37 40
\pinlabel $-2$ at 73 40
\pinlabel $-1$ at 109 40
\pinlabel $0$ at 145 40
\pinlabel $1$ at 181 40
\pinlabel $2$ at 217 40
\pinlabel $3$ at 253 40
\pinlabel $-3$ at 397 40
\pinlabel $-2$ at 433 40
\pinlabel $-1$ at 469 40
\pinlabel $0$ at 505 40
\pinlabel $1$ at 541 40
\pinlabel $2$ at 577 40
\pinlabel $3$ at 613 40

\endlabellist
\centering
\includegraphics[height=2in]{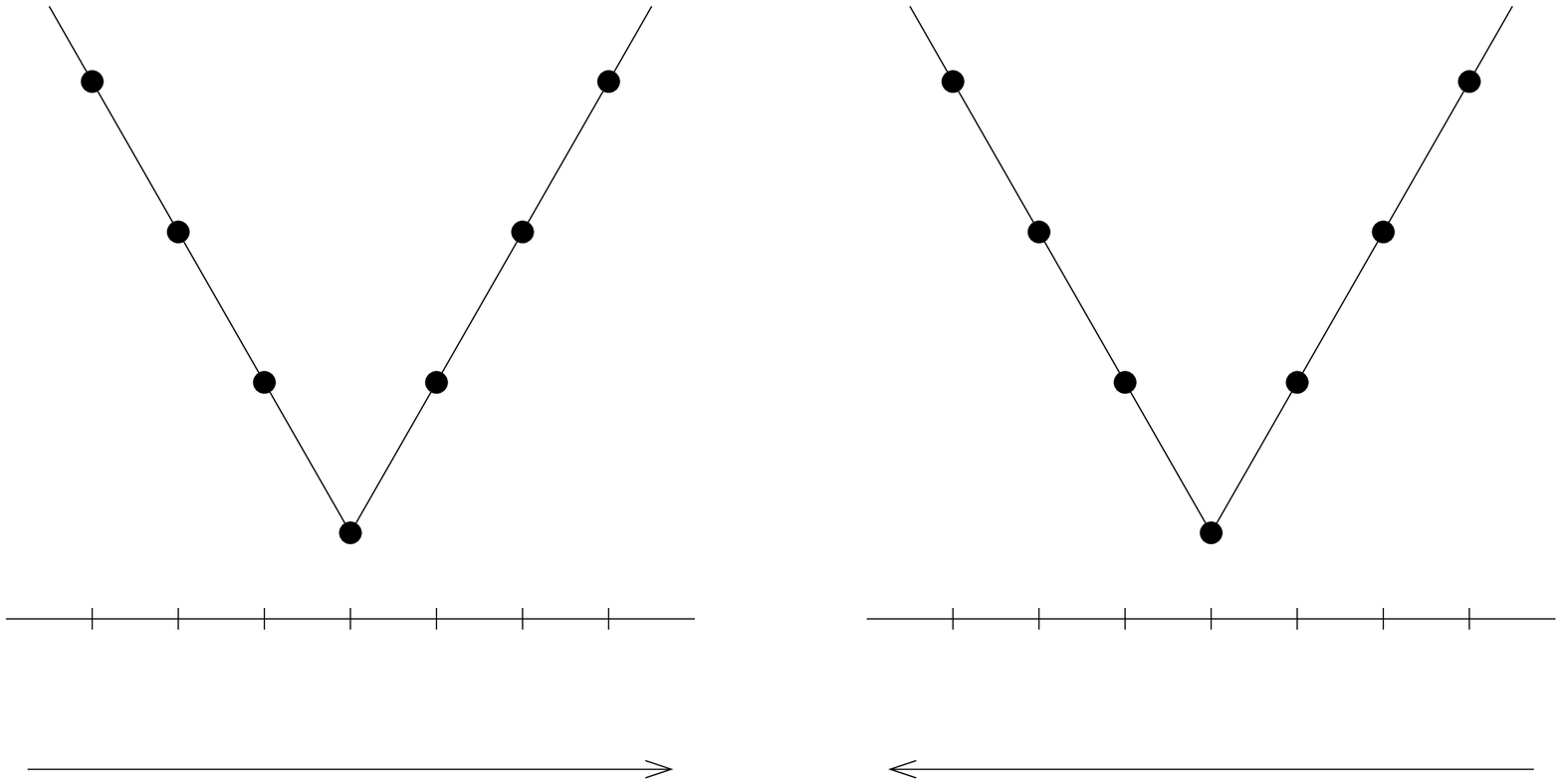}
\caption{} \label{orientation_reversal_fig}
\end{figure}

\begin{lemma}\label{gn_homeo_extension_lemma}
If $X$ admits an orientation reversing homeomorphism which extends
to $\til{X}$, then $\gn(\til{X},X) = 0$.
\end{lemma}
\begin{proof}
We must show that $\gn$ is both chiral and an invariant. That is,
\begin{enumerate}
\item{If $h: (\til{X},X) \rightarrow (\til{Y},Y)$ is an orientation-preserving
homeomorphism, then $\gn(\til{X},X) = \gn(\til{Y},Y) \in \Z[\frac{1}{2}]$.}
\item{Reversing the orientation on $(\til{X},X)$ reverses the sign
of $\gn(\til{X},X)$.}
\end{enumerate}
The lemma follows from (1) and (2). We remark that this is somewhat
reminiscent of the classical observation that the signature of a
$4k$-dimensional manifold must vanish if it admits an orientation
reversing homeomorphism.

Property~(1) is more or less proved already, and is just a restatement
of the fact that $\gn$ is well-defined. Property~(2) is more interesting.  The marking of the
$\Z[\frac{1}{2}]$-torsor associated to a cusp $C$ of $X$ depends on
whether the adjacent JSJ piece $X'$ is Seifert fibered or hyperbolic,
and it depends on the unoriented class of $f'$ in $C'$ in the first
case, and the geometry of $C'$ in the second case. The key point
is that in neither case do the functions $g_{SF}$ and $g_{hyp}$
depend on the {\em orientation} of $X,C'$.  Reversing the
orientation of $X$ reverses the action of $\Z$ on $\tau_C$ but it does not
change the basepoint(s), so it just reverses the sign of the invariant;
see Figure~\ref{orientation_reversal_fig}.
\end{proof}

We are now almost ready to modify the tensors associated to
Seifert fibered JSJ pieces. First we define the property of {\em domination}
for Hermitian forms.
\begin{definition}
\label{dominates_defn} Let $\Oh, \Oh'$ be two Hermitian
forms on the same vector space. We say $\Oh$ {\em dominates} $\Oh'$ if
$$\langle x, \Oh x \rangle > \langle x, \Oh' x \rangle$$
for all nonzero $x$. Equivalently, $\Oh$ dominates $\Oh'$ if $\Oh - \Oh'$ is positive.
\end{definition}

Now we can state the properties that we want our modified tensors to have.
\begin{definition}\label{suitable_definition}
For $X$ Seifert fibered, a tensor $T_{X,\gn}$ is {\em suitable} if it
satisfies the following properties:
\begin{enumerate}
\item It has the same symmetry and positivity properties as $T_X$, and (compare
equation~\ref{eqn4}) if $\overline{X}$
is the orientation reverse of $X$, then $T_{\overline{X}} =
\overline{T}_X$.
\item The collection is otherwise generic, as was $\{\widehat{T}_X\}$.  In particular,
$T_{X,n} \neq T_{X,m}$ for $n \neq m \in \Z[\frac{1}{2}]$.
\item When $\gn$ assumes the value $0 \in \Z$ or $\Z[\frac{1}{2}]$,
then $T_{X,0} > T_{X,n}$, for $n \neq 0$, in the sense that for
any reflection symmetry $r$ and for any given index inputs for the
horizontal cusps, the contracted operator $\Oh_{X,0}$ dominates
the contracted operator $\Oh_{X,n}$ in the sense of Definition~\ref{dominates_defn}.
\end{enumerate}
\end{definition}

A simple way to satisfy condition (3) (as well as (1) and (2)), that
is, to promote $\widehat{T}_X$ to $T_X$, is to choose a function
$\theta : \Z[\frac{1}{2}] \rightarrow \R^+$ with a single maximum at
$0 \in \Z[\frac{1}{2}]$. First define $\check{T}_{X,n} = \theta(n)
\widehat{T}_{X}$ and then obtain $T_{X,n}$ from $\check{T}_{X,n}$ by
a fully generic (but tiny) perturbation of the entries subject only
to the constraint that it should be block positive.  That is, the
entries (for each fixed $n$) depend only on the images of the
indices under the quotient map $\pi$. So suitable tensors exist.

In the next subsection, we will insure that for each of the
manifolds $AB$, $AA$ and $BB$, the ``middle third'' consists of the
same collection of JSJ pieces which have been glued together in ways
which differ potentially only by fiber slip homeomorphisms. This
combinatorial pattern of gluings (in each case) determines a
Hermitian form, obtained by suitably contracting various tensors
$T_{X_i,n_i}$ for Seifert fibered $X_i$ and $T_{X_i}$ for hyperbolic
$X_i$ defined as above. The next lemma compares the operators
obtained by different gluings.  We suppress the parameter $m$ as a
superscript to the operator $H$.

\begin{lemma}
\label{T_X,0_dominates_lemma} The Hermitian form $H$ obtained by
contracting the various $T_{X_i,0}$, $\tau$'s, and identity
operators in the AMT dominates all forms $H'$ obtained by replacing
at least one $T_{X_i,0}$ by $T_{X_i,n_i}$, $n_i \neq 0$, for some SF
$X_i$.
\end{lemma}
\begin{proof}
This follows formally from condition~(3) above, provided that the
factor $\displaystyle \theta(0) / \max_{j \neq 0} \theta(j)$ overwhelms
the product of the perturbations.  This may be ensured by choosing
perturbations by genuine infinitesimals (in the sense of
nonstandard analysis), or by extending the coefficient
ring from $\C$ to $\C [x]$, with lexicographical norm $(|a_0|,
|a_1|, \dots)$ on $\sum a_i x^i \in \C [x]$, and taking
perturbations only in the linear term (and fixing $\displaystyle \theta(0)
/ \max_{j \neq 0} \theta(j) > 1$).
\end{proof}

The importance of ``domination'' comes from the following lemma:

\begin{lemma}
\label{O_dominance_lemma} If $\Oh$ dominates $\Oh'$ then for all non-zero
$x,y$ there is an inequality
$$\langle x, \Oh' y \rangle < \max(\langle x, \Oh x \rangle , \langle y, \Oh y \rangle)$$
\end{lemma}

\begin{proof}
By the definition of domination, there are inequalities
$\langle x, \Oh x \rangle > \langle x, \Oh' x \rangle$ and
$\langle y, \Oh y \rangle > \langle y, \Oh' y \rangle$, so
$$\max(\langle x, \Oh x \rangle , \langle y, \Oh y \rangle) >
\max(\langle x, \Oh' x \rangle , \langle y, \Oh' y \rangle)$$
But the Cauchy-Schwarz inequality gives
$$\max(\langle x, \Oh' x \rangle , \langle y, \Oh' y \rangle) \ge \langle x, \Oh' y \rangle$$
and we are done.
\end{proof}

\subsection{Assembling the assembly complexity} \label{acasect}

The assembly complexity is a lexicographic pair $(c_r, c_t)$, where:
\begin{itemize}
\item $c_r$ is a simple symmetry based complexity, described below.
Its job is to ensure the the middle thirds of $AA$, $BB$ and $AB$
are homeomorphic up to fiber slips.
\item $c_t$ is an infinite lexicographic tuple $(c_t^m)$, $m \ge 2$.
\item $c_t^m$ is a number obtained by contracting configurations of
tensors of the kind described in \S\ref{actjsjsect} and \ref{better_tensor_subsection}.
For each Seifert fibered JSJ piece $X$ the tensor $T_{X,\gn}$ depends
on the way in which $X$ is glued to its neighboring JSJ pieces, through
the invariant $\gn$. The
purpose of this modification is to correct for the failure of $\widehat{T}_X$
to detect fiber slips.
\end{itemize}

First we define $c_r$ and prove an accompanying lemma.
Let $M$ be a closed irreducible $3$-manifold and let $T$ be a JSJ
torus of $M$.
\begin{definition}\label{reflective_torus_definition}
A JSJ torus $T$ of a closed irreducible $3$-manifold $M$ is
{\em reflective} if there is an orientation reversing involution
of $T$ with two fixed circles, which extends to the JSJ pieces
on either side of $T$.
\end{definition}
Note that for the doubles $AA$ and $BB$, each JSJ torus transverse to $S$
(hereafter {\em vertical JSJ torus}) is reflective. The two fixed circles
of the reflective involution interchanging the two sides are exactly the two
components of $T \cap S$.

\begin{definition}
Let $M$ be a closed irreducible $3$-manifold. Define $c_r(M)$ to be the
number of reflective JSJ tori in $M$.
\end{definition}

\begin{lemma}[$c_r$-Lemma schema]\label{c_r_lemma_schema}
Let $AB$, $AA$ and $BB$ be closed, irreducible and connected, and
suppose that their $c_S$ and $c_h$ complexities are equal. Then
either $c_r(AB) < \max(c_r(AA), c_r(BB))$ or $c_r(AB) = c_r(AA) =
c_r(BB)$ and the geometric middle thirds (GMTs) of $AB$, $AA$ and
$BB$ are homeomorphic (rel boundary) up to fiber slip homeomorphism.
\end {lemma}

\begin{proof}
First we show that $c_r(AB) \le \max(c_r(AA), c_r(BB))$. Define
$c_r(A)$ to be the number of reflective JSJ tori in $A$ (reflective annuli
don't count). Define $c_r(B)$ similarly. For $M = AB, AA$ or $BB$,
let $c_m(M)$ be the number of reflective tori in the middle third.
Recall that these middle thirds are identical except for the gluing
maps at the vertical JSJ tori, which might differ by Dehn twists parallel to
$S$. We have
\begin{eqnarray*}
c_r(AB) &=& c_r(A) + c_r(B) + c_m(AB) \\
c_r(AA) &=& c_r(A) + c_r(A) + c_m(AA) \\
c_r(BB) &=& c_r(B) + c_r(B) + c_m(BB) .
\end{eqnarray*}
As noted above $c_m(AA)$ and $c_m(BB)$ are both equal to the total
number of vertical JSJ tori in the middle level. Thus
$c_m(AB) \le c_m(AA) = c_m(BB)$ and it follows that $c_r(AB) \le \max(c_r(AA), c_r(BB))$.

The equality $c_r(AB) = \max(c_r(AA), c_r(BB))$ can hold only if
both $c_r(A) = c_r(B)$ and $c_m(AB) = c_m(AA) = c_m(BB)$. The latter
equality means that every vertical JSJ torus in the middle third of
$AB$ is reflective. We will show that this in turn implies that the
middle third $AB$ differs from the middle third of $AA$ (or $BB$)
only by fiber slips.

Let $T$ be a vertical JSJ torus in the middle third of $AB$, and
denote the corresponding torus in $AA$ also by $T$. On each side of
$T$ one finds one of the following three kinds of pieces:
\begin{itemize}
\item a hyperbolic JSJ piece,
\item a Seifert fibered JSJ piece with fibers parallel to $S \cap T$
(horizontally Seifert fibered or HSF), or
\item a Seifert fibered JSJ piece with fibers intersecting each
component of $S \cap T$ once
(vertically Seifert fibered or VSF).
\end{itemize}
Note that the existence of the global reflection in $AA$ means that
the fibers of a Seifert fibered JSJ piece are either parallel to $S
\cap T$ or have intersection number one with each component. This
property is preserved under Dehn twists parallel to $S$, and
therefore also holds in $AB$.

Suppose a vertical JSJ torus $T$ in $AB$, $AA$ or $BB$ abuts two
Seifert fibered pieces $X,X'$. Let $f,f'$ be the class of the fibers
in $X,X'$ in the torus $T$. Since every vertical JSJ torus in each
space is reflective, there is an involution $\iota$ of $T$ which
extends to $X,X'$. By uniqueness of fiber structures, the involution
$\iota$ takes $f$ to $\pm f$ and $f'$ to $\pm f'$. Note that the
pieces $X,X'$ and whether they are HSF or VSF are the same in $AB$,
$AA$ and $BB$, though not necessarily the way they are glued. If
$X,X'$ are both VSF in $AA$ then the involution $\iota$ may be
chosen to fix $S$. This implies that $f = \pm f'$, which means that
$X$ and $X'$ can be glued along $T$ to make a bigger JSJ piece,
contrary to the defining properties of the JSJ decomposition. It
follows that a vertical JSJ torus can bound a VSF piece on at most
one side in $AA$, and therefore also in $AB$.

The gluing map at $T$ gives a bijection between two copies of
$\GL(2,\Z)$, one from each of the two adjacent JSJ pieces. Inside
each of these copies is the subset of reflections which extend over
the full JSJ piece. The fact that $T$ is reflective in $AA$ means
that the under the isomorphism of the two copies of $\GL(2,\Z)$
induced by the gluing map for $AA$, the two sets of reflections have
at least one element in common. We know that the gluing maps for
$AB$ along the middle third differ from those of $AA$ by Dehn twists
parallel to $T\cap S$. We can ask whether any of the homeomorphisms
(other than the gluing map for $AA$) in this affine $\Z$ coset give
a non-empty intersection of the two sets of reflections. We will see
below that if neither adjacent JSJ piece is HSF, then the answer is
no, so in this case the only way for $T$ to be reflective in $AB$ is
for the gluing map for $AB$ to coincide with the gluing map for
$AA$. On the other hand, if $T$ is adjacent to a HSF JSJ piece, then
all of the possible gluing maps at $T$ are related by fiber slips,
so the lemma will be proved.

Let $T$ be a JSJ torus in the middle third, with adjacent JSJ pieces
either hyperbolic or VSF. Fix coordinates so that
$$\left[ \begin{array}{rr} 1 & 0 \\ 0 & -1 \end{array} \right] \text{ denotes the canonical
reflection of } AA$$
and
$$\left[\begin{array}{rr} 1 & n \\ 0 & 1 \end{array} \right] \text{ denotes a product of }n\text{ Dehn
twists in a component of }T\cap S .$$ In these coordinates, the set
of reflections of $T$ which extend over an adjacent VSF piece must
be a subset of
\[
{\mathcal R}_v := \left\{\left[ \begin{array}{rr} \pm 1 & 0 \\ 2m &
\mp 1 \end{array} \right] , m \in \Z \right\},
\]
since such reflections must preserve the (vertical) Seifert fibers
and be conjugate to the standard reflection by Dehn twisting.
The set of reflections of $T$ which extend over an adjacent
hyperbolic piece must be a subset of
\[
    {\mathcal R}_h := \left\{
        \left[\begin{array}{rr} 1 & 0 \\ 0 & -1 \end{array} \right],
        \left[\begin{array}{rr} -1 & 0 \\ 0 & 1 \end{array} \right],
        \left[\begin{array}{rr} 0 & 1 \\ 1 & 0 \end{array} \right],
        \left[\begin{array}{rr} 0 & -1 \\ -1 & 0 \end{array} \right]
    \right\},
\]
since the geometry of the cusp must be rectangular, possibly square
(note that the last two matrices occur only in the case of a square cusp).

Let the gluing map at $T$ in $AB$ be the product of $n$ Dehn twists
parallel to $T\cap S$. Then we
must have $P \in {\mathcal R}_h$ and
$Q \in {\mathcal R}_h$ or $Q \in {\mathcal R}_v$ such that
\[
    \left[
    \begin{array}{rr} 1 & n \\ 0 & 1 \end{array} \right] P \left[
    \begin{array}{rr} 1 & -n \\ 0 & 1 \end{array} \right] = Q .
\]
It is straightforward to check that the only solution is $n=0$. Thus
the gluing maps at $T$ for $AA$ and $AB$ coincide. As remarked
above, if one side of $T$ is HSF, any gluing map is a fiber slip,
and the lemma is proved.
\end{proof}

We now show that the complexity $c_p: = (c_S,c_h,c_a)$, where $c_a =
(c_r,c_t)$, satisfies the conclusion of
Theorem~\ref{prime_complexity_theorem}. This is the final step in
the proof of the Topological Cauchy-Schwarz inequality in the guise
of Theorem~\ref{c3_theorem}, and therefore also of the proof of
Theorem~A.

\vskip 12pt

We assume without loss of generality that the terms $c_S,c_h,c_r$
are equal for $AA$, $AB$ and $BB$ and immediately focus on $c_t =
(c_t^m)$.

Contracting, as appropriate, copies of tensors $T_{X,gn}$ for $X$
SF, $T_W$ for $W \subset X_{\text{hyp}}$ a fundamental domain, the
various $\tau_{H(C)}$ tensors associated to horizontal hyperbolic
cusps, and the identity operator on $V$ for vertical cusps, we
obtain the three pieces of $c_t^m(AB)$ (respectively $c_t^m(AA)$ and
$c_t^m(BB)$).  Writing the formula, just for $AB$, we have:
$$
    c_t^m(AB) = \langle c_t^m(A^\sim) | H_{AB}^m | c_t^m(B^\sim) \rangle
$$
This is the tensor contraction corresponding to $A^\sim \cup
\text{AMT} \cup B^\sim$. Similar formulas hold for $c_t^m(AA)$ and
$c_t^m(BB)$. (Note also that our notation here differs slightly from
our notation for graphs where we wrote $G = A \cup H\cup B$. Here we
write $AB = A^\sim \cup \text{AMT} \cup B^\sim$.) If $c_r(AB) =
c_r(AA) = c_r(BB)$ then by Lemma~\ref{c_r_lemma_schema},
$\text{GMT}_{AB}$, $\text{GMT}_{AA}$ and $\text{GMT}_{BB}$ differ
only by fiber slips. Furthermore,
Lemma~\ref{gn_homeo_extension_lemma} shows that the $\gn$ invariant
of Seifert fibered pieces in $M_{AA}$ (and similarly in $M_{BB}$)
must vanish. Therefore by Lemma~\ref{T_X,0_dominates_lemma}, for
each $m$ the operator $c_t^m(M_{AA})$ dominates the operator
$c_t^m(M_{AB})$, unless $\text{GMT}_{AA} = \text{GMT}_{AB}$ (and
similarly for $\text{GMT}_{BB}$).

By Lemma~\ref{O_dominance_lemma} it follows that there is a {\em
strict} inequality of lexicographic strings
\begin{multline*}
(\langle c_t^m(A^\sim) | H_{AB}^m | c_t^m(B^\sim) \rangle) < \\
\max( (\langle c_t^m(A^\sim) | H_{AA}^m | c_t^m(A^\sim) \rangle),
(\langle c_t^m(B^\sim) | H_{BB}^m | c_t^m(B^\sim) \rangle) )
\end{multline*}
unless $\text{GMT}_{AB} = \text{GMT}_{AA} = \text{GMT}_{BB}$ and
$c_t^m(A^\sim) = c_t^m(B^\sim)$.

So to complete the proof we just need to check that $c_t^m(A^\sim) =
c_t^m(B^\sim)$ implies $A^\sim = B^\sim$ (rel their gluings to the
AMT). Of course, $A^\sim$ and $B^\sim$ are attached to different
sides of $\text{AMT}_{AA} = \text{AMT}_{AB}$; we really mean that
the pairs $(A^\sim,A^\sim \cup \text{AMT}_{AB})$ and $(B^\sim,
B^\sim \cup \text{AMT}_{AB})$ are diffeomorphic by a diffeomorphism
which restricts to the canonical involution on $\text{AMT}_{AB}$
inherited from the identification of $\text{AMT}_{AB}$ with
$\text{AMT}_{AA}$.

Recall that $H^m$ is positive on the subspace $Y$ (Lemma
\ref{GMT_sym_lemma}). Let $m$ be larger than the number of edges in
the fine scale tensor network in $A$ or $B$. Consider the most
general possible basis element $y$ belonging $Y$. $y$ is of course
symmetrized on the $S^\pi$ factors and also has a dictated form on
the $\mathcal{U}_C$ factors. This form will repeat a basis element
whenever $h$ and $h'$ belonging to $H(C)$ are homotopic near the
cusp. But as far as these constraints allow, $y$ should be a tensor
product of distinct basis elements (symmetrized as required in
$S^\pi$). Call these two features of $y$ ``deficiency 1''  (its
symmetrization) and ``deficiency 2'' (the repeated indices that may
occur in $\mathcal{U}_C$). We need to show that the $y$ components
of $c_t^m(A^\sim)$ and $c_t^m(B^\sim)$ are (formally/generically)
equal only if $A^\sim$ is equivalent to $B^\sim$ in the sense of the
preceding paragraph. We call these features of $y$ ``deficiencies''
because they complicate the most straight forward strategy of
producing a geometric isomorphism by first matching monomials of the
partition function.

Deficiency 2 is actually immaterial. To faithfully assemble copies
of the fundamental domain $W$ back into $X$ (or the $r$-upper/lower
half of $X$), we actually already have plenty of gluing information
given by the face crossing legs of the small scale tensors and do
not need to rely on distinctly labeled 1/2 legs of $\tau$ tensors at
the cups.

Regarding deficiency 1, the situation is quite similar to Lemma
\ref{grpos3}. Recall that the broad idea would be to look at a
monomial where all edge colors which occur twice are distinct from
each other and from the colorings occurring in $y$, and use it as a
recipe to reconstruct $A^\sim$ and equivalently $B^\sim$. But just
as in Lemma \ref{grpos3} we should not expect to reconstruct
$A^\sim$ ( or $B^\sim$) rel identity on its lower boundary but
rather should expect to encounter homeomorphisms of the lower
boundary which extend to admissible SF cusp permutations (which
extend to the AMT). Such permutations correspond (as in the proof of
Lemma \ref{grpos3}) to a corresponding permutation of the monomials.
So deficiency 1 is, for the most part, the familiar and expected
result of the geometric symmetry of the SF cusps of the AMT.
However, $S^\pi$ also projects out the fiber coordinate of each SF
cusp so with the original tensors (unmodified by the $gn$ invariant)
we could only hope to identify $A^\sim$ and $B^\sim$ up to fish
equivalence. This ambiguity is resolved in the next paragraph.

Away from the GMT we conclude that $A^-$ and $B^-$ are fiber slip
homeomorphic (rel their gluings to the GMT). But properties (1) and
(2) from Definition~\ref{suitable_definition} (especially property
(2)) show that the gluing invariants $t_i$ of a fiber slip
homeomorphism from $A^-$ to $B^-$ can be recovered from the
properties of $T_{X,\gn}$. This is the last piece of ambiguity, and
therefore we conclude $A^- = B^-$ rel their gluings to the GMT, and
therefore that $(A,S) = (B,S)$. This completes the proof of the main
theorem.

\section{Extensions of the Main Theorem}

In this section we discuss some natural variations on our main
positivity theorem, and sketch proofs.
We write ``theorem'' in lower case to acknowledge that important points
may have been overlooked.
We hope that one or more of our readers
will supply full proofs for these statements.

In a similar vein, although the proofs we give of graph positivity in \S\ref{acgcssect}
are complete, a broader formulation of the statements there would be desirable.
Graduate students take note.

\subsection{Unoriented case}

\begin{lctheorem} \label{M_S^u_positivity_theorem}
Given $S$, an unoriented, compact (possibly disconnected) surface
without boundary, let $\Mdot_S^u$ be the set of unoriented
$3$-manifolds with boundary $S$ and let $\M_S^u$ be the $\C$ vector
space spanned by $\Mdot_S^u$.
The natural pairing $\M_S^u \times M_S^u \rightarrow \M^u$ is positive.
\end{lctheorem}

\begin{proof}
The key is again to define a suitable diagonally dominant complexity
function $c^u$. If $N$ is a closed $3$-manifold, it is almost enough
to take $c(\til{N})$ for $c^u(N)$ where $\til{N}$ denotes the
orientation cover of $N$ (which consists of two copies of $N$ when $N$
is already oriented) and $c$ is the complexity function which is
diagonally dominant on oriented $3$-manifolds.

This definition only requires two modifications. First recall
the terms of the form $c_\iota$ defined in \S\ref{c3_subsection}.
These terms count distinct elements $P$ in the prime decomposition
of $\til{N}$. The problem is that the same prime $P$ in $\til{N}$
may arise as the oriented double cover of several distinct non-orientable
manifolds $Q_i$. An orientation-reversing involution
$\theta_i:P \to P$ may be free, or have isolated fixed points giving
rise to $\RP^3$ summands in $Q_i:=P/\theta_i$.  This
problem is solved by keeping track of {\em pairs} $(P,\theta_i)$ in
the definition of $c_\iota$.

Second, for a Seifert fibered JSJ
piece $X \subset \til{N}$ we make an additional modification to the
associated tensor $T_X$, which is involved in the construction of
the $c_t$ complexity term, defined in \S\ref{acasect}. Let $b_\theta(X)$
denote the number of boundary components of $X/\theta$. Then multiply
$T_X$ by the factor $(1+b_\theta)$.

It remains to check diagonal dominance of $c^u$. Suppose $(A,S)$ and
$(B,S)$ are elements of $\Mdot_S^u$. The orientation covers $\til{AB}$, $\til{AA}$
and $\til{BB}$ are each assembled from $\til{A}$ and $\til{B}$ glued
along $\til{S}$. The only
way for diagonal dominance to fail is if $(\til{A},\til{S}) = (\til{B},\til{S})$
but $(A,S) \ne (B,S)$. One way this might happen is for different
free prime factors in $A$ and $B$ (or factors which intersect $S$ only in $S^2$ or
$\RP^2$ components) to be covered by the same prime factor in $\til{A}$ and
$\til{B}$; the modification of $c_\iota$ takes care of this case.

So to complete the proof, assume that $(A, S_0)$ is nonorientable,
$\RP^2$-irreducible, and boundary irreducible, and that $\partial A =
S_0$ is a (possibly noncompact) surface of finite type other than a disk
or M\"obius band. It
suffices to show that there is no distinct $(A', S_0)$ with
identical orientation covers $\til{A}' \cong \til{A}$ (rel $S_0$).
We look at the orientation involution $\theta: \til{A} \rightarrow
\til{A}$. The involution $\theta$ preserves the relative JSJ decomposition of
$\til{A}$. Note that the restriction of $\theta$ to $\til{S}_0$ does not
depend on $A$ or $A'$, but only on $S$. The restriction of this involution
propagates in a locally rigid way across the relative
JSJ pieces which make up $A$ and $A'$.  This is clear for a
hyperbolic JSJ piece $X$, since an isometry is
determined by its germ near some invariant subset.
For a Seifert fibered JSJ piece $X$ there may be several
conjugacy classes (with respect to the mapping class group action)
of extension of an involution at the cusps over $X$.
But our modified tensors $T_X$ give rise to $c_t$ terms which favor
involutions on both $\til{A}$ and $\til{B}$ which
self-pair boundary components of $X$ rather than pair distinct
components. So strict diagonal dominance holds unless $b_\theta(X) = b_\theta(X')$
for all corresponding pairs $X,X'$ of Seifert fibered JSJ pieces in $\til{A}$
and $\til{A}'$. But the data of $b_\theta(X)$ is enough information to propagate
$\theta$ uniquely over $X$, and inductively one sees that the covering
involutions act on the identical pieces $\til{A}$ and $\til{A}'$ in the same way
(rel $\til{S}_0$) and therefore $A$ and $A'$ are the same rel $S_0$. This concludes our
sketch of the proof.
\end{proof}

\subsection{Orbifolds}

For the basic definitions and properties of orbifolds (particularly in low dimensions),
see Thurston \cite{Thurston_notes}, Chapter~13. We remind the reader of the
important distinction between good and bad orbifolds:
\begin{definition}
An orbifold $M$ is {\em good} if it is (orbifold) covered by a
smooth manifold.
\end{definition}

The following important theorem of Thurston shows that bad orbifolds can be
``localized'' to a $2$-dimensional stratum.
\begin{theorem}[Thurston \cite{Thurston_notes}]
A $3$-orbifold, all of whose $2$-dimensional suborbifolds are good, is itself
good.
\end{theorem}
In particular, the result of gluing two good $3$-orbifolds (with boundary)
along a good $2$-orbifold is itself good.

In what follows, let $S$ denote a closed oriented good
$2$-orbifold, and $\Mdot_S^{or}$ the set of
isomorphism classes (rel $S$) of compact oriented good
$3$-orbifolds with boundary marked by an identification with $S$.
As before, let $\M_S^{or}$ denote the complex vector space spanned by the
set $\Mdot_S^{or}$.

\begin{lctheorem} \label{M_S^or_theorem}
The natural pairing $$\M_S^{or} \times \M_S^{or} \rightarrow \M^{or}$$ is positive.
\end{lctheorem}

\begin{proof}
We sketch the modifications to $c$ required to produce a
diagonally dominant complexity function $c^{or}$ on $\M^{or}$.

In the definition of $c_1$, replace fundamental groups by
orbifold fundamental groups. The Dijkgraaf-Witten TQFTs are perfectly
well defined on orbifolds, and the corresponding partition functions
$Z_F$ as a family, are maximize only when the kernels of the maps
$\pi_1^{or}(S) \to \pi_1^{or}(A)$ and $\pi_1^{or}(S) \to \pi_1^{or}(B)$
induced by inclusion agree. By the equivariant loop theorem and
equivariant Dehn's Lemma  (Meeks-Yau \cite{MY}), one obtains matching
families of compressing elliptic orbifolds with boundary on $S$
in $A$ and $B$.

By the equivariant sphere theorem (Meeks-Simon-Yau \cite{MSY}), good
$3$-orbifolds have prime decompositions along essential spherical
$2$-orbifolds. We can construct a sum graph $G^{or}$ associated to
the decomposition exactly as in \S\ref{c_2_section} with similar
conventions as in the manifold case, where only edges marked by
actual $2$-spheres and vertices by actual punctured $3$-spheres are
treated as {\em thin}.  To define $c_2^{or}$, little needs to be
changed. Define $r^{or}(M)$ to be the number of irreducible orbifold
summands other than $S^1 \times S^2$, and $s^{or}(M)$ to be the
number or $S^1 \times S^2$ summands.  Close analogs of
Theorem~\ref{Laudenbach_theorem} and Lemma~\ref{c^2_Lemma_Schema}
continue to hold (thin edges may be slid over themselves and over
``fat edges'' whereas ``fat edges'' are rigid and cannot be slid).
As in the original proof of Lemma~\ref{c^2_Lemma_Schema}, {\em
thick} glues to {\em thick} and {\em thin} to {\em thin}; similar
algebra lets us establish a suitable diagonal dominance for
$c_2^{or}$.

The next terms to consider are $c_S$ and $c_h$ (inside $c_p$).
Irreducible (in the orbifold sense) good orbifolds have canonical
JSJ decompositions along Euclidean $2$-orbifolds into $3$-orbifolds
which are hyperbolic or covered by Seifert-fibered manifolds
(Bonahon-Siebenmann \cite{Bonahon_Siebenmann}). These second class
of Seifert ``fibered'' orbifolds are somewhat annoying. A
fiber-preserving involution on a Seifert fibered manifold (even an
orientation-preserving one) may collapse some circle fiber into an
interval with mirror endpoints. Even worse are the Euclidean
$3$-orbifolds, which might admit no natural fiber structure at all
(even a singular one). In fact, the ``miracle'' of Euclidean
orbifolds is that the manifolds amongst them actually {\em do} admit
Seifert fibered structures; this is explained very carefully in
Peter Scott's well-known article on the geometries of $3$-manifolds
\cite{Scott}. Mirror interval fibers do not cause much trouble, but
the Euclidean orbifolds are a headache (in fact for our purposes we
only need to consider orbifolds which are quotients of sufficiently
large manifolds; hence by Theorem~\ref{mapping_class_exceptions} the
Euclidean orbifolds are the only ones which might not admit a fiber
structure).

One can deal with this issue in one of two ways. The first is to
restrict the definition of $\M_S^{or}$ to allow only fibered or
hyperbolic sub-orbifolds. The second is to deal with the unfibered
Euclidean orbifolds on their own, as a separate class, distinct from
the Seifert fibered and hyperbolic JSJ pieces. The good news is that
unfiberable Euclidean $3$-orbifolds --- even those with boundary ---
have a rigid geometric structure, up to similarity. If we modify
$c_S$ so that the terms $m,m'$ count the maximum number of
independent Euclidean $2$-orbifolds, and the number of JSJ Euclidean
$2$-orbifolds respectively, then maximizing $(m,-m')$ ensures that
fibered Euclidean $3$-orbifolds with boundary on $S$ are glued to
compatibly fibered Euclidean $3$-orbifolds on the other side, and
unfibered Euclidean $3$-orbifolds on $S$ are glued to unfibered
Euclidean $3$-orbifolds compatibly with their rigid geometric
structures. A term, analogous to the length spectrum term in the
definition of $c_h$, ensures that unfibered $3$-orbifolds are
identified in matching pairs by doubling along their intersection
with $S$.

The $c_h$ term requires no modification. The analogue of
Theorem~\ref{relative_AST} holds in the orbifold context with
essentially the same proof. Since all orbifolds in question are
good, one can pass to manifold covers and use equivariance of Ricci
flow with surgery in the cover to deduce the conclusion. Similarly,
the length spectrum is defined on conjugacy classes in the orbifold
fundamental group, and the analysis in Lemma~\ref{spectrum_lemma}
generalizes to the orbifold context.

Very little in \S\ref{c3_subsection} requires modification. The
unfibered Euclidean $3$-orbifolds need special treatment in
\S\ref{c_a_section}, but the trick is to use their rigid geometric
structures, and treat them in the same way that hyperbolic JSJ
pieces are treated. Namely, given a Euclidean $3$-orbifold $X$, let
$W$ be a fundamental domain for the action of the group of
isometries of $X$, and build a tensor $T_X$ made up from copies of
tensors $T_W$ which has exactly the same symmetries as $X$.
Unfibered JSJ $2$-orbifolds have finite symmetry groups which can be
dealt with by hand. The invariant $\gn$ makes sense, and is natural
and chiral, for $3$-orbifolds bounded by fibered Euclidean
$2$-orbifolds, and a suitable generalization of
Lemma~\ref{gn_homeo_extension_lemma} continues to hold.
\end{proof}

\subsection{Tangles}

Finally, let $S$ denote a closed, oriented surface with finitely
many marked points $\{p\}$.  Let $\Mdot^l_S$ denote the set of
(relative isomorphism classes of) oriented compact manifolds $A$
containing a tangle $l$ with $\partial A = S$ and $\partial l =
\{p\}$, and let $\M^l_S$ denote the complex vector space it formally
spans.

\begin{lctheorem} \label{tangle_theorem}
The natural pairing $$\M^l_S \times \M^l_S \rightarrow \M^l$$ is
positive.
\end{lctheorem}

\begin{proof}
Define a complexity function $c^l$ on pairs $(M,l)$ where $M$ is a
closed oriented $3$-manifold and $l$ is a link in $M$, as follows.
First define a complexity $c^l_c$ when $M$ is connected, and then
define
$$c^l(M,L) = (c_0(M), B(M,l), \{B(M_i, l_i)\},
\{c_c^l(M_i, l_i)\})$$ where the last term in brackets is a
lexicographical list of the complexities of the connected
components.  Here $B$ is a ``bad orbifold'' weight: $B(M,l)$ denotes
the number of connected summands of $(M,l)$ of the form $(S^1 \times
S^2, S^1 \times \ast)$, where $\ast \subset S^2$ is a base point and
the connected sum is taken {\em away} from the link $S^1 \times
\ast$, and $\{B(M_i,l_i)\}$ is the component-wise lexicographical
list counting summands of this type.

To complete the definition, we must define $c_c^l$.  The idea is
that once all $(S^1 \times S^2, S^1 \times \ast)$ summands have been
removed, the remainder $(M_-,l_-)$ will become a {\em good} orbifold
once any {\em fixed} positive integer $n \ge 2$ is attached to {\em
each} component of $l_-$.  Thus, it suffices to define $c_c^l(M,l) =
c^{or}(M_-,l_{-,2})$ where $(M_-, l_{-,2})$ denotes the good
orbifold obtained by marking all components of $l_-$ as cone
geodesics with angle $\pi$.

The number of $(S^1 \times S^2, S^1 \times \ast)$ in $(M,l)$ and
$\{(M_i, l_i)\}$ is maximized in the double.  By
Theorem~\ref{M_S^or_theorem}, the complexity $c^{or}$ is diagonally
dominant, so
$$c^{or}((A_-,l_{A_-})(B_-,l_{B_-})) <
\max(c^{or}((A_-,l_{A_-})(A_-,l_{A_-})),
c^{or}((B_-,l_{B_-})(B_-,l_{B_-})))$$ unless
$(A_-,l_{A_-})(B_-,l_{B_-}) \cong (A_-,l_{A_-})(A_-,l_{A_-}) \cong
(B_-,l_{B_-})(B_-,l_{B_-})$ as orbifolds. Since an isomorphism of
orbifolds necessarily induces an isomorphism of the orbifold loci,
we obtain an isomorphism in the category of pairs ($3$-manifold,
link). This completes the sketch.
\end{proof}

\section{Acknowledgments}

We would like to thank Ian Agol, Daryl Cooper, Nathan Dunfield,
Cameron Gordon, Alexei Kitaev, Sadayoshi Kojima, Marc Lackenby,
Darren Long, Shigenori Matsumoto, John Morgan, Marty
Scharlemann, and the anonymous referee for numerous helpful comments, suggestions and corrections.
We would especially like to single out Alexei Kitaev
for thanks for posing questions which were the initial stimulus for
this work.

Danny Calegari was partially funded by NSF grants DMS 0405491 and
DMS 0707130.

\appendix

\section{Finite group TQFT's}\label{TQFT_appendix}

\def\bd{\partial}

This is a very brief summary of finite group TQFTs. For more detail,
see \cite{Dijkgraaf_Witten}, \cite{Freed_Quinn}, \cite{Quinn},
\cite{Freed} and \cite{Walker}. Here we consider only the
$2{+}1$-dimensional untwisted theories, and we ignore gluing with
corners, higher codimension gluing, etc.

\medskip

Let $G$ be a finite group and $BG$ be a classifying space for $G$.
All constructions in this appendix depend of $G$, but that is left
notationally implicit.

For $X$ a manifold (of any dimension), let $\M(X)$ denote the space
of all continuous maps $X \to BG$. There is a restriction map
(fibration)
\[
    r: \M(X) \to \M(\bd X) .
\]
For $c \in \M(\bd X)$, let $\M(X; c)$ denote $r^{-1}(c)$.

Let $Y$ be a closed 2-manifold. Define the (finite dimensional)
Hilbert space $V(Y)$ to be the space of all locally constant
functions $f: \M(Y) \to \C$. In other words, $f$ assigns a number to
each homotopy class of map $Y \to BG$ (equivalently, to each
isomorphism class of principal $G$-bundle over $Y$). If $Y$ is the
empty 2-manifold (e.g.\ the boundary of a closed 3-manifold), then
we have a canonical identification $V(Y) = \C$.

(For twisted theories, one replaces locally constant maps with maps
which vary according to $f - g = \omega(h)$, where $\omega$ is a
fixed 3-cocycle on $BG$ and $h:Y\times I \to BG$ is a homotopy from
$g$ to $f$. To make sense of this, manifolds must be equipped with
fundamental classes; see \cite{Freed_Quinn}.)

Let $\alpha \in \pi_0(\M(Y))$ and let $\chi_\alpha$ denote the
function which is 1 on the path component $\alpha$ and 0 on the
other path components of $\M(Y)$. The various $\chi_\alpha$ form a
basis of $V(Y)$. The inner product on $V(Y)$ is given by
\[
    \langle \chi_\alpha , \chi_\beta \rangle = \left\{ \begin{array}{ll}
            0 & \alpha \ne \beta \\
            {\displaystyle \frac{1}{|\pi_1(\M(Y); \alpha)|}} & \mbox{otherwise} .
        \end{array} \right.
\]
In particular, $\{\chi_\alpha\}$ is an orthogonal basis of $V(Y)$.
Note that if $Y$ is connected and $\rho: \pi_1(Y) \to G$ is a
representative of the conjugacy class of group homomorphisms
corresponding to the path component $\alpha$, then $|\pi_1(\M(Y);
\alpha)| = |\stab(\rho)|$, where $\stab(\rho)$ denotes the
stabilizer of $\rho$ under the outer action of $G$ by conjugation
--- $\stab(\rho)$ is the set of all elements of $G$ which commute
with the image of $\rho$. In terms of bundles, $\pi_1(\M(Y);
\alpha)$ and $\stab(\rho)$ can be identified with the automorphisms
of the $G$-bundle over $Y$ corresponding to $\alpha$.

Next we define the path integral $Z(M) \in V(\bd M)$ for a
3-manifold $M$. For $a \in \M(\bd M)$, define
\[
    Z(M)(a) = \sum_{\beta \in \pi_0(\M(M; a))} \frac{1}{|\pi_1(\M(M; a); \beta)|} .
\]
In other words, we sum over the path components of the extensions of
$a$ to all of $M$, and each path component counts as the inverse of
the size of its fundamental group. (In terms of bundles,
$|\pi_1(\M(M; a); \beta)|$ is the number of automorphisms of the
bundle corresponding to $\beta$ which restrict to the identity over
$\bd M$.)

A basic result for finite group TQFTs is the following gluing formula.
Let $A$ and $B$ be 3-manifolds with common boundary $S$.
Then we have $Z(A), Z(B) \in V(S)$ and
\[
    Z(A \cup_S B) = \langle Z(A), Z(B) \rangle .
\]
The proof considers the homotopy long exact sequence of the
fibration $\M(A \cup_S B) \to \M(S)$, whose fiber at $x \in \M(S)$
is $\M(A; x) \times \M(B; x)$. All homotopy groups in the sequence
are finite, and only the six $\pi_1$ and $\pi_0$ terms are
non-trivial.  See \cite{Quinn} for details.

\section{Application to Heegaard Genus}\label{heegaard_app}

This appendix illustrates how the Dijkgraaf-Witten TQFTs can be
used to obtain nontrivial lower bounds on Heegaard genus.

For $i=1,2$ let $\gamma_i:S^1 \to S^1 \times S^2$ be disjoint
nonparallel embeddings so that the complement of the images is a
homology ($\text{torus} \times I$). Let $S_1,S_2$ be disjoint open
solid torus neighborhoods of the images
$\gamma_1(S^1),\gamma_2(S^1)$. Let $C = S^1 \times S^2 - S_1$ and
$C' = S^1 \times S^2 - (S_1 \cup S_2)$.

\begin{figure}[htpb]
\labellist
\small\hair 2pt
\pinlabel $C_1=C$ at 20 245
\pinlabel $DC=C\overline{C}$ at 220 245
\pinlabel $C_2=C'\overline{C}$ at 20 142
\pinlabel $DC_2=C\overline{C}'C'\overline{C}$ at 255 142
\pinlabel $C_3=C'\overline{C}'C'\overline{C}$ at 30 36
\pinlabel $DC_3=C\overline{C}'C'\overline{C}C'\overline{C}'C'\overline{C}$ at 300 36
\pinlabel $\cdots$ at 70 0
\endlabellist
\centering
\includegraphics[height=2in]{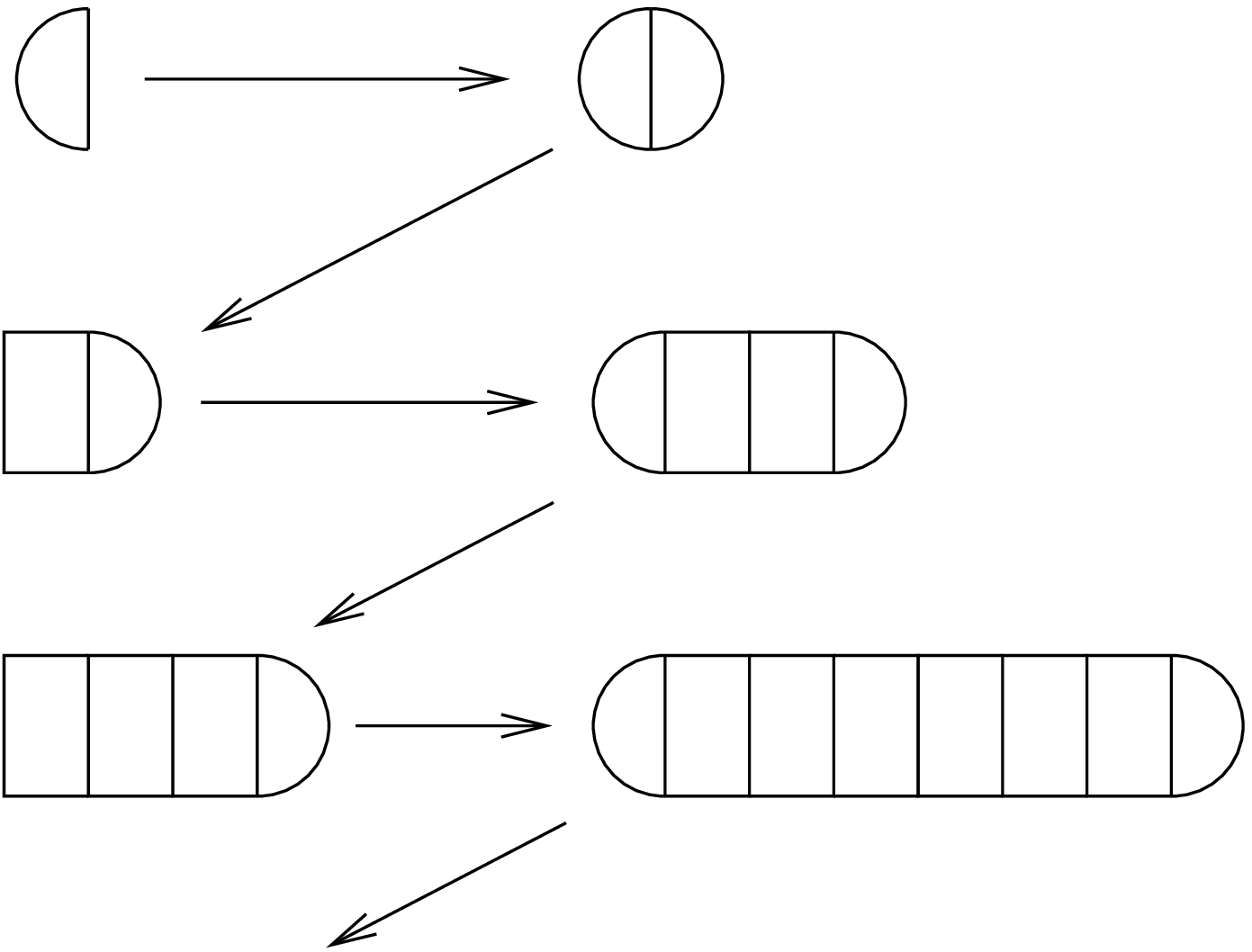}
\caption{Schematic picture of manifolds constructed in proof of Theorem~\ref{heegaard_thm}} \label{tqftappfig}
\end{figure}

For any integer $n$, let $DC_n$ be obtained from two copies of $C$
(with opposite orientations) and $2^n - 2$ copies of $C'$ (coming in
orientation-reversed pairs) glued together in the pattern
$$DC_n = C \overline{C}' C' \overline{C}' \cdots C' \overline{C}$$
where $C$ is glued to $\overline{C}'$ along their common boundary
component $\partial S_1$, where a $\overline{C}'C'$ term glues the
first to the second along a $\partial S_2$ component, and a
$C'\overline{C}'$ term glues the first to the second along a
$\partial S_1$ component. Notice that $DC_n$ can be built from
$2^n\cdot\text{const.}$ simplices, by choosing fixed triangulations
of $C$ and $C'$ which agree on their common boundary component; see
Figure~\ref{tqftappfig}.

Apart from the two copies of $C$ at either end, $DC_n$ is almost
equal to the cyclic cover of order $2^{n-1}$ of the double $DC'$,
where the cover unwraps the circle obtained by doubling an arc from
one boundary component of $C'$ to the other.

This appendix shows that the Heegaard genus of $DC_n$ is of size at
least $2^n\cdot\text{const.}$ where the (positive) constant depends
only on $\gamma_1,\gamma_2$. By the assumption on the $\gamma_i$,
the manifolds $DC_n$ are all homology $S^1 \times S^2$'s, so one
cannot obtain a lower bound on their Heegaard genus from homological
invariants. The point is not that this result cannot be obtained by
other means (although for some cases of the above construction, the
method of \cite{Schultens-2001} also yields a lower bound which is
optimal up to a constant); the point is that this estimate is
obtained as a formal consequence of the properties of
Dijkgraaf-Witten TQFTs. Classical $3$-manifold topology enters only implicitly
(but importantly) in the fact that $\pi_1(C)$ is residually finite.

For a closed $3$-manifold $M$ triangulated with $T$ tetrahedra,
there is a natural Heegaard splitting, where one side is taken to
be a regular neighborhood of the $1$-skeleton. The genus of the
splitting surface is $T+1$, which we take as a naive upper bound
for the Heegaard genus $h(M)$ of $M$.  We construct a family of
irreducible examples of fixed homological type where this naive
bound is sharp up to a constant.

\begin{theorem} \label{heegaard_thm}
There exists a constant $c > 0$ and a sequence of irreducible
3-manifolds $M_1$, $M_2, \cdots$ each of which is
$\mathbb{Z}$-homology equivalent to $S^1 \times S^2$, such that the
Heegaard genus satisfies $h(M_i) > c \cdot T(M_i)$, where $T$ denotes the
minimal number of tetrahedra in a triangulation of $M_i$, and
$T(M_i)$ approaches infinity.
\end{theorem}

\begin{proof}
Let $\gamma: S^1 \times D^2 \hookrightarrow S^1 \times S^2$ be an
embedding with incompressible complement $C$, so that $S^1 \times
S^2 = S \cup C$, where $S = S^1 \times D^2$.  Using residual
finiteness of $\pi_1(C)$ [H], let $F$ be a finite group so that
the conjugacy class of the meridian of $S$ in $\pi_1(C)$ maps nontrivially under some homomorphism
$\theta: \pi_1(C) \rightarrow F$.  Let $Z$ be the partition function
(or relative partition function) of the TQFT associated to $F$, which is unitary by
Appendix~\ref{TQFT_appendix}. Since the meridian survives in at
least one $F$-representation of $\pi_1(C)$, we have  $Z(S) \neq Z(C)$.
From the definition of $Z$, we have
\begin{equation}
Z(S^1 \times S^2) = \tr(\id_{V(S^2)}) = \tr(\id_\mathbb{C}) = 1
\end{equation}

Now, $\innerprod{Z(S)}{Z(S)} = Z(S^1 \times S^2) = 1$, so $|Z(S)| =
1$.  But $\innerprod{Z(C)}{Z(S)} = Z(S^1 \times S^2) = 1$ as well.
Since $Z(C) \neq Z(S)$, we conclude that $|Z(C)| = l > 1$ for some $l$.
Furthermore, if $DC$ denotes the double, then $|Z(DC)| = \innerprod{Z(C)}{Z(C)} = l^2$.

Now let $\omega : S^1 \times D^2 \hookrightarrow C$ be a fixed
embedding inducing an isomorphism on integer homology.
Define $C_2 = DC - \image(\omega)$, where we take $\image(\omega)$ to lie
in the positive copy of $C$ in $DC$. Then
$\innerprod{Z(S_\omega)}{Z(C_2)} = Z(DC) = l^2$ where $S_\omega$
denotes the solid torus embedded via $\omega$.  Thus, $|Z(C_2)| \geq
l^2$ and so $Z(DC_2) > l^4$.

Using another copy of $C \subset DC_2$, compose $\omega_2 : S^1
\times D^2 \hookrightarrow C \hookrightarrow DC_2$ and set $C_3 =
DC_2 - \image(\omega_2)$.  We conclude
$\innerprod{Z(S_{\omega_2})}{Z(C_3)} \geq l^4$ and so $|Z(DC_3)|
\geq l^8$.

Find $\omega_3 : S^1 \times D^2 \hookrightarrow DC_3$ (again, by
composing $\omega$ with an inclusion $C \hookrightarrow DC_3$.)
Define $C_4$ and similarly conclude that $|Z(C_4)| > l^{16}$.  By
induction, we obtain a sequence $\{DC_i\}$ of $\mathbb{Z}$-homology
$S^1 \times S^2$'s with
\begin{equation}
\label{eqn_b_1} Z(DC_i) > l^{2^i}, \;\;\;\; l > 1
\end{equation}

Associated to a finite group TQFT is a positive integer called the
{\em total quantum dimension}, which we denote $\mathcal{D}$.
This number is just the cardinality of the finite group $G$ on which
the TQFT is based. The only property we use is that $\mathcal{D} > 1$
for nontrivial groups, and the following Lemma:
\begin{lemma}
\label{appendix2a} $|Z(\#_g S^1 \times S^2)| = \mathcal{D}^{g-1}$
\end{lemma}
See e.g. \cite{Dijkgraaf_Witten}.

\begin{corollary}[to lemma]
\label{appendixcorr} $h(DC_i) > const \cdot 2^i$ for some positive
constant.
\end{corollary}

\begin{proof}
Let $h_g$ be a handlebody of genus $g$. By the Cauchy-Schwarz inequality,
$$Z(\#_g S^1 \times S^2) = \innerprod{Z(h_g)}{Z(h_g)} \geq Z(M)$$
when $h(M) \leq g$. Together with equation~\ref{eqn_b_1} this implies:
\begin{equation}
\label{eqn_b_2} l^{2^i} < \mathcal{D}^{h(DC_i) - 1}
\end{equation}
So, $\frac{log(l)}{log(\mathcal{D})}2^i < h(DC_i) - 1$, and thus,
$h(DC_i) > c'\cdot 2^i$ for come constant $c' > 0$.
\end{proof}

By choosing fixed triangulations of $C$ and $C \backslash \omega(S^1
\times D^2)$ with common boundaries, we see that there is a constant
$c''$ so that
\begin{equation}
\label{eqn_b_3} T(DC_i) < c'' \cdot 2^i
\end{equation}

Combining (\ref{appendixcorr}) and (\ref{eqn_b_3}), we have
\begin{equation}
\label{eqn_b_4} h(DC_i) > \frac{c'}{c''} T(DC_i)
\end{equation}
\end{proof}

Combinatorial methods for obtaining lower bounds on Heegaard genus
include extensions due to Schultens \cite{Schultens-2001} 
of the connect sum formula $h(M \# N) =
h(M) + h(N)$, and Lackenby's
sweepout method \cite{Lackenby-2002}, partially summarized
below.

Suppose a closed hyperbolic manifold $M$ has a family of finite
covers $M_i$ of degree $d_i$ with a uniform {\em spectral gap}; i.e.
$$\lambda_1(M_i) \ge \epsilon > 0$$
for some $\epsilon$, where $\lambda_1$ denotes the first positive
eigenvalue of the Laplacian on functions (see \cite{Lubotzky}
especially Chapter~4 for definitions and an introduction). Obtain an
index one minimal Heegaard surface $\Sigma_i \subset M_i$ by using a
minimax ``sweepout'' of $M_i$. A minimal surface in a hyperbolic
manifold has curvature bounded above by $-1$, so by Gauss-Bonnet,
$-\chi(\Sigma_i) \geq \frac{\area(\Sigma_i)}{2\pi}$. Moreover, since
$\Sigma_i$ is the maximal area surface in a minimax sweepout, there
is an inequality
$$\area(\Sigma_i) \ge \area(\Sigma_i')$$
where $\Sigma_i'$ is one of the sweepout surfaces which evenly
divides the volume of $M_i$. Cheeger's isoperimetric constant $c$
(see \cite{Lubotzky}, Chapter~4), defined by
$$c = \inf_{\text{separating } \Sigma'} \frac{\text{area}(\Sigma')}{\min(\text{vol}( \text{component }
M-\Sigma'))}$$ satisfies $\lambda_1\cdot\text{const.}  < c <
\sqrt{\lambda_1}\cdot\text{const}$. Therefore we may conclude that
$-\chi(\Sigma_i) > d_i\cdot\text{const} $. In other words, the ratio
of Heegaard genus to volume (or to $T(M_i)$) is bounded below
uniformly for the family $\lbrace M_i \rbrace$; one says the
$\{M_i\}$ have a positive {\em Heegaard gradient}.

These three methods for bounding $h(M)$ from below appear quite
different. It would be interesting to see if they can be usefully
combined.  There is no shortage of interesting objectives for
sufficiently powerful lower bounds on $h$.  Lackenby has pointed out
that any hyperbolic $3$-manifold exhibiting positive Heegaard
gradient for all finite covers would be a counterexample to
Thurston's virtually fibered conjecture.

\section{Conjecture} \label{appendix_conjecture}
First consider the following conjecture in linear algebra.  We call
it the ``quantum maxflow/mincut''.  Let $N$ be a $v$-valent graph
with $k$ input cut edges and $l$ output cut edges.  At each vertex
of $N$ we presume that the incoming edges are locally ordered
$1,2,\dots,k$.  Given this data on $N$, picking a $v$-index tensor
$W$ (on some vector space $V$) transforms $N$ into a linear map
$\theta_{N,W}: V^{\otimes k} \rightarrow V^{\otimes l}$.  Simply
insert a copy of $W$ at each vertex according to the local edge
labeling and contract where appropriate.

\begin{conjecture}[QMF/MC] \label{qmf/mc_conjecture}
$\theta_{N,W}$ is an injection for generic $W$ if and only if $N$
admits $k$ edge-disjoint paths from input to output.  Furthermore,
$image(\theta_{N,W})$ is generic with respect to any fixed basis for
$V$.  That is, if $h \leq l-k$ and basis vectors $\vec{b}$ are
inserted into $h$ output slots not on the $k$ edge disjoint paths,
then the resulting map $\theta_{N,W,\vec{b}}: V^{\otimes k}
\rightarrow V^{\otimes (l-h)}$ is also an injection for generic $W$.
\end{conjecture}

\begin{remark}
\begin{enumerate}
\item The ``only if'' part follows form the ordinary maxflow/mincut
theorem.  If the $k$ edge-disjoint paths do not exist, then there is
a $j$-cut, $j < k$, separating input from output and factoring
$\theta_{N,W}$ through $V^{\otimes j}$
\item The nonsingularity of $\theta_{N,W}$ is an algebraic
condition, so, fixing $N$, it will hold generically for $W$ provided
it holds for a single $W$.
\item The condition that $W$ becomes an injective operator for every
division of the index (slots) into $v = a+b$, $a \leq b$, is a
finite intersection of generic conditions and therefore generic.  As
a consequence, the conjecture is true for {\em monotone} networks
$N$.  These are networks $N$ admitting a function to the interval
$\phi : N \rightarrow [0,1]$, $\phi(\text{input}) = 0$,
$\phi(\text{output})=1$, with $|\phi^{-1}(y)| \geq |\phi^{-1}(x)|$
for all $0 \leq x < y \leq 1$, provided $\phi(x) \neq \phi(w)$, $w$
a vertex of $N$.
\end{enumerate}
\end{remark}

In section \ref{c_a_section} we constructed a tensor $T_X$ with
positivity property $P$ for $X$ SF.  In our proof, the distinction
between algebraic and geometric middle thirds $\text{AMT} \subset
\text{GMT}$ was necessitated by a lack of a similar tensor for the
hyperbolic JSJ pieces.  (In an early draft of this paper, we thought
we had constructed such tensors but relied on a mistaken ``proof''
of QMF/MC.)  The existence of such tensors is, for hyperbolic JSJ
pieces, an independently interesting question.

Given a finite volume hyperbolic three manifold $X$, consider a
tensor $T_X$ living on $\displaystyle Q^m = \bigotimes_{\text{cusps
} C} V_m(C)$, where $V_m(C)$ is as in section \ref{c_a_section}. The
important properties of $V_m(C)$ are that (1) it is canonically
identified to the cusp so that actions on the cusps induce actions
on $Q$, (2) the action is ``eventually faithful'' as $m \rightarrow
\infty$ and also has many orbits ($\geq m$), and (3) $V_m(C)$ has a
natural basis which we regard as generating a positive cone.  We
assume $V_m(C)$ has these properties:

\begin{definition} \label{prop_ph_def}
$T_X$ has property $P_h$ if
\begin{enumerate}
\item every symmetry of $T_X$ is induced by a unique symmetry of $X$
\item for each reflection symmetry $r$ (of $X$), each operator
$\mathcal{O}$ from ``lower'' to ``upper'' cusp spaces induced by
inserting (positive) basis vectors into the cusps meeting $fix(r)$
is positive: $\langle x | \mathcal{O} | x \rangle > 0$ for
$|x\rangle \neq 0$ in the lower cusp space and $\langle x | = r_* |
x \rangle$ the corresponding element of the upper cusp space.
\end{enumerate}
\end{definition}

\begin{lemma} \label{QMF/MC_lemma}
QMF/MC implies that for each hyperbolic finite volume $3$-manifold
$X$, there exists a tensor $T_X$ satisfying property $P_h$.
\end{lemma}

\begin{proof}[Proof sketch]
Let $G$ be the symmetry group of $X$ and $W$ a fundamental domain.
Build a tensor $T_W$ as indicated in Figure \ref{twpic}, which has a
$V_m$ index (drawn as an edge) exiting each interior cusp of $W$ and
also a $V_m$ index exiting each fractional cusp located on a
singular stratum of $W$. At every codimension one face of $W$ the
tensor $T_W$ should have $l$ indices, each taking values in $V_m$.
Here, the number $l$ of ``legs'' is chosen sufficiently large; we
will say just how large in what follows. Pictorially, we think of
$T_W$ as a high valence vertex in $W$ with $l$ legs exiting every
face and one leg exiting each cusp and cusp fraction.

\begin{figure}[htpb]
\labellist
\small\hair 2pt
\pinlabel $V_m$ at 385 175
\pinlabel $V_m$ at 480 40
\pinlabel $V_m$ at 270 230
\pinlabel $V_m$ at 300 205
\pinlabel $l\;V_m\text{'s}$ at 410 80
\pinlabel $l\;V_m\text{'s}$ at 120 -10
\pinlabel $l\;V_m\text{'s}$ at -20 125
\endlabellist
\centering
\includegraphics[height=1.5in]{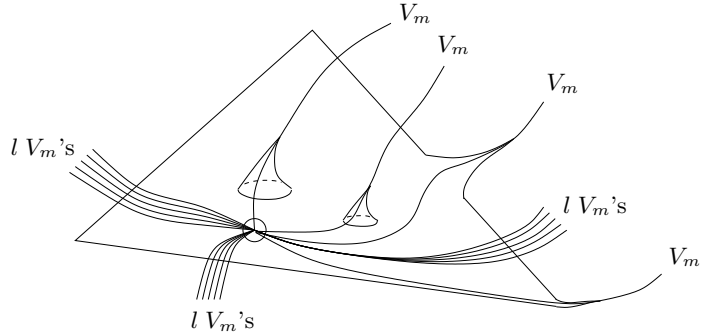}
\caption{This diagram represents an alternative construction of
$T_W$ given \ref{qmf/mc_conjecture}.  One leg exits every cusp and
fractional cusp, $l$ legs exit every face.} \label{twpic}
\end{figure}

Use $r$ to define upper and lower halves of $X$. Consider the
network of tensors coming from the copies of $W$ in the upper half
of $X$; see Figure~\ref{tensors}. We think of this network of
tensors as an operator $U$ from the upper (possibly fractional) cusp
edges to the ``middle'' edges along the codimension 1 fixed
submanifold of $r$. Choose $l$ large enough so that we can draw
mutually edge-disjoint paths from each upper cusp edge to the middle
level.  The input from $r$-invariant cusps corresponds to inserting
some positive basic vectors at the nodes closest to the mid level in
Figure \ref{tensors}.  The QMF/MC implies that a generic choice of
$T_W$ yields an injective operator $\mathcal{U}$ when inserted into
the network.  Thus, $\mathcal{O}=\mathcal{U}^\dagger \mathcal{U}$ is
strictly positive.

\begin{figure}[htpb]
\labellist
\small\hair 2pt
  \pinlabel $\text{The tensors }T_W=$ at 30 163
  \pinlabel $\text{and }T_W^\dagger$ at 165 163
  \pinlabel $\text{contract to yield}$ at 24 143
  \pinlabel $\text{an injective morphism}$ at 42 123
  \pinlabel $\text{from ``top'' to ``mid level''}$ at 54 103
  \pinlabel $\text{cusp and}$ at 480 420
  \pinlabel $\text{cusp fractions}$ at 495 400
  \pinlabel $\text{at top}$ at 470 380
  \pinlabel $\text{edge-disjoint}$ at 492 280
  \pinlabel $\text{paths in bold}$ at 495 260
  \pinlabel $\text{mid level}$ at 490 82
\endlabellist
\centering
\includegraphics[scale=0.7]{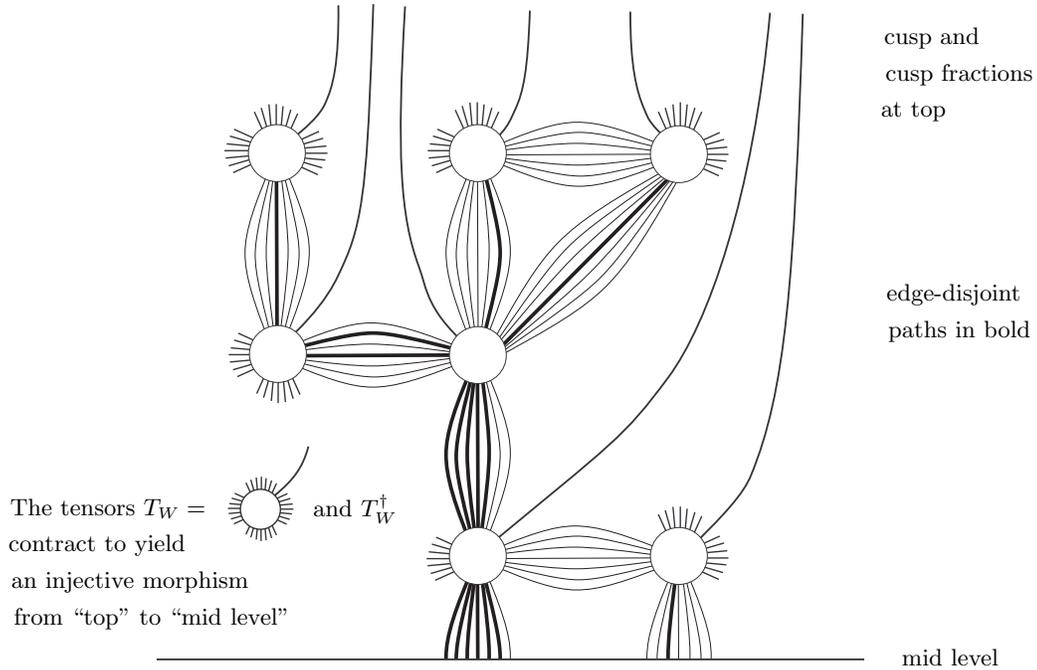}
\caption{In this figure $l=7$} \label{tensors}
\end{figure}

\end{proof}

\bibliographystyle{plain}
\bibliography{references}

\end{document}